\documentclass[12pt]{article}

\input amssym.def
\input amssym.tex

\usepackage{euscript}

\begin{document}

\renewcommand{\theequation}{\thesection.\arabic{equation}}
\newcommand{\renorm}{{ \setcounter{equation}{0} }}

\title{Complete sets of relations in the cohomology rings of moduli spaces of
holomorphic bundles and parabolic bundles over a Riemann surface}
\author{Richard Earl and Frances Kirwan}



\newtheorem{prop}{Proposition}[section]
\newtheorem{lem}[prop]{Lemma}
\newtheorem{cor}[prop]{Corollary}
\newtheorem{thm}[prop]{Theorem}
\newtheorem{guess}{Conjecture}
\newtheorem{dff}[prop]{Definition}
\newenvironment{df}{\begin{dff}\rm}{\end{dff}}
\newtheorem{REM}[prop]{Remark}
\newenvironment{rem}{\begin{REM}\rm}{\end{REM}}
\newtheorem{examplit}[prop]{Example}
\newenvironment{example}{\begin{examplit}\rm}{\end{examplit}}

\newcommand{\aff}{{\Bbb A }}
\newcommand{\RR}{{\Bbb R }}
\newcommand{\CC}{{\Bbb C }}
\newcommand{\ZZ}{{\Bbb Z }}
\newcommand{\PP}{ {\Bbb P } }
\newcommand{\QQ}{{\Bbb Q }}
\newcommand{\UU}{{\Bbb U }}

\newcommand{\Q}{\QQ}
\newcommand{\nn}{\tilde{n}}
\newcommand{\DD}{\tilde{d}}
\newcommand{\mnd}{{\cal M}(n,d)}
\newcommand{\mndl}{{\cal M}_{\Lambda}(n,d)}
\newcommand{\Mnd}{{\cal M}(n,d)}
\newcommand{\dai}{\frac{\partial}{\partial a_{i}}}
\newcommand{\daij}{\frac{\partial^{2}}{\partial a_{i} \partial a_{j}}}
\newcommand{\HS}{H^{*}}
\newcommand{\ar}{a_{r}}
\newcommand{\bjr}{b_{r}^{j}}
\newcommand{\fr}{f_{r}}
\newcommand{\Z}{\ZZ}
\newcommand{\Res}{\mbox{Res}}
\newcommand{\ch}{\mbox{ch}}
\newcommand{\mtc}{M_{\T}(c)}
\newcommand{\n}{\hat{n}}
\newcommand{\D}{\hat{d}}
\newcommand{\E}{\hat{E}}
\newcommand{\V}{\hat{V}}
\newcommand{\M}{{\cal M}(\n,\D)}
\newcommand{\G}{{\cal G}}
\newcommand{\GG}{\overline{\cal G}}
\newcommand{\pa}{\overline{\partial}}
\newcommand{\g}{(g-1)}
\newcommand{\Cmu}{{\cal C}_{\mu}}
\newcommand{\HG}{H^{*}_{\cal G}}
\newcommand{\Css}{{\cal C}^{ss}}
\newcommand{\C}{{\cal C}}
\newcommand{\Cnu}{{\cal C}_{\nu}}
\maketitle

\indent The cohomology of the moduli space $\mnd$ of stable holomorphic
vector bundles of coprime rank $n$ and degree $d$ over a fixed compact
Riemann surface $\Sigma$ of genus $g \geq 2$ has been intensely
studied over many years. In the case when $n=2$ we now
 have a very thorough understanding of its
structure \cite{B,KN,ST,Z}. For arbitrary $n$ it is known that the cohomology has no
torsion \cite{AB}, and inductive and closed
formulas for computing the Betti
numbers have been obtained \cite{AB,DB,DR,HN,LR}, as well as a set of
generators for the cohomology ring  \cite{AB}. 

When $n=2$ the
relations between these generators can be explicitly described
\cite{B,HS,KN,ST,T,Z} and in particular a conjecture of Mumford, that a
certain set of relations is a complete set, is now known to be true
\cite{K,Z}. 
Less is known about the relations between the
generators when $n>2$, although in principle these can be obtained
from the formulas given in \cite{JK2,WI} for the evaluation of polynomials
in the generators on the fundamental class of $\mnd$. It was shown in \cite{E} that the most
obvious generalisation of Mumford's conjecture to the cases when $n>2$
is false, although a modified version of the conjecture (using `dual
Mumford relations' together with the original Mumford relations) is
true for $n=3$. In this paper we generalise the
concept of the Mumford relations somewhat further and show that these
generalised Mumford relations form a complete set for arbitrary rank $n$. 
\\
\indent The generators for $\HS(\mnd)$ given by Atiyah and Bott in \cite{AB} are
obtained from a (normalised) universal bundle $V$ over $\mnd \times
\Sigma$. With respect to the K\"{u}nneth decomposition of 
\[
\HS(\mnd \times \Sigma)
\]
the $r$th Chern class $c_{r}(V)$ of $V$ can be
written as
\[
c_{r}(V) = a_{r} \otimes 1 + \sum_{j=1}^{2g} b_{r}^{j} \otimes
\alpha_{j} + f_{r} \otimes \omega
\]
where $\{1\}, \{\alpha_{j}:1 \leq j \leq 2g\},$ and $\{ \omega\}$ are
standard bases for $H^0(\Sigma),H^1(\Sigma)$ and $H^2(\Sigma)$, and
\begin{equation}
a_{r} \in H^{2r}(\mnd), \quad b_{r}^{j} \in H^{2r-1}(\mnd), \quad
f_{r} \in H^{2r-2}(\mnd), \label{gen}
\end{equation}
for $1 \leq r \leq n$ and $1\leq j \leq 2g$. It was shown in
\cite[Prop. 2.20 and 
p.580]{AB} that the classes $a_{r}$ and $f_{r}$ (for $2 \leq r \leq
n$) and $b_{j}^r$ (for $1 \leq r \leq n$ and $1 \leq j \leq 2g$)
generate the rational cohomology ring of $\mnd$.\\
\indent Since tensoring by a fixed holomorphic line bundle of degree
$e$ gives an isomorphism between the moduli spaces $\mnd$ and ${\cal
M}(n,d+ne)$, we may assume without loss of generality that
\[
(2g-2)n < d < (2g-1)n.
\]
This implies that $H^1(\Sigma ,E) = 0$ for any stable bundle of rank
$n$ and degree $d$ \cite[Lemma 5.2]{N2}, and hence that $\pi_{!}V$ is
a bundle of rank $d-n(g-1)$ over $\mnd$, where 
\[
\pi: \mnd \times \Sigma \to \mnd 
\]
is the projection onto the first component and $\pi_!$ is the
K-theoretic direct image map. It follows that 
\[
c_{r} (\pi_! V) = 0,
\]
for $r>d-n(g-1)$. Via the Grothendieck-Riemann-Roch theorem we can
express the Chern classes of $\pi_! V$ as polynomials in the
generators $a_r, b_r^j, f_r$ described above, and hence their
vanishing gives us relations between these generators. 
Mumford
conjectured \cite[p.582]{AB} that when $n=2$ these relations can be
used to obtain a complete set. This is indeed true when $n=2$ (see
\cite{K,Z}) but not in general for $n>2$ \cite[p.19]{E}. However we
will show that a similar construction can be used to obtain a complete set of relations,
as follows.\\
\indent Suppose that $0<\n<n$, and that $\D$ is coprime to $\n$. Then
we have a universal bundle $\V$ over $\M \times \Sigma$, and both $V$
and $\V$ can be pulled back to $\M \times \mnd \times \Sigma$. If
\[
\frac{\D}{\n} > \frac{d}{n}
\]
then there are no nonzero holomorphic bundle maps from a stable
bundle of rank $\n$ and degree $\D$ to a stable bundle of rank $n$ and
degree $d$, and hence, if 
\[
\pi: \M \times \mnd \times \Sigma \to \M \times \mnd,
\]
is the projection onto the first two components, it follows that 
$ -\pi_! (\V^* \otimes V) $
is a bundle of rank $n\n(g-1)-d\n+\D n$ over $\M \times \mnd$. Thus 
\[
0 = c_{r} (-\pi_! (\V^* \otimes V)) \in  \HS(\M \times \mnd)
\]
if $r>n\n(g-1)-d\n+\D n$ and hence the slant product
\[
c_{r} (-\pi_! (\V^* \otimes V)) \backslash \gamma \in \HS(\mnd)
\]
of $c_{r} (-\pi_! (\V^* \otimes V))$ with any homology class $\gamma
\in H_*(\M)$ vanishes when 
\[
r>n\n(g-1)-d\n+\D n.
\]
The main result of this paper (see Theorem \ref{thm} in $\S$2
below) is that the relations obtained
between the generators $a_r,b_r^j,f_r$ in this way\footnote{A little 
more care must be taken when $\n$ and $\D$ are not coprime.} 
for $0<\n<n$ and 
\[
\frac{d}{n}+1 > \frac{\D}{\n} > \frac{d}{n}
\]
and 
\[
n\n(g-1)-d\n+\D n <r < n\n(g+1)-d\n+\D n
\]
form a complete set of relations (which are essentially Mumford's
relations when $n=2$ and are essentially the relations described in
\cite{E} when $n=3$.) This is generalised to moduli spaces of
parabolic bundles in Theorem \ref{Par2.1}.

In principle all the relations between the Atiyah-Bott generators can be
obtained by Poincar\'{e} duality from the formulas given in \cite{JK2,WI}
for the evaluation of polynomials in these generators on the fundamental
class $[\Mnd ]$, but the linear algebra involved would be extremely
laborious, and our approach is a more
geometric one. However we use the formulas of \cite{JK2}  in the final section of the paper 
to give explicit formulas for
the generalised Mumford relations in terms of the Atiyah-Bott
generators.

 The first section of this paper describes background needed
from \cite{AB} and \cite{E} and the second gives a careful statement
of the main result, Theorem \ref{thm}. $\S$3 introduces modified
completeness criteria for relations between the generators similar to those used in \cite{E,K}. $\S$4
contains the heart of the proof, which involves a diagonal argument and the 
(excess) Porteous formula, and $\S$5 refines this to complete the case when 
$\D$ is coprime to $\n$. $\S$6 generalises what has been done so far to
cover parabolic bundles, and $\S$7 proves a weak version of the main result in the case of
parabolic bundles with full flags. This is extended to the general case in  
$\S$8 using the method of reduction to a maximal
torus (cf. \cite{Knew}). $\S$9 completes the proof of the main theorem by restricting the 
range of $\D$ and $r$ required for a complete set of relations, and finally
$\S$10 uses nonabelian localisation \cite{JK2} to
provide some explicit formulas for the generalised Mumford relations.

\section{Background}

\renorm

There are many descriptions, from very different points of view, of the
moduli space $\mnd$ of stable holomorphic bundles of coprime rank $n$
and degree $d$ over a compact Riemann surface $\Sigma$ of genus $g
\geq 2$: see for example \cite{AB,G2,J,M,NR2,NR3,N2}. The construction
we shall use is due to Atiyah and Bott in \cite{AB}; it is
closely analogous to the construction of quotients in geometric
invariant theory \cite{LP,MFK,N2} but involves a space and a group which
are both infinite dimensional. There are other constructions of $\mnd$
as genuine geometric invariant theoretic quotients of nonsingular
quasi-projective varieties by reductive group actions, which are in
some sense finite dimensional approximations to the construction of
Atiyah and Bott (see \cite{K3}), and the arguments we shall use could
all be rewritten in terms of these finite-dimensional quotient
constructions.\\
\indent In \cite{AB} Atiyah and Bott obtained generators for the
cohomology ring of $\mnd$. In the next section we shall state our
main theorem giving a complete set of relations between these
generators. In fact it will be enough to describe the relations
between generators over the rationals, since the cohomology ring of
$\mnd$ is torsion-free \cite[p.578]{AB}. All cohomology groups in this
paper will therefore have rational coefficients unless specified otherwise.\\
\indent Recall that a holomorphic vector bundle $E$ over $\Sigma$ is
called {\em  semistable} (respectively {\em stable}) if every holomorphic
subbundle $D$ of $E$ satisfies
\[
\mu (D) \leq \mu(E) , \indent (\mbox{respectively } \mu(D) < \mu(E)),
\]
where $\mu(D) =$ degree($D$)/rank($D$) is the {\em slope} of
$D$. Nonsemistable bundles are said to be unstable. Note that
semistable bundles of coprime rank and degree are stable.\\
\indent Let ${\cal E}$ be a fixed $C^{\infty}$ complex vector bundle
of rank $n$ and degree $d$ over $\Sigma$. Let $\C$ be the space of all holomorphic
structures on ${\cal E}$ and let $\G_{c}$ denote the group of all
$C^{\infty}$ complex automorphisms of ${\cal E}$. Atiyah and Bott
\cite{AB} identify the moduli space
$\mnd$ with the quotient $\Css/{\cal G}_{c}$ where $\Css$ is the open
subset of $\C$ consisting of all semistable holomorphic structures on
${\cal E}$. The group $\G_{c}$ is the
complexification of the gauge group $\G$ which consists of all smooth automorphisms
of ${\cal E}$ which are unitary with respect to a fixed Hermitian
structure on ${\cal E}$ \cite[p.570]{AB}. We shall write
$\overline{\G}$ for the quotient of $\G$ by its $U(1)$-centre and
$\overline{\G}_{c}$ for the quotient of $\G_{c}$ by its ${\bf 
C}^{*}$-centre. There are natural isomorphisms \cite[9.1]{AB}
\[
H^{*}(\Css/\G_{c}) = H^{*}(\Css/\overline{\G}_{c}) \cong
H^{*}_{\overline{\G}_{c}}(\Css) \cong H^{*}_{\overline{\G}}(\Css)
\]
since the centre of $\G_{c}$ acts trivially on $\Css$, while
$\overline{\G}_{c}$ acts freely on $\Css$ and $\overline{\G}_{c}$ is
the complexification of $\overline{\G}$. Atiyah and Bott \cite[Thm.
7.14]{AB} show that the restriction map $H^{*}_{\overline{\G}}(\C)
\rightarrow H^{*}_{\overline{\G}}(\Css)$ is surjective. Further
$H^{*}_{\overline{\G}}(\C) \cong
H^{*}(B\overline{\G})$ since $\C$ is an infinite-dimensional affine
space \cite[p.565]{AB}. So putting this all together we find that
there is a natural surjection
\begin{equation} \Phi:
H^{*}(B\overline{\G}) \cong H^{*}_{\overline{\G}}(\C)
\rightarrow H^{*}_{\overline{\G}}(\Css) \cong H^{*}(\mnd). \label{14}
\end{equation}
Thus generators of the cohomology ring
$\HS(B\overline{\G})$ give generators of the cohomology ring of the
moduli space $\mnd$. 

The cohomology rings of the classifying spaces
$B\G$ and $B\overline{\G}$ are described in \cite[$\S\S$ 2, 9]{AB}.
It is shown in \cite[Prop. 2.4]{AB} that the classifying space $B\G$
can be identified with the space $\mbox{Map}_{d}(\Sigma,BU(n))$ of all
smooth maps $f:\Sigma \rightarrow BU(n)$ such that the pullback to $\Sigma$ of
the universal vector bundle over $BU(n)$ has degree $d$. If we
pull back this universal bundle using the evaluation map
\[
\mbox{Map}_{d}(\Sigma,BU(n)) \times \Sigma \rightarrow BU(n): (f,m) \mapsto f(m)
\]
then we obtain a rank $n$ vector bundle ${\cal V}_0$ over $B\G \times
\Sigma$. If further we restrict the pullback bundle induced by the maps
\[
\Css \times E\G \times \Sigma \rightarrow \C \times E\G \times \Sigma \rightarrow \C
\times_{\G} E\G \times \Sigma \stackrel{\simeq}{\rightarrow} B\G \times \Sigma
\]
to $\Css \times \{e\} \times \Sigma$ for some $e \in E\G$ then we obtain
a $\G$-equivariant holomorphic bundle on $\Css \times \Sigma$. The
centre of
$\G$ acts as scalar multiplication on the fibres,
and the associated projective bundle descends to a holomorphic
projective bundle over $\mnd \times \Sigma$.\\
\indent In fact this projective bundle is the projective bundle of a
holomorphic vector bundle $V$ over $\mnd \times \Sigma$
\cite[pp. 579-580]{AB}. This vector bundle $V$ has the universal
property that, for any $[E] \in \mnd$ representing a bundle $E$ over
$\Sigma$, the restriction of $V$ to $\{[E]\} \times \Sigma$ is
isomorphic to $E$.\\
\indent Before constructing such a universal bundle $V$ we note that
tensoring with the pullback to $\mnd \times \Sigma$ of any holomorphic
bundle $K$ over $\mnd$ preserves the universal property and the
associated projective bundle remains unchanged. If we fix $x \in
\Sigma$ and identify $\Css$ with the subset $\Css \times \{x\}$ of
$\Css \times \Sigma$ then the ${\bf C}^*$-centre of $\G_c$ acts on
the fibres of $\left. \det({\cal V}_0 \right|_{\Css})$ with weight $n$
and acts on the fibres of $\det (\pi_!{\cal V}_0)$ with weight
$d-n(g-1)$ where $\pi:\Css \times \Sigma \to \Css$ is the first
projection and $\pi_!$ is the direct image map of K-theory. Since $n$
and $d$ are coprime there exist integers $u$ and $v$ such that 
\[
u n + v(d- n(g-1)) = 1.
\]
Thus the ${\bf C}^*$-centre of $\G_c$ acts on the fibres of 
\[
K =  \det(\left. {\cal V}_0 \right|_{\Css})^u \otimes \det (\pi_!{\cal
V}_0)^v
\]
as scalar multiplication, and hence acts trivially on the bundle
\[
{\cal V} = {\cal V}_0 \otimes \pi^*(K^{-1})
\]
over $\Css \times \Sigma$. This bundle ${\cal V}$ then descends to a
universal bundle $V$ over $\mnd \times \Sigma$.\\
\indent Tensoring a universal bundle $V$ by the pullback to $\mnd
\times \Sigma$ of any holomorphic bundle $K$ over $\mnd$ changes the
generators $a_r,b_r^j$ and $f_r$ of $\HS(\mnd)$. In particular it
changes $a_1$ by $nc_1(K)$ and $c_1(\pi_!(V))$ by $(d -
n(g-1))c_1(K)$. (Here we are using $\pi$ to
denote the first projection $\mnd \times \Sigma \to \mnd$ as well as
$\pi:\Css \times \Sigma \to \Css$.) Therefore tensoring $V$ by
$\pi^*(K)$ changes $ua_1 + v c_1(\pi_!V)$ by $c_1(K)$, with $u$ and $v$
chosen as above. Our choice of universal bundle $V$ satisfies
\begin{equation}
ua_1 + v c_1(\pi_!V)=0, \label{norm}
\end{equation}
i.e. it is normalised in the sense used by Atiyah and Bott
\cite[p.582]{AB}.\\
\indent By a slight abuse of notation we define elements $\ar, \bjr, \fr$
in $\HS(B\G;\Q)$ by writing
\[
c_{r}({\cal V}) = \ar \otimes 1 + \sum_{j=1}^{2g} \bjr \otimes
\alpha_{j} + \fr \otimes \omega \indent 1 \leq r \leq n.
\]
where, as before, $\omega$ is the standard generator of $H^{2}(\Sigma)$ and
$\alpha_{1},...,\alpha_{2g}$ form a fixed canonical cohomology basis for
$H^{1}(\Sigma)$.
Atiyah and Bott show \cite[Prop. 2.20]{AB} that the ring
$H^{*}(B\G)$ is freely generated 
as a graded algebra over $\Q$ by the elements 
\begin{equation} \label{genabf}
\{\ar : 1 \leq r \leq n\} \cup \{\bjr : 1 \leq r \leq n, 1 \leq j \leq 2g\} \cup 
\{ \fr : 2 \leq r \leq n\}.
\end{equation}
The only relations among these generators are that the 
$\ar$ and $\fr$ commute with everything else and that the
$\bjr$ anticommute with each other.\\
\indent The fibration $BU(1) \rightarrow B\G \rightarrow
B\overline{\cal G}$ induces an isomorphism \cite[p.577]{AB}
\[
H^{*}(B\G) \cong H^{*}(B\overline{\G}) \otimes H^{*}(BU(1)),
\]
where $H^*(BU(1))$ is a polynomial ring on one generator in degree 
two, which we can identify with $a_1$. The generators $\ar,\bjr$ and
$\fr$ of $\HS(B\G)$ can be pulled back via a section of this 
fibration to give rational generators of the cohomology ring of
$B\overline{\cal G}$. We may if we wish omit $a_{1}$ since its image
in $H^{*}(B\overline{\G})$ can be expressed in terms of the other
generators. The only other relations are again the commuting of the $\ar$ and
$\fr$ with all other generators, and the anticommuting of the $\bjr$. These generators for $H^{*}(B\overline{\G})$
restrict to the generators $\ar, \bjr, \fr$ (\ref{gen}) for
$\HS(\mnd)$ under the surjection (\ref{14}).\\
\indent The relations between these generators for $\HS(\mnd;\Q)$ are
then given by the kernel of the restriction map (\ref{14}) which in
turn is determined by the restriction map
\begin{equation}
\HG(\C) \cong  H^{*}_{\overline{\G}}(\C) \otimes \HS(BU(1)) \rightarrow
H^{*}_{\overline{\G}}(\Css) \otimes
\HS(BU(1)) \cong \HG(\Css). \label{surj}
\end{equation}
\indent  The bundle $V$ is universal in
the sense that its
restriction to $\{[E]\} \times \Sigma$ is isomorphic to $E$ for each
semistable holomorphic bundle $E$ over $\Sigma$ of rank $n$ and degree $d$,
and $\pi_! V$ is the alternating sum of the
$i$th higher direct image sheaves $R^{i}\pi_{*}V$ (see \cite[$\S
3.8$]{H}), whose stalks at $[E]$ are given by
\[
H^{i}(\pi^{-1}([E]),V_{|\pi^{-1}([E])}) = H^{i}(\Sigma,V_{|[E] \times
\Sigma}) \cong H^{i}(\Sigma,E).
\]
Tensoring $E$ with a holomorphic line bundle over $\Sigma$ of
degree $e$ gives an isomorphism between $\Mnd$ and ${\cal M}(n,d+ne)$.
Since $n$ and $d$ are coprime we may assume without any loss of
generality that 
\[
(2g-2) n < d <(2g-1) n.
\]
From \cite[lemma 5.2]{N2} we know that $H^{1}(\Sigma,E)=0$
for any semistable holomorphic bundle $E$ of slope greater than
$2g-2$, and $H^{i}(\Sigma,E)=0$ if $i>1$ since the complex
dimension of $\Sigma$ is one.
  Thus $\pi_{!}V$ is in fact a vector bundle over $\Mnd$ with fibre
$H^{0}(\Sigma,E)$ over $[E] \in \Mnd$ and, by the Riemann-Roch theorem, of
rank $d-n(g-1)$.\\
\indent In particular, if we express the Chern classes
$c_{r}(\pi_{!}V)$ in terms of the generators $a_{r},b_{r}^{j}$ and
$f_{r}$ of $\HS(\mnd)$, then the fact that they vanish for $r>d-n(g-1)$ gives
us relations between the generators. Mumford conjectured
\cite[p.582]{AB} that when $n=2$ these relations can be used to obtain
a complete set and this was proved in \cite{K} and \cite{Z}. For $n
\geq 3$ the relations obtained from the vanishing of the Chern classes
of $\pi_{!}V$ above its rank do not in general form a complete set
\cite[p.19]{E}. However one can also obtain {\em dual Mumford
relations} from vanishing Chern classes of $\pi_{!}(V^*)$, because
$H^1(\Sigma,E^*)=0$ for any semistable bundle $E$ over $\Sigma$ of
rank $n$ and degree $d$ satisfying $(1-2g)n<d<(2-2g)n$, and in
\cite{E} it is shown that when $n=3$ the Mumford relations together
with the dual Mumford relations form a complete set. In this paper we
shall further modify Mumford's relations to obtain a complete set for
general $n$.

\begin{rem} We could also choose $d$ such that $-n<d<0$, which
implies that $H^0(\Sigma,E)=0$ for any semistable bundle of rank $n$
and degree $d$ over $\Sigma$. Then $-\pi_!V$ is a bundle over
$\mnd$ of rank $n(g-1)-d$, so that $c_{r}(-\pi_!V)=0$ if
$r>n(g-1)-d$. By Serre duality the relations obtained are equivalent
to the dual Mumford relations, while the relations obtained by choosing $d$
satisfying $0<d<n$, so that $-\pi(V^*)$ is a bundle of rank $n(g-1)+d$,
are equivalent to the original Mumford relations.
\end{rem}

\section{Statement of Results}

\renorm

Suppose now that $\n$ and $\D$ are integers satisfying $0<\n<n$ and 
\[
\frac{\D}{\n} > \frac{d}{n}.
\]
Then, as in $\S$1, we have a $\G(n,d)$-equivariant bundle ${\cal V}$
over $\C(n,d) \times \Sigma$ and a $\G(\n,\D)$-equivariant bundle
$\hat{{\cal V}}$ over $\C(\n,\D) \times \Sigma$. (Note that as it is
now important to specify the values of $n$ and $d$ we denote $\C$
and $\G$ by  $\C(n,d)$ and $\G(n,d)$.) Both ${\cal V}$ and $\hat{{\cal
V}}$ can then be pulled back to $\C(\n,\D) \times \C(n,d) \times \Sigma$.\\
\indent Let  $\pi:\C(\n,\D) \times \C(n,d) \times \Sigma \to \C(\n,\D)
\times \C(n,d)$ be the projection. As $\D/\n > d/n$ there are no
nonzero holomorphic bundle maps from a semistable bundle $\hat{E}$ of
rank $\n$ and degree $\D$ to a semistable bundle $E$ of
rank $n$ and degree $d$, and hence by the Riemann-Roch theorem we have
\[
H^0(\Sigma, \hat{E}^* \otimes E) = 0 \mbox{ and } \dim H^1(\Sigma,
\hat{E}^* \otimes E) = \n n (g-1) - d \n + \D n.
\]
Therefore the restriction to $\C(\n,\D)^{ss} \times \C(n,d)^{ss}$ of
$-\pi_!(\hat{\cal V}^* \otimes {\cal V})$ is a bundle of rank 
\[\n n (g-1) - d \n + \D n.
\]
This means that its equivariant Chern classes
\begin{equation}
\label{chcl}
c_{r}(-\pi_!(\hat{\cal V}^* \otimes {\cal V})) \in
H^*_{\G(n,d)}(\C(n,d)) \otimes H^*_{\G(\n,\D)}(\C(\n,\D))
\end{equation}
restrict to zero in $\C(n,d)^{ss} \times \C(\n,\D)^{ss}$ when $r > \n n
(g-1) - d \n + \D n$. Therefore if we restrict the equivariant Chern
class $c_{r}(-\pi_!(\hat{\cal V}^* \otimes {\cal V}))$ to
$\C(n,d)^{ss} \times \C(\n,\D)^{ss}$ for such a value of $r$, and take its slant product with
any element of $H_*^{\G(\n,\D)}(\C(\n,\D)^{ss})$, this gives us an
element of the kernel of the restriction map
\[
H^*_{\G(n,d)}(\C(n,d)) \to H^*_{\G(n,d)}(\C(n,d)^{ss}) \cong
\HS(\mnd) \otimes \HS(BU(1)).
\]
The main aim of this paper is to show that the kernel of this
restriction map is generated as an ideal in $H^*_{\G(n,d)}(\C(n,d))$
by the elements obtained in this way. Equivalently these relations,
together with the relation (\ref{norm}) due to the normalisation of
the universal bundle, give a complete set of relations among the
generators $\ar,\bjr$ and $\fr$ of $\HS(\mnd)$.\\
\indent We need to bear in mind that, although by assumption $n$ and
$d$ are coprime, we are not assuming that $\n$ and $\D$ are
necessarily coprime. When $\n$ and $\D$ are coprime we can express the
relations obtained above in terms of universal bundles $V$ and $\V$
over the moduli spaces $\mnd$ and $\M$ as in the introduction, and
take slant products with elements of $H_*(\M)$. However when $\n$ and
$\D$ are not coprime we need to take care about the distinction
between stable and semistable bundles, and the bundle $\hat{\cal V}$
over $\C(\n,\D) \times \Sigma$ does not in general induce a universal
bundle $\V$, even over the product with $\Sigma$ of the moduli space 
${\cal M}^s(\n, \D) \cong \C(\n,\D)^s
/\G(\n,\D)$ of stable bundles over $\Sigma$ \cite[Chapter 5]{N2}. 
Thus the precise statement of our main result is phrased in terms
of ${\cal C}(\n, \D)$, not ${\cal M}(\n,\D)$, as follows.

\begin{thm}
\label{thm}
If $n$ and $d$ are coprime then the kernel of the restriction map 
\[
H^*_{\G(n,d)}(\C(n,d)) \to H^*_{\G(n,d)}(\C(n,d)^{ss})
\]
is generated as an ideal in $H^*_{\G(n,d)}(\C(n,d))$ by slant
products of the form 
\[
c_{r}(-\pi_!(\hat{\cal V}^* \otimes {\cal V})) \backslash \gamma
\] 
for integers $\n,\D$ and $r$ satisfying $0<\n<n$ and $\frac{d}{n} <
\frac{\D}{\n} < \frac{d}{n}+1$ and 
\[
n\n(g-1)-d\n+\D n <r < n\n(g+1)-d\n+\D n
\]
where 
\[
c_{r}(-\pi_!(\hat{\cal V}^* \otimes {\cal V})) \in
H^*_{\G(\n,\D)}(\C(\n,\D)) \otimes H^*_{\G(n,d)}(\C(n,d))
\]
is the $r$th equivariant Chern class of the equivariant virtual bundle
$-\pi_!(\hat{\cal V}^* \otimes {\cal V})$ over $\C(\n,\D) \times
\C(n,d)$, and $\gamma$ lies in the image of the natural map
\[ H_*^{\G(\n,\D)}(\C(\n,\D)^{ss}) \to H_*^{\G(\n,\D)}(\C(\n,\D)). \]
\end{thm}
We shall rephrase this theorem more explicitly in terms of a complete set
of relations among the generators $\ar,\bjr,\fr$ for $\HS(\mnd)$ 
 in $\S$ 10 (see Theorems \ref{explicit} and
\ref{explicit2} and Remark \ref{martin}).

\section{Completeness Criteria}

\renorm

We have observed that the relations between the generators $\ar,\bjr$
and $\fr$ for $\HS(\mnd)$ are determined by the kernel of the
restriction map
\[
\HG(\C) \cong  H^{*}_{\overline{\G}}(\C) \otimes
\HS(BU(1)) \rightarrow H^{*}_{\overline{\G}}(\Css) \otimes
\HS(BU(1)) \cong \HG(\Css).
\]
In order to describe this kernel we consider the Harder-Narasimhan-Shatz stratification of
the space $\C$ of holomorphic structures on a fixed $C^\infty$
complex vector bundle ${\cal E}$ of rank $n$ and degree $d$ over
$\Sigma$ \cite{AB,HN,Sh}. The stratification $\{\Cmu : \mu \in {\cal M} \}$ is
indexed by the partially ordered set ${\cal M}$ consisting of all the
types of holomorphic bundles of rank $n$ and degree $d$, as follows.\\
\indent Any holomorphic bundle $E$ over $M$ of rank $n$ and degree $d$ has a
filtration (or flag) \cite[p.221]{HN},\cite{LP}
\[
0 = E_{0} \subset E_{1} \subset \cdot \cdot \cdot \subset E_{s} = E
\]
of subbundles such that the quotient bundles $Q_{p}= E_{p}/E_{p-1}$
are semistable for $1 \leq p \leq s$ and satisfy 
\[
\mu(Q_{p}) = \frac{d_{p}}{n_{p}} > \frac{d_{p+1}}{n_{p+1}} =
\mu(Q_{p+1})
\] 
where $d_{p}$ and $n_{p}$ are respectively the degree and rank of
$Q_{p}$, and $\mu(Q_p) = d_p/n_p$ is its slope. This filtration is canonically 
associated to $E$ and is called the Harder-Narasimhan filtration of $E$. We define the type of $E$ to be
\[
\mu = (\mu(Q_{1}),...,\mu(Q_{s})) \in \Q^{n}
\]
where the entry $\mu(Q_{p})$ is repeated $n_{p}$ times. The
semistable bundles have type 
\[
\mu_{0} = (d/n,...,d/n)
\]
and form the
unique open stratum. The set ${\cal M}$ of all possible types of holomorphic
vector bundles over $\Sigma$ provides an indexing set for the stratification, and 
if $\mu \in {\cal M}$ then the subset $\Cmu \subseteq {\cal C}$ is
defined to be the set of all holomorphic vector bundles over $\Sigma$ of type
$\mu$. A partial order on
${\cal M}$ is defined as follows. Let $\sigma=(\sigma_{1},...,\sigma_{n})$ and
$\tau=(\tau_{1},...,\tau_{n})$ be two types; we define 
\begin{equation} \label{po}
\sigma \geq \tau \mbox{ if and only if }
\sum_{j \leq i} \sigma_{j} \geq \sum_{j \leq i} \tau_{j} \mbox{ for } 1 \leq i
\leq n-1.
\end{equation}
The stratification has the following properties:-\\

\indent (i) The stratification is smooth and $\G_c$-invariant. That is, each stratum $\Cmu$ is a
locally closed $\G_{c}$-invariant complex submanifold of $\C$. Furthermore, by
\cite[7.8]{AB}, the closure in ${\cal C}$ of the stratum ${\cal C}_{\mu}$ for any
$\mu \in {\cal M}$ satisfies\\
\begin{equation}
\overline{\Cmu} \subseteq \bigcup_{\nu \geq \mu} \C_{\nu}. \label{11}
\end{equation}
\indent (ii) Each stratum $\Cmu$ is connected and has finite (complex)
codimension $d_{\mu}$ in $\C$, given by the formula \cite[7.16]{AB}
\begin{equation}
d_{\mu}= \sum_{i>j} (n_{i}d_{j}-n_{j}d_{i}+n_{i}n_{j}(g-1)), \label{12}
\end{equation}\\ 
when $\mu = (d_1/n_1,...,d_P/n_P)$ as above. Moreover given any
integer $N$ there are only finitely many $\mu \in {\cal M}$ such that
$d_{\mu} \leq N$.\\
\indent (iii) The gauge group $\G$ acts on $\C$ preserving the
stratification which is equivariantly
perfect with respect to this action \cite[Thm. 7.14]{AB}. In
particular there is an isomorphism of vector spaces
\[
H^{k}_{\G}(\C) \cong \bigoplus_{\mu \in {\cal M}}
H^{k-2d_{\mu}}_{\G}(\Cmu) = H^{k}_{\G}(\Css) \oplus \bigoplus_{\mu \neq
\mu_{0}} H^{k-2d_{\mu}}_{\G} (\Cmu).
\]
The restriction map $H^{k}_{\G}(\C) \rightarrow H^{k}_{\G}(\Css)$ is the
projection onto the summand $H^{k}_{\G}(\Css)$ and so the kernel is
isomorphic {\em as a vector space} to
\begin{equation}
\bigoplus_{\mu \neq \mu_{0}} H^{k-2d_{\mu}}_{\G} (\Cmu). \label{13}
\end{equation}
\indent (iv) Let $\mu=(d_{1}/n_{1},...,d_{s}/n_{s})$. Atiyah and
Bott show \cite[Prop. 7.12]{AB} that the map
\[
\prod_{p=1}^{s} \C(n_{p},d_{p})^{ss} \rightarrow \Cmu,
\]
which sends a sequence of semistable bundles $(F_{1},...,F_{s})$ to
the direct sum $F_{1} \oplus \cdot \cdot \cdot \oplus F_{s}$, induces an
isomorphism
\[
\HG(\Cmu) \cong \bigotimes_{1 \leq p \leq s} H^{*}_{{\cal
G}(n_{p},d_{p})}({\cal C}(n_{p},d_{p})^{ss}).
\]
\begin{rem}
\label{two}
It follows from the isomorphism (\ref{13}) that for any set of relations in
$\HG(\Css)$ to be
complete it is necessary that the least degree of a nontrivial relation in the set must be
equal to the smallest real codimension of an unstable
stratum. When $n>2$ and $n/2 < d < n$ then the least real codimension
of an unstable stratum is $2(n-d +(n-1)\g)$ and the smallest degree
of any of the original Mumford relations is $2(d + (n-1)\g)$ \cite[p.19]{E}, and
consequently the set of original Mumford relations is not generally complete.
\end{rem}
\begin{rem}
\label{three}
By property (i) above, the subset
\[
U_\mu = \C - \bigcup_{\nu > \mu} \Cnu = \bigcup_{\nu \not > \mu} \Cnu
\]
is an open subset of $\C$ containing $\Cmu$ as a closed
submanifold. Hence there is a long exact sequence (the Thom-Gysin sequence)
\[
\cdots \to H_{\G}^{j-2d_\mu}(\Cmu) \to H^j_{\G}(U_\mu) \to
H^j_{\G}(U_\mu-\Cmu) \to \cdots .
\]
From (iii) we know that the stratification is equivariantly perfect,
which means that each of these long exact sequences break up into short
exact sequences
$ 0 \to H_{\G}^{j-2d_\mu}(\Cmu) \to H^j_{\G}(U_\mu) \to
H^j_{\G}(U_\mu-\Cmu) \to 0$, and so 
the maps $\HG(\C) \to \HG(U_\mu)$ are all surjective. In fact
Atiyah and Bott show that the stratification is equivariantly perfect
by considering the composition 
\[
H_{\G}^{j-2d_\mu}(\Cmu) \to H^j_{\G}(U_\mu) \to H^j_{\G}(U_\mu-\Cmu)
\]
which is multiplication by $e_\mu$, the equivariant Euler class of the
normal bundle to $\Cmu$ in $\C$. Atiyah and Bott show \cite[p.569]{AB}
that $e_\mu$ is not a zero-divisor in $\HG(\Cmu)$ and deduce that the
Thom-Gysin maps 
\[
H_{\G}^{j-2d_\mu}(\Cmu) \to H^j_{\G}(U_\mu)
\]
are all injective.
\end{rem}

\indent This leads to certain criteria for a set of relations to be
complete, described in \cite[Prop. 1]{K} (see also \cite[Prop. 4]{E} where
the criteria are slightly modified and the proof of completeness is
corrected). First we introduce a total order $\preceq$ on ${\cal M}$
extending the partial order $\leq$ defined at (\ref{po}) above. This
total order $\preceq$ is a lexicographic ordering: given
$\mu=(\mu_{1},...,\mu_{n})$ and
$\nu=(\nu_{1},...,\nu_{n})$ in ${\cal M}$ we define
\[
\mu \prec \nu \mbox{ if there exists $q \in \{1,...,n\}$ with } 
\mu_{q} < \nu_{q} \mbox{ and } \mu_{i} = \nu_{i} \mbox{
for } 1 \leq i < q.
\]
We then write $\mu \preceq \nu$ if $\mu \prec \nu$ or
$\mu=\nu$. \\
\indent Following \cite{Sh} and \cite[$\S$ 7]{AB} we associate with
$\mu = (d_{1}/n_{1}, \ldots, d_{s}/n_{s})$ the convex polygon
$P_{\mu}$ with vertices 
\[
(0,0),(n_1,d_1),(n_1+n_2,d_1+d_2), \ldots ,
(n_1 + \cdots +n_s,d_1+\cdots +d_s)= (n,d). 
\]

We see that $\mu \leq \nu$ if and only if $P_\mu$ lies below $P_\nu$
whereas $\mu \preceq \nu$ if and only if $P_\mu$ and $P_\nu$ agree to
the left of some vertex and $P_\mu$ lies below $P_\nu$
immediately to the right of this vertex. 
Notice that the total order $\preceq$ introduced here
differs from that of \cite[p.20]{E} and is in a sense its dual;
 in \cite{E} $\preceq$ is defined by $\mu \preceq \nu$
if and only if $P_{\mu}$ and $P_{\nu}$ agree to the {\em right} of
some vertex and $P_{\mu}$ lies below $P_{\nu}$ immediately to the
left of this vertex.\\
\indent By \cite[Prop. 4]{E} we have:

\begin{prop}
\label{KCC}
(Completeness Criteria) Let ${\cal R}$ be a
subset of the kernel of the restriction map
\[
\HG(\C) \rightarrow \HG(\Css).
\]
Suppose that for each unstable type $\mu \neq \mu_{0}$ there is a subset ${\cal
R}_{\mu}$ of the ideal generated by ${\cal R}$ in $\HG(\C)$ such that
the image of ${\cal R}_{\mu}$ under the restriction map
\[
\HG(\C) \rightarrow \HG(\C_{\nu})
\]
is zero when $\nu \prec \mu$ and when $\nu = \mu$ contains the ideal of
$\HG(\C_{\mu})$ generated by the equivariant
Euler class  $e_{\mu}$ of the normal bundle ${\cal N}_{\mu}$ to the stratum $\C_{\mu}$ in
$\C.$ Then ${\cal R}$ generates the kernel of the restriction map
\[
\HG(\C) \rightarrow \HG(\Css)
\]
as an ideal of $\HG(\C).$
\end{prop}

\indent In order to simplify later arguments, we now observe that
Atiyah and Bott could have used a coarser stratification of
$\C.$
\begin{df}
\label{defn}
For any integers $n_1$ and $d_1$ let $S_{n_1,d_1}$ be the
subset of $\C$ consisting of all those holomorphic structures with
Harder-Narasimhan filtration 
\[
0 = E_0 \subset E_1 \subset \cdots \subset E_s = E
\]
\end{df}
where $E_1$ has rank $n_1$ and degree $d_1$. We shall say that such a
holomorphic structure has {\em coarse type} $(n_1,d_1)$.
When it is necessary to specify the dependence of $S_{n_1,d_1}$ on
$n$ and $d$ we shall write $S_{n_1,d_1}^{(n,d)}$.
Then
\[
S_{n,d} = \Css = \C_{\mu_0}
\]
and each $S_{n_1,d_1}$ is a union of finitely many strata from the
stratification $\{\Cmu: \mu \in {\cal M}\}$ of $\C$ by type. Moreover since
\[
\overline{\Cmu} \subseteq \bigcup_{\mu' \geq \mu} \C_{\mu'}
\]
we have 
\[
\overline{S_{m,k}} \subseteq \left( \bigcup_{\frac{k'}{m'} > \frac{k}{m}}
S_{m',k'} \right) \cup \left( \bigcup_{\frac{k'}{m'} = \frac{k}{m}, k'
\geq k} S_{m',k'} \right). 
\]
Now $\Cmu$ is locally a submanifold
of finite codimension 
\[
d_{\mu} = \dim H^1(\Sigma, \mbox{End}''E) = \sum_{i>j} n_i d_j -n_j d_i
+ n_i n_j (g-1) 
\]
where $E$ is any holomorphic bundle of type $\mu$ and $\mbox{End}''E$
is the quotient of $\mbox{End}E$ by the subbundle $\mbox{End}'E$ of
holomorphic endomorphisms of $E$ which preserve its Harder-Narasimhan
filtration (see \cite[$\S$7]{AB}). The proof of this is based on the
observation that  
\[
H^0(\Sigma, \mbox{End}''E) = 0
\]
which is a direct consequence of the fact that if $D_1$ and $D_2$ are
semistable bundles of rank $n_1$ and $n_2$ and degrees $d_1$ and $d_2$
satisfying $d_1/n_1 > d_2/n_2$ then there are no nonzero bundle maps
from $D_1$ to $D_2$. If we replace $\mbox{End}'E$ by the bundle of
holomorphic endomorphisms of $E$ which preserve the first proper
subbundle $E_1$ in its Harder-Narasimhan filtration and replace $\mbox{End}''E$
by the quotient of $\mbox{End}E$ by this subbundle, then we still have
$H^0(\Sigma, \mbox{End}''E) = 0,$ and the argument of \cite[$\S$
7]{AB} that $\Cmu$ is locally a submanifold of finite dimension
$d_\mu$ in $\C$ applies equally well to show that $S_{n_1,d_1}$ is
locally a submanifold of finite codimension 
\begin{equation} \label{codim*}
\begin{array}{rcl} \delta_{n_1,d_1} & = & (n-n_1)d_1 - n_1(d-d_1) + n(n-n_1)(g-1) \\
& = & n d_1 - n_1 d + n_1 (n-n_1)(g-1)
\end{array} \end{equation}
in $\C$.\\
\indent Note that if $d_1/n_1 > d/n$ then a holomorphic bundle $E$
has 
 coarse type
$(n_1,d_1)$ if and only if $E$ has a
semistable subbundle $E_1$ of rank $n_1$ and degree $d_1$ such that
the coarse type $(n_2,d_2)$ of the quotient bundle $E/E_1$ satisfies
$d_2/n_2 < d_1/n_1$. Thus the proof of \cite[Prop. 7.12]{AB} shows
that
\begin{equation}
\label{tensor}
\quad \HG(S_{n_1,d_1}) \cong \HS_{\G(n_1,d_1)}(\C(n_1,d_1)^{ss}) \otimes 
 \HS_{\G(n-n_1,d-d_1)}\left( U(n_1,d_1) \right)
\end{equation}
where 
\[ U(n_1,d_1) = 
\bigcup_{\frac{d_2}{n_2} < \frac{d_1}{n_1}} S^{(n-n_1,d-d_1)}_{n_2,d_2}
= \C(n-n_1,d-d_1) \backslash \bigcup_{\frac{d_2}{n_2} \geq \frac{d_1}{n_1}}
S^{(n-n_1,d-d_1) }_{n_2,d_2}
\]
is an open subset of $\C(n-n_1,d-d_1)$. 

The isomorphism (\ref{tensor}) is induced by inclusions 
\[
\C(n_1,d_1)^{ss} \times U(n_1,d_1) \hookrightarrow 
\C(n_1,d_1) \times \C(n-n_1,d-d_1) \hookrightarrow \C(n,d) = \C
\]
and
\[
\G(n_1,d_1) \times \G(n-n_1,d-d_1) \hookrightarrow \G(n,d) = \G,
\]
given by identifying our fixed $C^\infty$ bundle of rank $n$ and
degree $d$ over $\Sigma$ with the direct sum of fixed $C^\infty$ bundles of
ranks $n_1$ and $n-n_1$ and degrees $d_1$ and $d-d_1$
respectively. The gauge group $\G(n,d)$ has a constant central circle subgroup
acting trivially on $\C(n,d)$, and therefore $
\G(n_1,d_1) \times \G(n-n_1,d-d_1)$ has a central subgroup $T=(S^1)^2$
acting trivially on $\C(n_1,d_1) \times \C(n-n_1,d-d_1)$. Moreover at
a point of 
$ \C(n_1,d_1)^{ss} \times U(n_1,d_1)$
corresponding to a direct sum of bundles $E=D_1 \oplus D_2$, the
normal to $S_{n_1,d_1}$ is naturally isomorphic to 
\begin{equation}
H^1(\Sigma,\mbox{End}''E) = H^1(\Sigma,D_1^* \otimes D_2) \label{normal}
\end{equation}
(cf. \cite[$\S$7]{AB}) and an element $(t_1,t_2) \in (S^1)^2 = T$ acts on this as scalar
multiplication by $t_1^{-1}t_2.$ Thus the representation of $T$ on
the normal to $S_{n_1,d_1}$ is primitive, and so by Atiyah and Bott's
criterion \cite[13.4]{AB} the equivariant Euler class $e_{n_1,d_1}$ of
the normal to $S_{n_1,d_1}$ in $\C$ is not a zero divisor in
$\HG(S_{n_1,d_1})$. Therefore by \cite[1.9]{AB} we have proved
(cf. \cite[Thm. 7.14]{AB}),
\begin{prop}
\label{ep}
The stratification
\[
\left\{ S_{n_1,d_1}: 0 < n_1 < n, \frac{d_1}{n_1} > \frac{d}{n}
\right\} \bigcup \left\{ S_{n,d}\right\} 
\]
of $\C$ by coarse type is equivariantly perfect.
\end{prop} 
Let $P_{G}(X)(t)$ denote the equivariant Poincar\'{e} series
\[
P_{G}(X)(t) = \sum_{i \geq 0} t^i \dim H^i_{G}(X)
\]
of a topological space $X$ acted on by a group $G$. Then Proposition
\ref{ep} together with (\ref{tensor}) tells us that 
\[
P_{\G}(\C)(t) = P_{\G}(\Css)(t) + \sum_{{ \begin{array}{l}
0<n_1<n,\\ 
d_1/n_1 > d/n \end{array}}} t^{2
\delta_{n_1,d_1}} P_{\G}(S_{n_1,d_1})(t)
\]
where 
\begin{eqnarray*}
P_{\G}(S_{n_1,d_1})(t) =  P_{\G(n_1,d_1)}(\C(n_1,d_1)^{ss})(t)
P_{\G(n-n_1,d-d_1)}(U(n_1,d_1))(t)\\
= P_{\G(n_1,d_1)}(\C(n_1,d_1)^{ss})(t) \left[
P_{\G(n-n_1,d-d_1)}(\C(n-n_1,d-d_1))(t) \right.
\\ \left. - \sum_{\frac{d_2}{n_2} \geq
\frac{d_1}{n_1}} t^{2\delta(n-n_1,d-d_1)_{n_2,d_2}} P_{\G(n-n_1,d-d_1)}(
S^{(n-n_1,d-d_1)}_{n_2,d_2})(t) \right],
\end{eqnarray*}
with
\[
\delta_{n_1,d_1} = nd_1 - n_1 d+n_1(n-n_1)(g-1)
\]
and
\begin{eqnarray*}
\delta^{(n-n_1,d-d_1)}_{n_2,d_2} & = & (n-n_1-n_2)d_2 - n_2(d-d_1-d_2) + n_2
(n-n_1 - n_2)(g-1)\\
& = & (n-n_1)d_2 - n_2(d-d_1) +n_2(n-n_1 - n_2)(g-1).
\end{eqnarray*}
Since \cite[Thm. 2.15]{AB}
\[
P_{\G}(\C)(t) =
\frac{\prod_{k=1}^{n}(1+t^{2k-1})^{2g}}{\prod_{k=1}^{n}(1+t^{2k})
\prod_{k=2}^{n}(1+t^{2k-2})},
\]
this gives us an inductive formula for $P_{\G}(\Css)(t)$, and hence
for the Betti numbers of $\mnd$ when $n$ and $d$ are coprime. This is
easily seen to be equivalent to the formula 
\begin{equation}
P_{\G}(\C)(t) = P_{\G}(\Css)(t) + \sum_{\mu = (d_1/n_1, \ldots ,
d_s/n_s) \neq \mu_0} t^{2d_\mu} \prod_{1 \leq j \leq s}
P_{\G(n_j,d_j)}(\C(n_j,d_j)^{ss})(t) \label{p16}
\end{equation}
of Atiyah and Bott \cite[Prop. 7.12, Thm. 7.14]{AB} since 
\begin{eqnarray*}
d_{\mu} & = & \sum_{i>j} n_i d_j -n_j d_i + n_i n_j (g-1)\\
& = & \sum_{j=1}^{s-1} \left\{ \left( \sum_{i=j+1}^{s} n_i\right) d_j - 
 n_j \left( \sum_{i=j+1}^{s} d_i \right) + n_j \left( \sum_{i=j+1}^{s}
n_i\right) (g-1) \right\}\\
& = &\delta_{n_1,d_1} + \delta^{(n-n_1,d-d_1)}_{n_2,d_2} + \delta^{(n_{s-1} +
n_s, d_{s-1}+d_{s})}_{n_s,d_s}.
\end{eqnarray*}
This inductive formula is equivalent to the following closed
formula \cite{DB,LR,Z2}
$$P(\mnd)(t) = \sum_{n_1 + \cdots n_l = n} (-1)^{l-1} \frac{(1+t)^{2gl}}{(1-t^2)^{l-1}}
\prod_{j=1}^l \prod_{i=1}^{n_j-1} \frac{(1+t^{2i+1})^{2g}}{(1-t^{2i})(1-t^{2i+2})} \times$$
\begin{equation} \label{closedf} 
\times \prod_{j=1}^{l-1} \frac{t^{2 \sum_{i<j} n_in_j(g-1) + \sum_{i=1}^{l-1}(n_i + n_{i+1})
\langle -(n_1 + \cdots + n_j)d/n \rangle}}{1 - t^{2(n_j + n_{j+1})}}.
\end{equation}
The fact that the equivariant Euler class $e_{n_1,d_1}$ of the normal
to any stratum $S_{n_1,d_1}$ in $\C$ is not a zero-divisor in
$\HG(S_{n_1,d_1})$ also gives us an alternative (and slightly simpler)
version of the completeness criteria (Proposition \ref{KCC} above) as follows.

\begin{prop}
\label{MCC}
(Modified Completeness Criteria) Let ${\cal R}$ be a
subset of the kernel of the restriction map
\[
\HG(\C) \rightarrow \HG(\Css).
\]
Suppose that for each pair of integers $(\n,\D)$ with $0<\n<n$ and
$\frac{\D}{\n} > \frac{d}{n}$ there is a subset ${\cal
R}_{\n,\D}$ of the ideal generated by ${\cal R}$ in $\HG(\C)$ such that
the image of ${\cal R}_{\n,\D}$ under the restriction map
\[
\HG(\C) \rightarrow \HG(S_{n_1,d_1})
\]
is zero when $d_1/n_1<\D/\n$ or $d_1/n_1=\D/\n$ and $n_1<\n$ and when
$(n_1,d_1) = (\n,\D)$ equals the ideal of
$\HG(S_{\n,\D})$ generated by the equivariant Euler class $e_{\n,\D}$
of the normal bundle to the stratum $S_{\n,\D}$ in
$\C.$ Then ${\cal R}$ generates the kernel of the restriction map
\[
\HG(\C) \rightarrow \HG(\Css)
\]
as an ideal of $\HG(\C).$
\end{prop}

\noindent{\em Proof:} Suppose that $\eta \in H^j_{\G}(\C)$ lies in the kernel
of the restriction map $
\HG(\C) \rightarrow \HG(\Css)$. There are only finitely many pairs of integers
$(\n,\D)$ with $0<\n <n$ and $\D/\n > d/n$ for which the real codimension
$$2 \delta_{\n,\D} = 2(n\D - \n d + \n(n-\n)(g-1))$$
of the stratum $S_{\n,\D}$ is at most $j$; let these be
$$(\n_1,\D_1), (\n_2,\D_2), \ldots, (\n_k,\D_k)$$
ordered so that $(\n_1,\D_1) = (n,d)$ and
$$\D_i/\n_i < \D_{i+1}/\n_{i+1} \quad \mbox{ or } \quad \D_i/\n_i = \D_{i+1}/\n_{i+1} \mbox{ and } \n_i < \n_{i+1}$$
for $1 \leq i <k$. Then by Proposition 3.5 we have short exact sequences
\begin{equation} \label{TG*} 0 \to H_\G^{j-2\delta_{\n_i,\D_i}}(S_{\n_i,\D_i}) \to H_\G^j(U(\n_i,\D_i) \cup S_{\n_i, \D_i})
\to H^j_\G(U(\n_i,\D_i)) \to 0 \end{equation}
for $1 \leq i <k$, where $U(\n_1,\D_1) = \Css$ and $U(\n_i,\D_i) \cup S_{\n_i, \D_i}=
U(\n_{i+1},\D_{i+1}) $ and the restriction map
$$H^j_\G(\C) \to H^j_\G(U(\n_k,\D_k) \cup S_{\n_k, \D_k})$$
is an isomorphism.

By assumption $\eta$ lies in the kernel of the restriction map $
\HG(\C) \rightarrow \HG(\Css)$ which is the composition of the restriction maps
$$H^j_\G(\C) \to H^j_\G(U(\n_k,\D_k) \cup S_{\n_k, \D_k}) \to
H^j_\G(U(\n_k,\D_k)) $$
$$\to H^j_\G(U(\n_{k-1},\D_{k-1})) \to \cdots \to H^j_\G(U(\n_1,\D_1)) = H_\G^j(\Css).$$
Therefore its image $\eta_1$ in $H^j_\G(U(\n_2,\D_2)) = H^j_\G(U(\n_1,\D_1)\cup S_{\n_1, \D_1})$ lies
in the kernel of the restriction to $H^j_\G(U(\n_1,\D_1))$, which by (\ref{TG*}) is the image of the Thom-Gysin
map
$$TG: H_\G^{j-2\delta_{\n_1,\D_1}}(S_{\n_1,\D_1}) \to H_\G^j(U(\n_1,\D_1) \cup S_{\n_1, \D_1}).$$
Let $\zeta_1 \in H_\G^{j-2\delta_{\n_1,\D_1}}(S_{\n_1,\D_1})$ be such that $\eta_1 = TG(\zeta_1)$; then by 
hypothesis there is some $\sigma_1$ in the ideal generated by ${\cal R}$ whose restriction to $\Css$ is 0
and whose restriction to $S_{\n_1,\D_1}$ is $\zeta_1 e_{\n_1,\D_1)}$. Since $\sigma_1$ restricts to 0
in $\Css = U(\n_1,\D_1)$ it restricts to $TG(\xi_1)$ in $U(\n_2,\D_2)$ for some $\xi_1 \in H_\G^{j-
2\delta_{\n_1,\D_1}}(S_{\n_1,\D_1})$. But then the restriction of $\sigma_1$ to $S_{(\n_1,\D_1)}$ is
$\zeta_1 e_{\n_1,\D_1)} = \xi_1 e_{\n_1,\D_1)}$, and as $e_{\n_1,\D_1)}$ is not a zero-divisor in
$H_\G^*(S_{\n_1,\D_1})$ it follows that $\zeta = \xi$. Thus $\eta$ and $\sigma_1$ both restrict to
$TG(\zeta_1)$ in $H^*_\G(U(\n_2,\D_2))$, so the restriction of $\eta - \sigma_1$ to $H^j_\G(U(\n_3,\D_3))$
lies in the image of the Thom-Gysin map
$$TG: H_\G^{j-2\delta_{\n_2,\D_2}}(S_{\n_2,\D_2}) \to H_\G^j(U(\n_2,\D_2) \cup S_{\n_2, \D_2}).$$
Repeating the argument above with $S_{\n_2,\D_2}$ replacing $S_{\n_1,\D_1}$ we find that there is some
$\sigma_2$ in the ideal generated by ${\cal R}$ such that $\eta - \sigma_1 -\sigma_2$ restricts to
0 in $H^j_\G(U(\n_3,\D_3))$. Repeating the same argument $k$ times gives us $\sigma_1, \ldots, \sigma_k$
in the ideal generated by ${\cal R}$ such that $\eta - \sigma_1 - \cdots -\sigma_k$ restricts to 0 in
$$H^j_\G(U(\n_k,\D_k) \cup S_{\n_k, \D_k}) \cong H^j_\G(\C) ,$$
and so $\eta = \sigma_1 + \cdots + \sigma_k$ lies in the ideal generated by ${\cal R}$ as required.

\section{Diagonals}

\renorm

The modified completeness criteria of Proposition \ref{MCC} tell us
that for each pair of positive integers $(\n,\D)$ with $0<\n<n$ and
$\frac{\D}{\n} > \frac{d}{n}$ we need to look for a set of relations
whose restriction in $\HG(S_{n_1,d_1})$ is zero when $d_1/n_1<\D/\n$ or
$d_1/n_1=\D/\n$ and $n_1<\n$, and when $(n_1,d_1) = (\n,\D)$ equals
the ideal of $\HG(S_{\n,\D})$ generated by the equivariant Euler class
$e_{\n,\D}$ of the normal to $S_{\n,\D}$ in $\C$. 
 Recall from
(\ref{tensor}) that  
\[
\HG(S_{n_1,d_1}) \cong \HS_{\G(n_1,d_1)}(\C(n_1,d_1)^{ss}) \otimes
\HS_{\G(n-n_1,d-d_1)} (U(n_1,d_1))
\]
where
\[
U(n_1,d_1) = \bigcup_{\frac{d_2}{n_2} < \frac{d_1}{n_1}}
S(n-n_1,d-d_1)_{n_2,d_2})
\]
is an open subset of $\C(n-n_1,d-d_1)$. Moreover it follows from
Proposition \ref{ep} and Remark \ref{three} that the restriction map 
\[
\HS_{\G(n-n_1,d-d_1)}(\C(n-n_1,d-d_1)) \to
\HS_{\G(n-n_1,d-d_1)}(U(n_1,d_1))
\]
is surjective. Thus given 
\[
\eta \in \HS_{\G(\n,\D)}(\C(\n,\D)^{ss}) \otimes
\HS_{\G(n-\n,d-\D)}(\C(n-\n,d-\D))
\]
it suffices to find a relation whose image under the restriction map
\[
\HG(\C) \to \HS_{\G(n_1,d_1)}(\C(n_1,d_1)^{ss}) \otimes
\HS_{\G(n-n_1,d-d_1)}(U(n_1,d_1)) 
\]
is zero when $d_1/n_1<\D/\n$ or $d_1/n_1=\D/\n$ and $n_1<\n$, and when
$(n_1,d_1)=(\n,\D)$ equals 
\[
\eta e_{\n,\D}
\]
(or, more precisely, the product of the  restriction of $\eta$ to
$\C(\n,\D)^{ss} \times U(\n,\D)$ with $e_{\n,\D}$.)\\
\indent By Lefschetz duality \cite[p.297]{Sp}, since $\C(\n,\D)^s
/\overline{\G}(\n,\D)= {\cal M}^s (\n,\D)$ is a manifold of dimension 
\[
D(\n,\D) = 2[(\n^2-1)(g-1) + g] = 2(\n^2 g - \n^2 +1)
\]
we have a natural map
\begin{equation} \label{ld}
LD:H_*^{\GG(\n,\D)} (\C(\n,\D)^s) \cong H_*({\cal M}^s(\n,\D)) \to 
H_{\GG(\n,\D)}^{D(\n,\D)-*}(\C(\n,\D)^{ss})
\end{equation}
such that if $\gamma \in H^{\GG(\n,\D)}_* (\C(\n,\D)^s)$ then $LD(\gamma)$ 
is the element of the dual of 
$$H^{\GG(\n,\D)}_{D(\n,\D)-*}(\C(\n,\D)^{ss})$$
 which 
takes a $\GG(\n,\D)$-equivariant  cycle on $\C(\n,\D)^{ss}$ to its intersection, 
modulo $\GG(\n,\D)$, with $\gamma$. This makes sense because 
$\GG(\n,\D)$ acts freely on $\C(\n,\D)^s$ with quotient the manifold ${\cal M}^s(\n,\D)$; 
thus modulo $\GG(\n,\D)$ we can make the intersection transverse and count the number 
of intersection points with signs in the usual way.\\ 
\indent In this section we shall prove that if $\eta = LD(\gamma)$ where
$\gamma \in H_*^{\GG(\n,\D)} (\C(\n,\D)^s)$ and we identify $ 
H_{\GG(\n,\D)}^{*}(\C(\n,\D)^{ss})$ with its image in 
\[
H^*_{\G(\n,\D)}(\C(\n,\D)^{ss}) \otimes
H^*_{\G(n-\n,d-\D)}(\C(n-\n,d-\D)^s),
\]
then the relations given by the slant products
\[
c_{r}(-\pi_!({\hat{\cal V}}^* \otimes {\cal V})) \backslash \gamma
\]
with $r \geq n \n (g-1) + \n d - n \D + 1$ will have the required
properties.
\begin{rem}
\label{six}
When $\n$ and $\D$ are coprime then $\C(\n,\D)^{ss}$ equals
$\C(\n,\D)^s$ and its quotient by $\GG(\n,\D)$, namely $\M$, is a compact
manifold. In this case the Lefschetz duality map is essentially
Poincar\'{e} duality and the map $LD$ is an isomorphism. When $\n$ and
$\D$ are not coprime we will have more work to do to find suitable relations for all $\eta \in
\HS_{\GG(\n,\D)}(\C(\n,\D)^{ss})$.
\end{rem}
\indent First we prove:
\begin{lem}
\label{seven}
If $\gamma \in H_*^{\GG(\n,\D)}(\C(\n,\D)^{s})$ where $\D/\n > d/n$, and if
$r \geq r(\n,\D)$ where $r(\n,\D) = n \n (g-1) + \n d - n \D +1$, then the image of the slant
product 
\[
c_{r}(-\pi_!(\hat{{\cal V}}^* \otimes {\cal V})) \backslash \gamma
\]
under the restriction map
\[
\HG(\C) \to \HG(S_{n_1,d_1}) \cong \HS_{\G(n_1,d_1)}(\C(n_1,d_1)^{ss})
\otimes \HS_{\G(n-n_1,d-d_1)}(U(n_1,d_1))
\]
is zero when $d_1/n_1 < \D/\n,$ and also when $d_1/n_1 = \D/\n$ and
$n_1<\n$.
\end{lem}
{\em Proof:} Let $\E$ be a stable bundle of rank $\n$ and degree $\D$,
and let $E=D \oplus C$ where $D$ is semistable of rank $n_1$ and
degree $d_1$ and $C$ has rank $n-n_1$ and degree $d-d_1$ and coarse
type $(n_2,d_2)$ where $d_2/n_2 < d_1/n_1$. If $d_1/n_1 < \D/\n$ there
are no nonzero bundle maps from $\E$ to $D$, since $\E$ and $D$ are
both semistable, and if $d_1/n_1 = \D/\n$ then, because $\hat{E}$ is stable, the only nonzero
bundle maps from $\E$ to $D$ are injections which cannot exist if
$n_1 < \n$ \cite[5.1, 5.2]{N2}. Thus $H^0(\Sigma,\E^* \otimes D)=
0$. Moreover since $\C$ has coarse type $(n_2,d_2)$ where
\[
\frac{d_2}{n_2}< \frac{d_1}{n_1} \leq \frac{\D}{\n},
\]
the Harder-Narasimhan filtration 
\[
0 = C_0 \subset C_1 \subset \cdots \subset C_s = C
\]
of $C$ has semistable subquotients $C_j/C_{j-1}$ of rank $m_j$ and
degree $c_j$ satisfying 
\[
\frac{c_j}{m_j} \leq \frac{c_1}{m_1} = \frac{d_2}{n_2} <
\frac{d_1}{n_1}.
\]
Therefore if $d_1/n_1 \leq \D/\n$, the projection of any bundle map
$\theta: \E \to C$ onto $C/C_{s-1}$ is zero, so $\theta$ maps $\E$
into $C_{s-1}$, but then the projection of $\theta$ onto $C_{s-1}/C_{s-2}$ 
is zero, so $\theta$ maps $\hat{E}$ into $C_{s-2}$, and inductively we find that $\theta=0$. Hence
$H^0(\Sigma,\E^* \otimes C) = 0$ and 
\[
H^0(\Sigma,\E^* \otimes E) = H^0(\Sigma,\E^* \otimes C) \oplus
H^0(\Sigma,\E^* \otimes D) = 0.
\]
Therefore by the Riemann-Roch theorem the restriction 
of $-\pi_{!}(\hat{{\cal V}}^* \otimes {\cal V})$ to
\[
\C(\n,\D)^{s} \times \C(n_1,d_1)^{ss} \times U(n_1,d_1)
\]
is a bundle of rank
$n\n(g-1)+\n d - n \D$; so if $r$ is strictly greater than this value
the image of the $r$th equivariant Chern class $c_r(-\pi_{!}(\hat{{\cal
V}}^* \otimes {\cal V}))$ in
\[
\HS_{\G(\n,\D)}(\C(\n,\D)^s) \otimes  \HS_{\G(n_1,d_1)}(\C(n_1,d_1)^{ss})
\otimes \HS_{\G(n-n_1,d-d_1)}(U(n_1,d_1))
\]
is zero. The result follows. 

\begin{df}
Let $\hat{a}_1, a_1^{(1)}$ and $a_1^{(2)}$ be the generators of the first, 
second and third copies of the polynomial ring $\HS(BU(1))$ in 
\[
\HS_{\G(\n,\D)}(\C(\n,\D)^s) \otimes
\HS_{\G(n_1,d_1)}(\C(n_1,d_1)^{ss}) \otimes
\HS_{\G(n-n_1,d-d_1)}(U(n_1,d_1)) 
\]
where
\begin{eqnarray*}
\HS_{\G(\n,\D)}(\C(\n,\D)^s) & \cong & \HS(BU(1)) \otimes \HS_{\GG(\n,\D)}(\C(\n,\D)^s),\\
\HS_{\G(n_1,d_1)}(\C(n_1,d_1)^{ss}) & \cong & \HS(BU(1) \otimes \HS_{\GG(n_1,d_1)}(\C(n_1,d_1)^{ss}),\\ 
\HS_{\G(n-n_1,d-d_1)}(U(n_1,d_1)) & \cong & \HS(BU(1) \otimes  \HS_{\GG(n-n_1,d-d_1)}(U(n_1,d_1)).
\end{eqnarray*}
Note that the restriction of $a_1 \in H^*_{\G(n,d)}(\C(n,d))$ to $\C(n_1,d_1)^{ss} 
\times U(n_1,d_1)$ is $a_1^{(1)} + a_1^{(2)}$.
\end{df}

\begin{df}
\label{as}
Let ${\cal V}_1$ and ${\cal V}_2$ be the universal bundles over $\C(n_1,d_1)$ 
and $\C(n-n_1,d-d_1)$, so that the restriction of ${\cal V}$ to $\C(n_1,d_1)^{ss} 
\times U(n_1,d_1)$ is ${\cal V}_1 \oplus {\cal V}_2$.
\end{df}

We need to consider the restriction of a relation of the
form $c_{r}(-\pi_!(\hat{{\cal V}}^* \otimes {\cal V})) \backslash \gamma$ to
$\HG(S_{\n,\D})$; that is, we consider the case when $(\n,\D) = (n_1,d_1)$.

\begin{prop}
\label{big}
If $\gamma \in H_*^{\GG(\n,\D)}(\C(\n,\D)^s)$ where $\D/\n > d/n$, and
if $r = r(\n,\D)+j$ where
\[
r(\n,\D) = n\n(g-1)+ \n d - n \D +1,
\]
then the image of the slant product 
 $c_{r}(-\pi_!(\hat{\cal V}^* \otimes {\cal V})) \backslash \gamma$ under
the restriction map 
\[
\HG(\C) \to \HG(S_{\n,\D}) \cong \HS_{\G(\n,\D)}(\C(\n,\D)^{ss})
\otimes \HS_{\G(n-\n,d-\D)}(U(\n,\D))
\]
is the product 
\[
(-a_1^{(1)})^j LD(\gamma) e_{\n,\D}
\]
of the equivariant Euler class $e_{\n,\D}$ of the normal bundle to
$S_{\n,\D}$ in $\C$ with the image of $\gamma$ under the Lefschetz
duality map
\[
LD:H_*^{\GG(\n,\D)} (\C(\n,\D)^s) \cong \HS({\cal M}^s(\n,\D)) \to 
H_{\GG(\n,\D)}^{D(\n,\D)-*}(\C(\n,\D)^{ss})
\]
and $j$ copies of minus the generator $a_1^{(1)} \in H^*_{\G(\n,\D)}(\C(\n,\D)^{ss})$.
\end{prop}

The proof of Proposition \ref{big} is based on the following three lemmas.

\begin{lem}
\label{lemone}
The restriction of the equivariant Chern polynomial
\[
c(-\pi_!(\hat{{\cal V}}^* \otimes {\cal V}))(t)  = \sum_{k \geq 0} c_k(-\pi_!(\hat{{\cal V}}^* \otimes {\cal V}))t^k
\]
to $\C(\n,\D)^{s} \times \C(\n,\D)^{ss} \times U(\n,\D)$ equals the product 
\[
c(-\pi_!(\hat{{\cal V}}^* \otimes {\cal V}_1))(t)  c(-\pi_!(\hat{{\cal V}}^* \otimes {\cal V}_2))(t)
\]
where $c(-\pi_!(\hat{{\cal V}}^* \otimes {\cal V}_2))(t)$ is a polynomial in $t$ of degree at most
\[
\n(n-\n)(g-1) + \n d - n \D
\]
as predicted by the Riemann-Roch theorem. Moreover the restriction of 
the Chern polynomial 
$c(-\pi_!(\hat{{\cal V}}^* \otimes {\cal V}_1))(t)$ to the complement in 
$\C(\n,\D)^{s} \times \C(\n,\D)^{ss} \times U(\n,\D)$ of $\Delta^s \times U(\n,\D)$
where
\[
\Delta^s = \{(F_1,F_2) \in \C(\n,\D)^{s} \times \C(\n,\D)^{ss} : F_1 \cong F_2 \}
\]
is a polynomial in $t$ of degree at most $\n^2(g-1)$ as predicted by Riemann-Roch.
\end{lem}

\begin{lem}
\label{lemtwo}
There is a line bundle ${\cal L}$ over $\Delta^s$ with first Chern class
\[
c_1({\cal L}) = a_1^{(1)} - \hat{a}_1
\]
such that 
\[
\left. \hat{{\cal V}} \right|_{\Delta^s \times U(\n,\D) \times \Sigma} \otimes \pi^*{\cal L}
\cong \left. {{\cal V}}_1 \right|_{\Delta^s \times U(\n,\D) \times \Sigma} 
\]
and hence 
\[
\left. \pi_!(\hat{{\cal V}}^* \otimes {\cal V}_2) \right|_{\Delta^s \times U(\n,\D)} \cong 
{\cal L} \otimes \pi_!(\left. {{\cal V}}^*_1 \otimes {\cal V}_2) \right|_{\Delta^s \times U(\n,\D)}. 
\]
Moreover the codimension of $\Delta^s$ in $\C(\n,\D)^s \times \C(\n,\D)^{s}$ is $\n^2 (g-1)+1$ 
and its normal bundle has equivariant Chern polynomial
\[
c(-\pi_!({\cal V}_1^* \otimes {\cal V}_1))(t) = 
c(-\pi_!(\hat{{\cal V}}^* \otimes \hat{{\cal V}}))(t) =
c(-{\cal L}^* \otimes \pi_!(\hat{{\cal V}}^* \otimes {\cal V}_1))(t).
\]
\end{lem}

\begin{lem}
\label{lemthree}
If $j \geq 0$ then the restriction to $\C(\n,\D)^s \times \C(\n,\D)^{ss} \times U(\n,\D)$ of 
\[
c_{\n^2(g-1)+j+1}(-\pi_!(\hat{{\cal V}}^* \times {\cal V}_1))
\]
is 
$$(\hat{a}_1 - a_1^{(1)})^j [\Delta^s \times U(\n,\D)]$$
where $[\Delta^s \times U(\n,\D)]$ is the equivariant cohomology
class in
$$H^{*}_{\G(\n,\D) \times \G(\n,\D) \times \G(n-\n,d-\D)}(\C(\n,\D)^s \times \C(\n,\D)^{ss} \times U(\n,\D))$$
dual to the closed subvariety $\Delta^s \times U(\n,\D)$ of $\C(\n,\D)^s \times \C(\n,\D)^{ss} \times U(\n,\D)$.
The restriction of
\[
c_{r(\n,\D)+j}(-\pi_!(\hat{{\cal V}}^* \otimes {\cal V}))
\]
is 
\[
e_{\n,\D} (\hat{a}_1 - a_1^{(1)})^j [\Delta^s \times U(\n,\D)]
\]
where $e_{\n,\D} = c_{\n(n-\n)(g-1)+\n d - n \D}(-\pi_!(\hat{{\cal V}}^*_1 \otimes {\cal V}_2))$
is the equivariant Euler class of the normal bundle to $S_{\n,\D}$ in $\C$.
\end{lem}

{\em Proof of Proposition \ref{big} given Lemma \ref{lemthree}:}
We need to show that if $r = r(\n,\D)+j$ and $\gamma \in  H_{*}^{\GG(\n,\D)}(\C(\n,\D)^{s})$ then the restriction to 
\[
\HS_{\G(\n,\D)}(\C(\n,\D)^{ss}) \otimes \HS_{\G(n-\n,d-\D)}(U(\n,\D))
\]
of the slant product 
\[
c_{r}(-\pi_!(\hat{\cal V}^* \otimes {\cal V})) \backslash \gamma
\] 
equals $(-a_1^{(1)})^j LD(\gamma) e_{\n,\D}$. By Lemma 
\ref{lemthree} it suffices to show that the slant product
\[
e_{\n,\D} (\hat{a}_1 - a_1^{(1)})^j [\Delta^s \times U(\n,\D)] \backslash \gamma
\]
equals $(-a_1^{(1)})^j LD(\gamma) e_{\n,\D}$, where by abuse of notation
$e_{\n,\D}$ is being used to denote cohomology classes on two different spaces.
By the definition of the slant product
\[
\HS(A \times B) \otimes H_*(B) \to \HS(A)
\]
we have
\[
\langle c \backslash b,a \rangle = \langle c, a \times b
\rangle
\]
for all $a \in H_*(A), b \in H_*(B)$ and $c \in \HS(A \times B)$,
where $\langle,\rangle$ denotes the natural pairing between cohomology
and homology \cite[p.287]{Sp}. Note that the slant product with any 
$\gamma \in  H_{*}^{\GG(\n,\D)}(\C(\n,\D)^{s})$ of any multiple of $\hat{a}_1$ 
is zero, because multiples of the generator $\hat{a}_1$ of $\HS(BU(1))$ in
\[
\HS_{\G(\n,\D)}(\C(\n,\D)^s) \cong \HS(BU(1)) \otimes \HS_{\GG(\n,\D)}(\C(\n,\D)^s)
\]
annihilate $H_*^{\GG(\n,\D)}(\C(\n,\D)^s)$. Thus we need to show that 
\[
((-a_1^{(1)})^j e_{\n,\D} 
[\Delta^s \times U(\n,\D)]  \backslash \gamma = (-a_1^{(1)})^j e_{\n,\D} LD(\gamma),
\]
or equivalently that
\[
\langle (-a_1^{(1)})^j e_{\n,\D} 
[\Delta^s \times U(\n,\D)] , a \times \gamma \rangle =
\langle (-a_1^{(1)})^j LD(\gamma) e_{\n,\D}, a \rangle
\]
for all $a \in H_*^{\G(\n,\D) \times \G(n-\n,d-\D)}(\C(\n,\D)^{ss}
\times U(\n,\D))$. But
$$
\langle (-a_1^{(1)})^j e_{\n,\D} 
[\Delta^s \times U(\n,\D)] ,a \times \gamma \rangle =
 \langle \left. (-a_1^{(1)})^j e_{\n,\D} \right|_{\Delta^s \times U(\n,\D)} ,(a \times
\gamma) \cap (\Delta^s \times U(\n,\D)) \rangle
$$
by the definition of the cohomology class $
[\Delta^s \times U(\n,\D)] $ dual to the closed subvariety $
\Delta^s \times U(\n,\D)$ of
$\C(\n,\D)^s \times \C(\n,\D)^{ss} \times U(\n,\D)
$. On the right
hand side here, we are working with the $\G(\n,\D) \times \G(\n,\D)
\times \G(n-\n,d-\D)$-equivariant cohomology and homology of $\Delta^s
\times U(\n,\D)$. As above we identify this with the tensor product of
the cohomology (respectively homology) of $BU(1) \times BU(1)$ and the 
 $\GG(\n,\D) \times \GG(\n,\D) \times \G(n-\n,d-\D)$-equivariant cohomology 
(respectively homology) of $\Delta^s \times U(\n,\D)$. The action of $\GG(\n,\D)
\times \GG(\n,\D)$ on $\Delta^s$ is free and its quotient, the diagonal
in ${\cal M}^s(\n,\D) \times {\cal M}^s(\n,\D)$, can be identified
naturally with ${\cal M}^s(\n,\D)$. Therefore we can identify the 
 $\G(\n,\D) \times \G(\n,\D)
\times \G(n-\n,d-\D)$-equivariant cohomology and homology of $\Delta^s
\times U(\n,\D)$ naturally with the $\G(n-\n,d-\D)$-equivariant
cohomology of ${\cal M}^s(\n,\D) \times U(\n,\D)$. The intersection of
a product $a \times b$ of cycles in $H_*(A)$ with the diagonal in $A
\times A$ becomes the intersection $a \cap b$ when the diagonal is
identified with $A$, and under Lefschetz duality the intersection
(or cap) product $\cap$ corresponds to the cup product of cohomology,
so we get
\[
\langle \left. (-a_1^{(1)})^j e_{\n,\D} \right|_{\Delta^s \times U(\n,\D)} ,(a \times
\gamma) \cap (\Delta^s \times U(\n,\D)) \rangle = 
\langle (-a_1^{(1)})^j LD(\gamma) e_{\n,\D},a \rangle
\]
for any $\gamma \in H_*^{\GG(\n,\D)} \cong H_*({\cal
M}^s(\n,\D))$. This completes the proof of Proposition
\ref{big}, given Lemmas \ref{lemone}, \ref{lemtwo}, \ref{lemthree}. 

\bigskip

{\em Proof of Lemma \ref{lemone}}: First note that the restriction of the equivariant Chern polynomial 
\[
c(-\pi_!(\hat{\cal V}^* \otimes {\cal V}))(t) = \sum_{k \geq 0}
c_k(-\pi_!(\hat{\cal V}^* \otimes {\cal V})) t^k
\]
to $\C(\n,\D)^s \times (\C(\n,\D)^{ss} \times U(\n,\D))$ equals the
product 
\[
c(-\pi_!(\hat{\cal V}^* \otimes {\cal V}_1))(t) c(-\pi_!(\hat{\cal V}^*
\otimes {\cal V}_2))(t) 
\]
where ${\cal V}_1 \oplus {\cal V}_2$ is the restriction of ${\cal V}$
to $\C(\n,\D)^{ss} \times U(\n,\D)$. Since $\D/\n > d/n$ we have
$\D/\n > (d-\D)/(n-\n)$, and so if $F_1$ and $F_2$ are semistable
bundles of ranks $\n$ and $n-\n$ and degrees $\D$ and $d-\D$
respectively then $H^0 (\Sigma,F_1^* \otimes F_2)=0$. Hence $
c(-\pi_!(\hat{\cal V}^* \otimes {\cal V}_2))(t)$ is a polynomial whose
degree is given by the Riemann-Roch theorem to be at most
\[
\n(n-\n)(g-1) + \n(d-\D) - (n-\n)\D = \n(n-\n)(g-1) +\n d - n \D.
\]
On the other hand if $F_1$ and
$F_2$ are both semistable of the same rank and degree, and one of them
is stable, then $H^0 (\Sigma,F_1^* \otimes F_2)=0$ unless $F_1 \cong
F_2$ in which case  $H^0 (\Sigma,F_1^* \otimes F_2)=\CC$
\cite[Lemma 5.1]{N2}. Thus the restriction of $c(-\pi_!(\hat{\cal V}^*
\otimes {\cal V}_1))(t)$ to the complement in
\[
\C(\n,\D)^s \times \C(\n,\D)^{ss} \times U(\n,\D)
\]
of $\Delta^s \times U(\n,\D)$ where
\begin{eqnarray}
\Delta^s & = & \{ (F_1,F_2) \in \C(\n,\D)^s \times \C(\n,\D)^{ss} : F_1
\cong F_2 \}, \nonumber\\
& = & \{(F_1,F_2) \in \C(\n,\D)^s \times \C(\n,\D)^s : F_1 \cong F_2 \}. \label{deltas} 
\end{eqnarray}
is a polynomial whose degree is given by the Riemann-Roch formula to
be at most $\n^2 (g-1)$. 

{\em Proof of Lemma \ref{lemtwo}:} By definition $\Delta^s$ is the inverse image of the diagonal
in ${\cal M}^s(\n,\D) \times {\cal M}^s(\n,\D)$ under the quotient map
\[
\C(\n,\D)^s \times \C(\n,\D)^{s} \to \frac{\C(\n,\D)^s \times \C(\n,\D)^{s}}{
\G(\n,\D) \times \G(\n,\D)} = {\cal M}^s(\n,\D) \times {\cal
M}^s(\n,\D).
\]
Therefore its codimension is $\dim {\cal M}^s(\n,\D) = \n^2(g-1)+1$ and its normal bundle in 
$\C(\n,\D)^s \times \C(\n,\D)^{s}$ is
the pullback of the normal bundle to the diagonal in ${\cal
M}^s(\n,\D)$, whose Chern polynomial is just the Chern polynomial of
the tangent bundle $T{\cal M}^s(\n,\D)$ when the diagonal is identified with ${\cal M}^s(\n,\D)$ in the obvious way.\\
\indent Now the constant central subgroup $U(1)$ of $\G(\n,\D)$ acts
trivially on $\C(\n,\D)^s$, the quotient $\GG(\n,\D) = \G(\n,\D)/U(1)$
acts freely with quotient ${\cal M}^s(\n,\D)$ and we have a canonical isomorphism \cite[p.577]{AB}
\begin{eqnarray*}
\HS_{\G(\n,\D)}(\C(\n,\D)^s) & \cong & \HS_{\GG(\n,\D)}(\C(\n,\D)^s)
\otimes \HS(BU(1))\\
& \cong & \HS({\cal M}^s(\n,\D)) \otimes \HS(BU(1)).
\end{eqnarray*}
Similarly the central subgroup  $U(1) \times U(1)$ of $\G(\n,\D)\times
\G(\n,\D)$ acts trivially on $\Delta^s$ (and on its normal bundle),
the quotient $\GG(\n,\D)\times 
\GG(\n,\D)$ acts freely with quotient ${\cal M}^s(\n,\D)$ and we have 
\begin{eqnarray}
\label{decomp}\HS_{\G(\n,\D)\times \G(\n,\D)}(\Delta^s) & \cong & 
\HS_{\GG(\n,\D)\times \GG(\n,\D)}(\Delta^s) \otimes \HS(BU(1))
\otimes \HS(BU(1))\nonumber \\
& \cong & \HS({\cal M}^s(\n,\D)) \otimes \HS(BU(1)) \otimes
\HS(BU(1)). 
\end{eqnarray}
By \cite[p.582]{AB}, in K-theory we may write 
\[
T{\cal M}^s(\n,\D) = 1 - \pi_!(\hat{\cal V}^* \otimes \hat{\cal V})
\]
if we note that $U(1)$ acts as scalar multiplication on the fibres of
$\hat{\cal V}$ and hence acts trivially on $\hat{\cal V}^* \otimes
\hat{\cal V}$, so that $\hat{\cal V}^* \otimes \hat{\cal V}$ descends
to a bundle over ${\cal M}^s(\n,\D)$. Thus the normal bundle to
$\Delta^s$ has equivariant Chern polynomial
\begin{equation}
c(1 - \pi_!(\hat{\cal V}^* \otimes \hat{\cal V}))(t) = 
c(- \pi_!(\hat{\cal V}^* \otimes \hat{\cal V}))(t), \label{ecp}
\end{equation}
which lies in the component 
\[
H^*_{\GG(\n,\D) \times \GG(\n,\D)}(\Delta^s) \cong \HS({\cal
M}^s(\n,\D))
\]
of $H^*_{\G(\n,\D) \times \G(\n,\D)}(\Delta^s)$ as decomposed at
(\ref{decomp}) above. Note that the rank of this normal bundle is
$\n^2(g-1)+1$.\\ 
\indent The only endomorphisms of a stable bundle are multiplication
by scalars \cite[Lemma 5.1]{N2}, so if $(F_1,F_2) \in \Delta^s$ then
$H^0(\Sigma, F_1^* \otimes F_2) \cong \CC$. It follows from this
that there is a line bundle ${\cal L}$ over $\Delta^s$ whose fibre at
$(F_1,F_2) \in \Delta^s$ is $H^0(\Sigma, \left. \hat{\cal V}^* \otimes
{\cal V}_1 \right|_{\{F_1\} \times \{F_2\} \times \Sigma}),$ and the
restriction of ${\cal L}$ to the diagonal in $\C(\n,\D)^s \times
\C(\n,\D)^s$ is trivial. The action of $\G(\n,\D) \times \G(\n,\D)$ on
$\Delta^s$ lifts naturally to ${\cal L}$, and if $t_1$ and $t_2$ lie
in the constant central subgroup $U(1)$ of $\G(\n,\D)$ then
$(t_1,t_2)\in \G(\n,\D) \times \G(\n,\D)$ acts on the fibres of ${\cal
L}$ as scalar multiplication by $t_1^{-1} t_2$. Thus the $\G(\n,\D)
\times \G(\n,\D)$-equivariant Chern class of ${\cal L}$ is given by 
\[
c_1({\cal L}) = a_1^{(1)} - \hat{a}_1, 
\]
where $\hat{a}_1, a_1^{(1)}$ and $a_1^{(2)}$ are the generators of the first,
second and third copies of the polynomial ring $\HS(BU(1))$ in 
\[
\HS_{\G(\n,\D) \times \G(\n,\D) \times \G(n-\n,d-\D)}(\Delta^s \times
U(\n,\D)),
\]
so that the restriction of $a_1$ equals $a_1^{(1)} + a_1^{(2)}$ (see Definition \ref{as} above).
There is a natural isomorphism 
\[
\left. \hat{\cal V} \right|_{\Delta^s \times U(\n,\D) \times \Sigma}
\otimes \pi^* {\cal L} \cong \left. {\cal V}_1 \right|_{\Delta^s
\times U(\n,\D) \times \Sigma}
\]
of bundles over $\Delta^s \times U(\n,\D) \times \Sigma$, where $\pi^*{\cal L}$ is the
pullback to $\Delta^s \times \Sigma$ of the line bundle ${\cal L}$
over $\Delta^s$. (By abuse of notation we are using $\pi$ to denote
the projection from $X \times \Sigma$ to $X$ for several different
spaces $X$.) Hence 
\[
\left. \pi_!(\hat{\cal V}^* \otimes {\cal V}_2) \right|_{\Delta^s
\times U(\n,\D)} \cong {\cal L} \otimes \left. \pi_!({\cal V}_1^* \otimes {\cal
V}_2)\right|_{\Delta^s \times U(\n,\D)},
\]
and
\[
\left. \pi_!(\hat{\cal V}^* \otimes \hat{\cal V}) \right|_{\Delta^s
\times U(\n,\D)} \cong \left. \pi_!({\cal V}_1^* \otimes {\cal
V}_1)\right|_{\Delta^s \times U(\n,\D)} \cong {\cal L}^* \otimes \left. 
\pi_!(\hat{\cal V}^* \otimes {\cal V}_1)\right|_{\Delta^s \times U(\n,\D)}.
\]
This completes the proof of Lemma \ref{lemtwo}. 

\bigskip

{\em Proof of Lemma \ref{lemthree}:} To prove Lemma \ref{lemthree} 
we use Porteous' Formula \cite[p.86]{ACGH}, \cite[$\S$14.4]{F}, 
as in Beauville's proof \cite{B} of the theorem of Atiyah and 
Bott that the classes $\ar, \bjr, \fr$ generate $\HS(\mnd)$ (cf. \cite{T2} 
and \cite{ES}). We can choose an effective divisor $D$ on $\Sigma$ 
with no multiple points such that for all semistable bundles $F_1$ 
and $F_2$ of rank $\n$ and degree $\D$ we have
\[
H^1(\Sigma,F^*_1 \otimes F_2 \otimes {\cal O}(D)) = 0.
\]
Then tensoring the exact sequence 
\[
0 \to {\cal O} \to {\cal O}(D) \to {\cal O}_D \to 0
\]
with $\hat{\cal V}^* \otimes {\cal V}_1$ over $\C(\n,\D)^{s} \times 
\C(\n,\D)^{ss} \times U(\n,\D) \times \Sigma$ gives us
\begin{equation}
\label{exseq}
\pi_!(\hat{\cal V}^* \otimes {\cal V}_1) \to  \pi_!(\hat{\cal V}^* 
\otimes {\cal V}_1 \otimes {\cal O}(D)) \to \pi_!(\hat{\cal V}^* \otimes {\cal V}_1 \otimes {\cal O}_D) 
\end{equation}
corresponding to the standard long exact sequence of cohomology. By our choice of $D$ the 
second and third terms in the sequence (\ref{exseq}) are vector bundles over 
$\C(\n,\D)^{s} \times \C(\n,\D)^{ss} \times U(\n,\D)$ which we will call ${\cal W}_1$ and ${\cal W}_2$. Then
\begin{equation}
c(-\pi_!(\hat{\cal V}^* \otimes {\cal V}_1))(t) = c({\cal W}_2 - {\cal W}_1)(t). \label{later}
\end{equation}
Moreover, by the argument above, the kernel of the map ${\cal W}_1 \to {\cal W}_2$ 
is a skyscraper sheaf supported on $\Delta^s \times U(\n,\D)$. Therefore 
$\Delta^s \times U(\n,\D)$ is the degeneracy locus of the map ${\cal W}_1 \to {\cal W}_2$, 
and by Riemann-Roch the expected codimension 
\[
\mbox{rank }{\cal W}_2 - \mbox{rank }{\cal W}_1 + 1
\]
of the degeneracy locus is equal to its actual codimension $\n^2 (g-1)+1$ (see Lemma 
\ref{lemtwo}). Therefore by Porteous' Formula the Poincar\'{e} dual of $\Delta^s \times U(\n,\D)$ in
\[
\HS_{\G(\n,\D)}(\C(\n,\D)^s) \otimes  \HS_{\G(\n,\D)}(\C(\n,\D)^{ss}) \otimes 
\HS_{\G(n-\n,d-\D)}(U(\n,\D))
\]
is
\[
c_{\n^2 (g-1) +1}({\cal W}_2 - {\cal W}_1) = c_{\n^2 (g-1) +1}(-\pi_!(\hat{\cal V}^* \otimes {\cal V}_1).
\]
In other words we have shown that the restriction of $c_{\n^2(g-1)+1}(-\pi_!(\hat{\cal V}^* \otimes {\cal V}_1))$ 
to $\C(\n,\D)^s \times \C(\n,\D)^{ss} \times U(\n,\D)$ is the equivariant cohomology class
 $[\Delta^s \times U(\n,\D)]$ defined by 
$\Delta^s \times U(\n,\D)$ in
$$ H_{\G(\n,\D) \times \G(\n,\D) \times
\G(n-\n,d-\D)}^{\n^2(g-1)+1}(\C(\n,\D)^s \times
\C(\n,\D)^{ss} \times U(\n,\D)).
$$
We can express the higher Chern classes of $-\pi_!(\hat{\cal V}^*
 \otimes {\cal V}_1)$ in a similar way using Fulton's excess Porteous
 formula \cite[p.258]{F}. We have been considering the degeneracy
 locus of the map ${\cal W}_1 \to {\cal W}_2$ (defined as the locus
 where the rank of the map is strictly less than $\mbox{rank}\, {\cal
 W}_1$). More generally for any $k \leq \min\{\mbox{rank } {\cal W}_1,
 \mbox{rank } {\cal W}_2\}\}$ we can consider the $k$th degeneracy
 locus where the rank of the map is at most $k$. But we know that if
 $F_1$ and $F_2$ are both semistable bundles of rank $\n$ and degree
 $\D$ and one of them is stable, then $H^0(\Sigma, F_1^* \otimes
 F_2)=0$ unless $F_1 \cong F_2$ in which case $H^0(\Sigma, F_1^*
 \otimes F_2) \cong \CC$. Thus the rank of the map ${\cal W}_1 \to
 {\cal W}_2$ is at least $\mbox{rank }{\cal W}_1 -1$ everywhere on
 $\C(\n,\D)^s \times \C(\n,\D)^{ss} \times U(\n,\D)$, so the $k$th
 degeneracy locus of the map ${\cal W}_1 \to {\cal W}_2$ is empty for
 any $k < \mbox{rank }{\cal W}_1 -1$. The same is of course true of
 the map ${\cal W}_1 \to {\cal W}_2 \oplus \CC^j$ given by
 composing ${\cal W}_1 \to {\cal W}_2$ with the inclusion ${\cal W}_2
 \to {\cal W}_2 \oplus \CC^j$ for any $j>0$. This means that we
 can apply the excess Porteous formula (see \cite{F}, Examples 14.4.7
 and 14.4.1) to this map ${\cal W}_1 \to {\cal W}_2 \oplus
 \CC^j$. Since the Chern polynomial of ${\cal W}_2 \oplus \CC^j$
is the same as that of ${\cal W}_2$, we obtain
\begin{equation}
\label{chern}
c_{\n^2(g-1)+j+1}({\cal W}_2 \ominus {\cal W}_1) = c_j({\cal L}^* \otimes ({\cal F} \oplus \CC^j) \ominus {\cal N})
[\Delta^s \times 
U(\n,\D)]
\end{equation}
where the restriction of ${\cal W}_1 \to {\cal W}_2$ to $\Delta^s \times U(\n,\D)$ fits into an exact sequence 
\begin{equation}
\label{es}
0 \to {\cal L} \to {\cal W}_1 \to {\cal W}_2 \to {\cal F} \to 0 
\end{equation}
and ${\cal N}$ is the normal bundle to $\Delta^s \times U(\n,\D)$ in $\C(\n,\D)^s \times \C(\n,\D)^{ss} 
\times U(\n,\D)$. Here the kernel ${\cal L}$ is of course the line bundle over $\Delta^s$ 
which first appeared in Lemma \ref{lemtwo}. From the exact sequence (\ref{es}) over 
$\Delta^s \times U(\n,\D)$ we obtain
\begin{eqnarray*}
c({\cal F})(t) & = & c({\cal W}_2 \ominus {\cal W}_1 \oplus {\cal L})(t)\\
& = & c(-\pi_!(\hat{{\cal V}}^* \otimes {\cal V}_1) \oplus {\cal L})(t). 
\end{eqnarray*}
Moreover by Lemma \ref{lemtwo} the normal bundle ${\cal N}$ to 
$\Delta^s \times U(\n,\D)$ has equivariant Chern polynomial
\[
c(-{\cal L}^* \otimes \pi_!(\hat{\cal V}^* \otimes {\cal V}_1))(t).
\]
Therefore
$ c({\cal L}^* \otimes ({\cal F} \oplus \CC^j) \ominus {\cal N})(t)$ is equal to
$$ c({\cal L}^* \otimes((-\pi_!(\hat{\cal V}^* \otimes {\cal V}_1))\oplus{\cal
L} \oplus \CC^j) \oplus {\cal L}^* \otimes \pi_!(\hat{\cal V}^* \otimes
{\cal V}_1))(t) $$
$$ =  c({\cal L}^* \otimes {\cal L} \oplus {\cal L}^* \otimes
\CC^j)(t)
 =  (c({\cal L}^*)(t))^j = (1 - c_1 ({\cal L})t)^j = (1+(\hat{a}_1 - a_1^{(1)})t)^j.
$$
Thus by (\ref{later}) and (\ref{chern})
\begin{eqnarray*}
(\hat{a}_1- a_1^{(1)})^j [\Delta^s \times U(\n,\D)] & = & c_{\n^2 (g-1) + j +1}({\cal W}_2 \ominus {\cal W}_1)\\
& = & c_{\n^2 (g-1) + j +1}(-\pi_!(\hat{\cal V}^* \otimes {\cal
V}_1)).
\end{eqnarray*}
This completes the proof of the first part of Lemma \ref{lemthree}.\\
\indent By Lemma \ref{lemone} the restriction  to $\C(\n,\D)^s \times \C(\n,\D)^{ss}
\times U(\n,\D)$ of $c(-\pi_!(\hat{\cal
V}^* \otimes {\cal V}))(t)$ equals
\[
c(-\pi_!(\hat{\cal V}^* \otimes {\cal V}_1))(t) c(-\pi_!(\hat{\cal
V}^* \otimes {\cal V}_2))(t)
\]
where $c(-\pi_!(\hat{\cal V}^* \otimes {\cal V}_2))(t)$ is a
polynomial of degree at most 
\[
\n(n-\n)(g-1) + \n d - n \D.
\]
If $r \geq r(\n,\D)$ then
\[
r - \n(n-\n)(g-1) + \n d - n \D \geq \n^2 (g-1) + 1
\]
and so the restriction of $c_r(-\pi_!(\hat{\cal V}^* \otimes {\cal
V}))$ is the coefficient of $t^r$ in
\[
\sum_{k \geq \n^2 (g-1) + 1} c_k(-\pi_!(\hat{\cal
V}^* \otimes {\cal V}_1)) t^k c(-\pi_!(\hat{\cal
V}^* \otimes {\cal V}_2))(t),
\]
which we have just seen is equal to 
\[
\sum_{j \geq 0} (\hat{a}_1 - a_1^{(1)})^j [\Delta^s \times U(\n,\D)]
t^{\n^2 (g-1) + j + 1} c(-\pi_!(\hat{\cal
V}^* \otimes {\cal V}_2))(t)
\]
\[
=
[\Delta^s \times U(\n,\D)] \left(\frac{t^{\n^2 (g-1)
+ 1}}{1 + (a_1^{(1)} - \hat{a}_1)t} c(-\pi_!(\hat{\cal
V}^* \otimes {\cal V}_2))(t) \right).
\] 
By Lemma \ref{lemtwo} we have
\[
 \left. \pi_!(\hat{\cal V}^* \otimes {\cal V}_2) \right|_{\Delta^s
 \times U(\n,\D)} \cong \left. {\cal L} \otimes \pi_!({\cal V}_1^* \otimes
 {\cal V}_2)  \right|_{\Delta^s \times U(\n,\D)}.
\]
If $F$ is a vector bundle of rank $k$ and $L$ is a line bundle then
the Chern polynomials $c(F)(t)$ and $c(F \otimes L)(t)$ are related by
the formula
\[
t^k c(F \otimes L)(t^{-1}) = (t+c_1(L))^k c(F) ((t+c_1(L))^{-1}) 
\]
\cite[21.10]{BT}. The same is therefore true if $F$ is a virtual vector bundle
of virtual rank $k$. Thus
\[
t^{\n(n-\n)(g-1)+\n d - n \D}c(- \pi_!(\left. \hat{\cal V}^* \otimes {\cal
V}_2)\right|_{\Delta^s \times U(\n,\D)})(t^{-1}) 
\]
equals
\[ 
(t+c_1({\cal L}))^{\n(n-\n)(g-1)+\n d - n \D}c(-\pi_!(\left. {\cal V}_1^* \otimes {\cal
V}_2)\right|_{\Delta^s \times U(\n,\D)})((t+c_1({\cal L}))^{-1}).
\]
Hence the restriction of $c_r(-\pi_!(\hat{\cal V}^* \otimes {\cal
V}))$ is the coefficient of $t^{-r}$ in the product of $
[\Delta^s \times U(\n,\D)]$ with 
$$ \frac{t^{-\n^2 (g-1)
- 1}}{1 + (a_1^{(1)} - \hat{a}_1)/t} c \left( -\pi_! \left(  {\cal
V}^* \otimes {\cal V}_2   
\right) \right)(t^{-1}) 
$$
$$=
 \frac{t^{-\n^2 (g-1)} (t + a_1^{(1)} - \hat{a}_1)^{\n(n-\n)(g-1)
+ \n d - n \D}      }{(t + a_1^{(1)} - \hat{a}_1)   t^{\n(n-\n)(g-1)
+ \n d - n \D}       }  c\left( -\pi_! \left(  {\cal
V}_1^* \otimes {\cal V}_2   
\right) \right) \left( \frac{1}{t + a_1^{(1)} - \hat{a}_1} \right)
$$
$$ =
 t^{-r(\n,\D)+1} (t +    a_1^{(1)} - \hat{a}_1 )^{\n(n-\n)(g-1) + \n d - n \D-1} c\left( -\pi_! \left(  {\cal
V}_1^* \otimes {\cal V}_2 
\right)\right)\left( \frac{1}{t + a_1^{(1)} - \hat{a}_1} \right).
$$
This coefficient is the coefficient of $t^{-r + r(\n,\D) -1}$ in 
\[
[\Delta^s \times U(\n,\D)]  (t +  a_1^{(1)} - \hat{a}_1
)^{\n(n-\n)(g-1) + \n d - n \D-1} c\left( -\pi_! \left(  {\cal
V}_1^* \otimes {\cal V}_2  
\right)\right)\left( \frac{1}{t + a_1^{(1)} - \hat{a}_1} \right).
\]
But $r \geq r(\n,\D)$, so $-r + r(\n,\D) -1$ is strictly negative, and
$c(-\pi_!({\cal V}_1^* \otimes {\cal V}_2 
))(t)$ is a polynomial of degree at most $\n(n-\n)(g-1) + \n d -
n \D$ in t, so that
\[
(t+a_1^{(1)} - \hat{a}_1)^{\n(n-\n)(g-1) + \n d - n \D -1} c\left( -\pi_!  \left(
{\cal V}_1^* \otimes {\cal V}_2 
\right) \right) \left( \frac{1}{t + a_1^{(1)} - \hat{a}_1} \right)
\]
is the sum of 
\[
\frac{ c_{\n(n-\n)(g-1) + \n d-n\D-1} \left(
-\pi_!\left( {\cal V}_1^* \otimes {\cal V}_2 
 \right) \right)}{t + a_1^{(1)} - \hat{a}_1} 
=
\frac{  \left. e_{\n,\D} \right|_{\Delta^s \times U(\n,\D)} }{t+a_1^{(1)} - \hat{a}_1}
\]
and a polynomial in $t$. Thus the coefficient we want is equal to the
coefficient of $t^{-r+r(\n,\D)-1}$ in
\[
[\Delta^s \times U(\n,\D)] \left( \frac{ e_{\n,\D}
}{t + a_1^{(1)} - \hat{a}_1} \right) =
e_{\n,\D} [\Delta^s \times U(\n,\D)] \left( \sum_{j \geq 0}
\frac{(-(a_1^{(1)} - \hat{a}_1))^j}{t^{j+1}} \right).
\]
This completes the proof of Lemma \ref{lemthree} and hence of
Proposition \ref{big}.

\section{Some refinements}

\renorm

Recall that given $\eta \in H^*_{\G(\n,\D)}(\C(\n,\D)^{ss}) \otimes
H^*_{\G(n-\n,d-\D)}(\C(n-\n,d-\D))$ we wish to find a relation whose
restriction to $\C(n_1,d_1)^{ss} \times U(n_1,d_1)$ is zero when
$d_1/n_1 < \D/\n$ or $d_1/n_1 = \D/\n$ and $n_1< n$, and when
$(n_1,d_1) = (\n,\D)$ equals $\eta e_{\n,\D}$. Lemma \ref{seven} and Proposition \ref{big} deal with $\eta$ in the image
of the Lefschetz duality map $LD$ which maps
$H_*^{\GG(\n,\D)}(\C(\n,\D)^s)$ to $H^*_{\GG(\n,\D)}(\C(\n,\D)^{ss})$ and
thus into 
\[
H^*_{\G(\n,\D)}(\C(\n,\D)^{ss}) = H^*_{\GG(\n,\D)}(\C(\n,\D)^{ss}) \otimes
\HS(BU(1)),
\]
and more generally with $\eta$ of the form
\begin{equation}
\label{extra}
\eta = (-a_1^{(1)})^j LD(\gamma),
\end{equation}
for some $\gamma \in H_*^{\GG(\n,\D)}(\C(\n,\D)^s)$ where $a_1^{(1)}$
is the generator of $\HS(BU(1))$. When $\n$ and $\D$ are coprime this
gives us all $\eta \in \HS_{\G(\n,\D)}(\C(\n,\D)^{ss})$.\\
\indent We know that $\HS_{\G(n,d)}(\C(n,d))$ is generated (as an
algebra over $\Q$) by the K\"{u}nneth components $a_r, b_r^j, f_r$ for
$1 \leq r \leq n$ and $1 \leq j \leq 2g$ of the equivariant Chern
classes
\[
c_r({\cal V}) = a_r \otimes 1 + \sum_{j=1}^{2g} b_r^j \otimes \alpha_j
+ f_r \otimes \omega \in  \HS_{\G(n,d)}(\C(n,d)) \otimes \HS(\Sigma)
\]
where $\{1\}, \{\alpha_1, \ldots , \alpha_{2g}\}, \{\omega\}$ are
standard bases for $H^0(\Sigma), H^1(\Sigma)$ and $H^2(\Sigma)$
respectively. The equivariant Chern classes of ${\cal V}$ are the
elementary symmetric polynomials in its (formal) equivariant Chern
roots $A_1,...,A_n$. All symmetric polynomials in the Chern roots
$A_1,...A_n$ represent elements of $\HS_{\G(n,d)}(\C(n,d)) \times
\HS(\Sigma)$, and in particular we have 
\[
\sum_{j=1}^n (A_j)^r \in \HS_{\G(n,d)}(\C(n,d)) \otimes \HS(\Sigma)
\]
for $r \geq 0$. Moreover all symmetric polynomials in $A_1,...,A_n$
can be expressed as polynomials in
\[
\{\sum_{j=1}^n (A_j)^r: 1 \leq r \leq n \} ,
\] 
(this follows by induction from Girard's formula;
see \cite[p 195]{MS})
and in particular $c_1({\cal V}),...,c_{n}({\cal V})$ can be expressed
thus.  Therefore the generators $a_r, b_r^j, f_r$ can all be expressed
as polynomials in the K\"{u}nneth components of $\{\sum_{j=1}^n
(A_j)^r: 1 \leq r \leq n \}$ and thus these K\"{u}nneth components
generate $\HS_{\G(n,d)}(\C(n,d))$ as a $\Q$-algebra. Similarly 
\[
H^*_{\G(\n,\D)}(\C(\n,\D)) \otimes H^*_{\G(n-\n,d-\D)}(\C(n-\n,d-\D))
\]
is generated as a $\Q$-algebra by the K\"{u}nneth components of
\[
\sum_{j=1}^{\n} (B_j)^r
\]
for $1 \leq r \leq \n$ and of 
\[
\sum_{j=\n+1}^{n} (B_j)^r
\]
for $1 \leq r \leq n-\n$, where $B_1,...,B_{\n}$ are the
$\G(\n,\D)$-equivariant Chern roots of the universal bundle ${\cal
V}_1$ over $\C(\n,\D)$ and  $B_{\n+1},...,B_n$ are the
$\G(n-\n,d-\D)$-equivariant Chern roots of the universal bundle ${\cal
V}_2$ over $\C(n-\n,d-\D)$. 
 We may assume that $B_j$ is the restriction of $A_j$ to
$\C(\n,\D) \times \C(n-\n,d-\D)$ for $1 \leq j \leq n$. Then
\[
\left. \sum_{j=1}^{\n} (B_j)^r + \sum_{j=\n+1}^{n} (B_j)^r = \sum_{j=1}^n
(A_j)^r \right|_{\C(\n,\D) \times \C(n-\n,d-\D)}
\]
for $1 \leq r \leq n$. This proves:
\begin{lem}
\label{twelve}
\[
\HS_{\G(\n,\D)}(\C(\n,\D)) \otimes \HS_{\G(n-\n,d-\D)}(\C(n-\n,d-\D))
\]
is generated as a module over $\HG(\C)$ by
$\HS_{\G(\n,\D)} (\C(\n,\D))$. 
\end{lem}

When $\n$ and $\D$ are coprime, we have now obtained the relations we
need:

\begin{prop}
\label{fourteen}
Let $\n$ and $\D$ be positive integers with $0<\n<n$ and
$\D/\n > d/n,$ and suppose 
\[
\eta \in H^*_{\G(\n,\D)}(\C(\n,\D)^{ss}) \otimes
H^*_{\G(n-\n,d-\D)}(U(\n,\D)) \cong \HG(S_{\n,\D}).
\]
If $\n$ and $\D$ are coprime there is a relation in the ideal in
$H^*_{\G}(\C)$ generated by the slant products
\[
\{ c_r(-\pi_!(\hat{\cal V}^* \otimes {\cal V})) \backslash \gamma: 
r \geq r(\n,\D) = n \n (g-1) + \n d - n \D + 1, \gamma \in
H_*^{\GG(\n,\D)}(\C(\n,\D)^{ss}) \}
\]
whose image under the restriction map
\[
H^*_{\G(n,d)}(\C(n,d)) \to H^*_{\G(n_1,d_1)}(\C(n_1,d_1)) \otimes
H^*_{\G(n-n_1,d-d_1)}(U(n_1,d_1))
\]
is zero when $d_1/n_1 < \D/\n$ and also when $d_1/n_1 = \D/\n$ and
$n_1<\n$, and when $(n_1,d_1) = (\n,\D)$ equals $\eta e_{\n,\D}$.
\end{prop}
{\em Proof:} This follows from Remark \ref{six}, Lemma \ref{seven},
Proposition \ref{big} and Lemma
\ref{twelve}, together with the fact (see Remark \ref{three} and
Proposition \ref{KCC}) that the restriction maps
\[
H^*_{\G(\n,\D)}(\C(\n,\D)) \to 
H^*_{\G(\n,\D)}(\C(\n,\D)^{ss})
\]
and 
\[
H^*_{\G(n-\n,d-\D)}(\C(n-\n,d-\D)) \to 
H^*_{\G(n-\n,d-\D)}(U(\n,\D))
\]
are surjective.

\section{Parabolic bundles}

\renorm

At this point we need to observe that all our arguments can be generalised from vector bundles
over a compact Riemann surface $\Sigma$ to apply to moduli spaces of parabolic vector bundles over 
$\Sigma$ equipped with finitely many marked points $p_1, \ldots, p_k$. For
more details on parabolic bundles see \cite{MeS,Ni,Se,Se2}. 

To simplify
the notation in this discussion we shall assume that there is only one marked point $p$.
A parabolic bundle on $(\Sigma,p)$ is a holomorphic bundle $E$ on $\Sigma$ with a flag
\begin{equation} \label{flag}
0 = E_p^{(m+1)} \subset E_p^{(m)} \subset \ldots \subset E_p^{(1)} = E_p
\end{equation}
of length $m \geq 1$ in the fibre $E_p$ of $E$ at $p$, together with a sequence of weights
$0\leq \alpha_1 < \alpha_2 < \ldots < \alpha_m < 1$. 

A holomorphic subbundle $F$ of the parabolic bundle $E$ inherits a parabolic structure
by intersecting the flag (\ref{flag}) with the fibre $F_p$ of $F$ at $p$ and taking the
corresponding weights, but omitting any intersections $F_p \cap E_p^{(k)}$ such
that $F_p \cap E_p^{(k)} = F_p \cap E_p^{(k-1)}$. A quotient bundle of $E$ becomes a
parabolic bundle in a similar way.

If $1 \leq k \leq m$ then the {\it multiplicity} of the
weight $\alpha_k$ is
$$j_k = \dim(E_p^{(k)}/E_p^{(k+1)}),$$
the {\it parabolic degree} of $E$ is
$$ \mbox{pardeg}(E) = \deg(E) + \sum_{k=1}^m \alpha_k \dim(E_p^{(k)}/E_p^{(k+1)})$$
and the {\it parabolic slope} of $E$ is
$$\mbox{par}\mu(E) = \frac{\mbox{pardeg}(E)}{\mbox{rank}(E)}.$$
Stability and semistability for parabolic bundles are then defined just as for
vector bundles, and each parabolic bundle has a canonical Harder-Narasimhan
filtration by parabolic subbundles.

\begin{rem} \label{Parss=Pars} It will be important for us to have sufficient conditions on 
$n,d$, $(\alpha_1,\ldots,\alpha_m)$ and $(j_1,\ldots,j_m)$ for all semistable
parabolic bundles with rank $n$, degree $d$, parabolic weights $\alpha_1, \ldots \alpha_m$
and multiplicities $j_1, \ldots, j_m$ to be stable. Such a semistable parabolic bundle
$E$ is not stable if and only if it has a proper subbundle $\hat{E}$ such that
$$ \frac{\mbox{pardeg}(\hat{E})}{\mbox{rank}(\hat{E})} =  \frac{\mbox{pardeg}(E)}{\mbox{rank}(E)},$$
and so there exist integers $\n, \D$ and $\hat{j}_1, \ldots, \hat{j}_m$ such that
$0< \n <n$ and $0 \leq \hat{j}_k \leq j_k$ and $\hat{j}_1 + \cdots + \hat{j}_m = \n$ and
$$ \frac{\D + \sum_{k=1}^m \alpha_k \hat{j}_k}{\n} = \frac{ d + \sum_{k=1}^m \alpha_k j_k}{n}.$$
It is easy to see that the set of parabolic weights $(\alpha_1,\ldots,\alpha_m)$ for which there exist
integers $\n, \D$ and $\hat{j}_1, \ldots, \hat{j}_m$ such that
$0< \n <n$ and $0 \leq \hat{j}_k \leq j_k$ and $\hat{j}_1 + \cdots + \hat{j}_m = \n$ and
$$ \frac{\D + \sum_{k=1}^m \alpha_k \hat{j}_k}{\n} = \frac{ d + \sum_{k=1}^m \alpha_k j_k}{n}$$
is closed in $\{(\alpha_1, \ldots, \alpha_m) \in 
\RR^m: 0 \leq \alpha_1 < \ldots < \alpha_m <1 \}$, since there are only finitely 
many choices of $\n$ and $\hat{j}_1, \ldots, \hat{j}_m$
and then
$$\hat{d} = \frac{\n d}{n} + \sum_{k=1}^m \alpha_k \left( \frac{\n j_k}{n} - \hat{j}_k \right)$$
is required to be an integer. This set is also nowhere dense in $\{(\alpha_1, \ldots, \alpha_m) \in 
\RR^m: 0 \leq \alpha_1 < \ldots < \alpha_m <1 \}$ {\em unless} 
\begin{enumerate}
\item[(i)] $d, n$ and $j_1, \ldots, j_m$
have a common factor $h>1$, 
\item[(ii)] $\n$ is an integer multiple of $n/h$ and
\item[(iii)] $\D = \n d/n$ and $\hat{j}_k = \n j_k/n$ for $1 \leq k \leq m$, 
\end{enumerate}
in which case there are no
constraints on $\alpha_1, \ldots, \alpha_m$ other than $0 \leq \alpha_1 < \ldots < \alpha_m <1$.
In particular this set is nowhere dense in $\{(\alpha_1, \ldots, \alpha_m) \in 
\RR^m: 0 \leq \alpha_1 < \ldots < \alpha_m <1 \}$ if $j_k = 1$ for some $k$.
\end{rem}

\newcommand{\PG}{{{Par}{\cal G}}}
\newcommand{\PGc}{{{Par}{\cal G}_c}}

Nitsure \cite{Ni} showed that it is possible to generalise to parabolic bundles the methods of Atiyah 
and Bott \cite{AB} for calculating the Betti numbers
of the moduli spaces $\mnd$,
and finding generators for their cohomology rings, 
provided that $n$, $d$, $(\alpha_1,\ldots,\alpha_m)$ and $(j_1, \ldots, j_m)$ are such that 
semistability coincides with stability. For this we fix a $C^\infty$ bundle ${\cal E}$
over $\Sigma$ with rank $n$ and degree $d$, and a filtration of the fibre of ${\cal E}$ at $p$ as
at (\ref{flag}) with associated weights $ 0 \leq \alpha_1 < \ldots < \alpha_m <1 $. As before we
let $\C$ be the space of holomorphic structures on ${\cal E}$, and we let $\PG$ 
(respectively $\PGc$) be the subgroup 
of the gauge group $\G$ (respectively the 
complexified gauge group $\G_c$) consisting of all complex $C^\infty$ automorphisms of 
${\cal E}$ which preserve the parabolic filtration. Then $\PGc$ acts on $\C$ and its orbits correspond 
to the isomorphism classes of parabolic bundles over $\Sigma$ of the given rank, degree and
parabolic structure. 

\newcommand{\Cparss}{{\cal C}^{parss}}

The parabolic Harder--Narasimhan filtration  defines for us a smooth $\PG$-equivariantly perfect
stratification of $\C$ whose unique open stratum $\Cparss$ consists of those
induced parabolic structures on ${\cal E}$ which are semistable. 
However the strata are in general not connected, and the connected
components may even have different dimensions. So Nitsure refines the stratification
by using the {\em compound type} $(\mu,I)$ of a parabolic bundle $E$. Here $\mu$ is the
parabolic Harder--Narasimhan type of $E$ given by
$$\mu= ( \mbox{par}\mu(E_1/E_0), \ldots, \mbox{par}\mu(E_s/E_{s-1}))$$
where $0=E_0 \subset E_1 \subset \ldots \subset E_s=E$ is the parabolic Harder--Narasimhan
filtration of $E$ and the entry $\mbox{par}\mu(E_j/E_{j-1})$ is repeated $\dim(E_j/E_{j-1})$ times,
while $I$ is the $s \times m$ intersection matrix $(I_{k,\ell})$ of $E$ with entries defined by
$$I_{1,m} = \dim((E_1)_p \cap E_p^{(m)})$$
and
$$ I_{k,\ell} = \dim ((E_k)_p \cap E_p^{(\ell)}) - \sum_{\begin{array}{c} 1\leq i\leq k \\ \ell \leq j \leq m \\ (i,j) \neq (k,\ell)
\end{array}} I_{i,j}$$
for $1\leq k \leq s$ and $1 \leq \ell \leq m$, where $(E_k)_p$ is the fibre at the marked point
$p$ of the subbundle $E_k$ of $E$, and $E^{(\ell)}_p$ is the $\ell$th term in the fixed flag
(\ref{flag}) in $E_p$. If $\C_{\mu,I}$ denotes the subset of $\C$ corresponding to all holomorphic structures
on ${\cal E}$ with the same compound type $(\mu,I)$ as $E$, then $\C_{\mu,I}$ is connected and
\begin{equation} \label{bod} H^*_\PG (\C_{\mu,I}) \cong \bigotimes_{i=1}^s 
H^*_{\PG({\cal E}_i/{\cal E}_{i-1})} ( \C({\cal E}_i/{\cal E}_{i-1})^{parss}) \end{equation}
where $0= {\cal E}_0  \subset {\cal E}_1 \subset \ldots \subset {\cal E}_s = {\cal E}$ is a fixed $C^\infty$
filtration of ${\cal E}$ of type $(\mu,I)$, while $\C({\cal E}_i/{\cal E}_{i-1})^{parss}$ and 
$\PG({\cal E}_i/{\cal E}_{i-1})$
are defined just as $\Cparss$ and $\PG$ except that ${\cal E}$ is replaced by ${\cal E}_i/{\cal E}_{i-1}$ with its
induced parabolic structure. Moreover by \cite{Ni} Proposition 3.6 the $\PG$-equivariant Euler class $e_{\mu,I}$ of the
normal bundle to $\C_{\mu,I}$ in $\C$ is not a zero divisor in $H^*_{\PG}(\C_{\mu,I})$, and hence by the criterion of
Atiyah and Bott \cite{AB} $\S$13 the stratification of $\C$ with strata $\C_{\mu,I}$ indexed by compound type
$(\mu,I)$ is $\PG$-equivariantly perfect (see \cite{Ni} $\S$3). Thus the restriction map
$$H^*_\PG(\C) \to H^*_\PG(\Cparss)$$
is surjective, and by \cite{Ni} Proposition 3.2 we have
$$H^*_\PG(\C) \cong H^*(B\G) \otimes H^*({\cal F})$$
where ${\cal F}$ is the flag variety of flags in $\CC^n$ of the type specified by the parabolic structure.
Furthermore, when semistability coincides with stability we have
$$H^*_\PG(\Cparss) \cong H^*({\cal M}) \otimes H^*(BU(1))$$
where ${\cal M}= \Cparss /\PG$ is the moduli space of parabolic bundles over $\Sigma$ of the given
parabolic data (that is, the given rank, degree, weights and multiplicities).

Thus generators of $H^*(B\G)$ and $H^*({\cal F})$ provide generators for the rational
cohomology of the moduli space ${\cal M}$, and the completeness criteria described in
Proposition 3.3 generalise immediately to give completeness criteria for relations between
these generators in terms of the strata $\C_{\mu,I}$.

Just as in the non-parabolic case, these completeness criteria can be simplified by using a 
coarser stratification of $\C$ indexed by the {\em coarse compound parabolic type} of
parabolic bundles. Here a parabolic bundle $E$ has coarse compound parabolic type
$(n_1,d_1,j_{1,1},\ldots,j_{m,1})$ if its parabolic Harder--Narasimhan filtration
$$0=E_0 \subset E_1 \subset \ldots \subset E_s  = E$$
is such that $E_1$ has rank $n_1$ and degree $d_1$ and
$$j_{\ell,1} = \dim \left( \frac{(E_1)_p \cap E_p^{(\ell)}}{(E_1)_p \cap E_p^{(\ell +1)}} \right)$$
for $1 \leq \ell \leq m$, where $(E_1)_p$ is the fibre of $E_1$ at $p$ and $E_p^{(\ell)}$ is
the $\ell$th term in the fixed flag (\ref{flag}) in the fibre $E_p$ of $E$ at the marked point
$p$.
We will denote by $S_{n_1,d_1,j_{1,1},\ldots,j_{m,1}}$ the subset of $\C$ consisting of all
parabolic structures with coarse compound parabolic type $(n_1,d_1,j_{1,1},\ldots,j_{m,1})$. 

The proof that the stratification of $\C$ by coarse type is smooth and $\G_c$-equivariantly
perfect (se $\S$3) carries over directly to give a proof that the stratification of $\C$ by coarse compound
parabolic type is smooth and $\PG$-equivariantly perfect, and we obtain

\begin{prop}

\label{MPCC}
(Modified Parabolic Completeness Criteria) Let ${\cal R}$ be a
subset of the kernel of the restriction map
\[
H^*_\PG(\C) \rightarrow H^*_\PG(\Cparss).
\]
Suppose that for each sequence of integers $(\n,\D,\hat{j}_1,\ldots, \hat{j}_m)$ with $0<\n<n$ and
$0 \leq \hat{j}_\ell \leq j_\ell$ for $1 \leq \ell \leq m$ and $\hat{j}_1 + \cdots + \hat{j}_m = \n$ and
$$\frac{\D + \sum_{\ell=1}^m \alpha_\ell \hat{j}_\ell }{\n} > \frac{d + \sum_{\ell = 1}^m \alpha_\ell j_\ell}{n},$$ 
there is a subset ${\cal
R}_{\n,\D,\hat{j}_1,\ldots, \hat{j}_m}$ of the ideal generated by ${\cal R}$ in $H^*_\PG (\C)$ such that
the image of ${\cal R}_{\n,\D,\hat{j}_1,\ldots, \hat{j}_m}$ under the restriction map
\[
H^*_{\PG} (\C) \rightarrow H^*_{\PG} (S_{n_1,d_1,j_{1,1},\ldots,j_{m,1}})
\]
is zero when either
$$\frac{d_1+ \sum_{\ell =1}^m \alpha_\ell j_{\ell,1}}{n_1}<\frac{\D + \sum_{\ell=1}^m 
\alpha_\ell \hat{j}_\ell }{\n}$$
or equality holds and $n_1 \leq \n$, 
except for the case when 
$$(n_1,d_1,j_{1,1},\ldots,j_{m,1}) = (\n,\D,\hat{j}_1,\ldots, \hat{j}_m),$$
and in the latter case 
the image of ${\cal R}_{\n,\D,\hat{j}_1,\ldots, \hat{j}_m}$ 
equals the ideal of
$H^*_\PG(S_{\n,\D},\hat{j}_1,\ldots, \hat{j}_m)$ generated by the equivariant Euler class $e_{\n,\D,
\hat{j}_1,\ldots, \hat{j}_m}$
of the normal bundle to the stratum $S_{\n,\D},\hat{j}_1,\ldots, \hat{j}_m$ in
$\C.$ Then ${\cal R}$ generates the kernel of the restriction map
\[
H^*_\PG(\C) \rightarrow H^*_\PG(\Cparss)
\]
as an ideal in $H^*_\PG(\C).$
\end{prop}

Furthermore, in the notation used at (\ref{bod}), we have
\begin{equation} \label{delta} H^*_\PG(S_{\n,\D,\hat{j}_1,\ldots, \hat{j}_m}) \cong 
H^*_{\PG({\cal E}_1)} (\C({\cal E}_1)^{parss}) \otimes H^*_{\PG({\cal E}/{\cal E}_1)} 
(U(\n,\D,\hat{j}_1,\ldots, \hat{j}_m)) \end{equation}
where $U(\n,\D,\hat{j}_1,\ldots, \hat{j}_m)$ is an open subset of $\C({\cal E}/{\cal E}_1)$
such that the restriction map from $H^*_{\PG({\cal E}/{\cal E}_1)}(\C({\cal E}/{\cal E}_1))$
to $H^*_{\PG({\cal E}/{\cal E}_1)} 
(U(\n,\D,\hat{j}_1,\ldots, \hat{j}_m))$ is surjective (cf. (3.6)).

Suppose now that $(\n,\D,\hat{j}_1,\ldots, \hat{j}_m))$ is a sequence of integers
satisfying the conditions of Proposition \ref{MPCC}. Since 
$$\frac{\D + \sum_{\ell=1}^m \alpha_\ell \hat{j}_\ell }{\n} > \frac{d + \sum_{\ell = 1}^m \alpha_\ell j_\ell}{n},$$
there are no nonzero parabolic bundle maps from a semistable parabolic bundle of rank $\n$
and degree $\D$ with parabolic weights $(\alpha_1,\ldots,\alpha_m)$ having multiplicities
$(\hat{j}_1, \ldots, \hat{j}_m)$ to a semistable parabolic bundle $E$ of rank $n$ and degree $d$
with parabolic weights $(\alpha_1,\ldots,\alpha_m)$ having multiplicities
$({j}_1, \ldots, {j}_m)$. Therefore if $ParHom(\hat{\cal V},{\cal V})$ denotes the subsheaf
of $\hat{\cal V}^* \otimes {\cal V}$ over $\C(n,d) \times \C(\n,\D) \times \Sigma$ consisting of
parabolic bundle maps, then the restriction to 
$\C(n,d)^{parss} \times \C(\n,\D)^{parss}$
of $- \pi_!(ParHom(\hat{{\cal V}},{\cal V})$ is a bundle (rather than just a virtual bundle) 
whose rank is given by the Riemann--Roch
theorem as follows:

\begin{lem} \label{rankParH}
The restriction to 
$$\C(n,d)^{parss} \times \C(\n,\D)^{parss}$$ 
of the virtual bundle 
$- \pi_!(ParHom(\hat{{\cal V}},{\cal V})$ is a bundle of rank
$$n \n (g-1) + \D n - d \n + \sum_{\ell=2}^m (j_1 + \cdots + j_{\ell -1})\hat{j}_\ell.$$
\end{lem}
\noindent {\em Proof} (cf. \cite{Ni} Proposition 1.17): There is an exact sequence of
sheaves
$$0 \to ParHom(\hat{{\cal V}},{\cal V}) \to \hat{{\cal V}}^* \otimes {\cal V} \to {\cal S} \to 0$$
over $\C(n,d) \times \C(\n,\D) \times \Sigma$, where ${\cal S}$ is supported on
$\C(n,d) \times \C(\n,\D) \times \{ p\}$ and is the pullback of the skyscraper sheaf on $\Sigma$
supported at $p$ with fibre 
$$\frac{\hat{{\cal V}}^*_p \otimes {\cal V}_p}{ParHom (\hat{{\cal V}}_p, {\cal V}_p)}$$
which has dimension 
$$\sum_{\ell=2}^m (j_1 + \cdots + j_{\ell - 1}) \hat{j}_\ell = \sum_{\ell = 1}^{m-1} j_\ell (\hat{j}_{\ell+1}
+ \cdots + \hat{j}_m).$$
The Riemann--Roch theorem tells us that the virtual rank of $-\pi_!(\hat{{\cal V}}^* \otimes
{\cal V})$ is $n\n(g-1) + \D n - d \n$ (cf. (\ref{codim*}) or \cite{AB} $\S$7) and the result follows.

\bigskip

\newcommand{\PhG}{{Par}\hat{\G}}

It follows from this lemma that the equivariant Chern classes
$$c_r(- \pi_!(ParHom(\hat{{\cal V}},{\cal V})) \in H^*_\PG(\C(n,d)) \otimes H^*_{{Par}\hat{\G}}
(\C(\n,\D))$$
restrict to zero in $\C(\n,\D)^{parss} \times \C(n,d)^{parss}$ when
$$r > r(\n,\D,\hat{j}_1,\ldots,\hat{j}_m) = n \n (g-1) + \D n - d \n + \sum_{\ell =2}^m (j_1 + \cdots + j_{\ell - 1}) \hat{j}_\ell,$$
and then their slant products with any $\gamma \in H_*^{\PhG} (\C(\n,\D)^{parss})$ give
us elements of the kernel of the restriction map
$$H^*_\PG(\C(n,d)) \to H^*_\PG(\C(n,d)^{parss}) \cong H^*({\cal M}) \otimes H^*(BU(1)).$$
The arguments of $\S$4 and $\S$5 now generalise directly to the case of parabolic bundles, and
the proof of Proposition 5.2 gives us

\begin{prop} \label{Par5.2}

Let $\n, \D$ and $\hat{j}_1, \ldots, \hat{j}_m$ be integers satisfying $0<\n<n$ and
$0\leq \hat{j}_\ell \leq j_\ell$ and $\hat{j}_1 + \cdots + \hat{j}_m = \n$ and
$$\frac{\D + \sum_{\ell = 1}^m \alpha_\ell \hat{j}_\ell}{\n} > \frac{d + \sum_{\ell=1}^m \alpha_\ell j_\ell}{n},$$
 and suppose that $\eta $
is any element of 
$$H^*_\PG(S_{\n,\D},\hat{j}_1,\ldots, \hat{j}_m) \cong 
H^*_{\PG({\cal E}_1)} (\C({\cal E}_1)^{parss}) \otimes H^*_{\PG({\cal E}/{\cal E}_1)} 
(U(\n,\D,\hat{j}_1,\ldots, \hat{j}_m)).$$
If every semistable parabolic bundle of rank $\n$, degree $\D$ and parabolic weights
$(\alpha_1, \ldots, \alpha_m)$ with multiplicities $(\hat{j}_1, \ldots, \hat{j}_m)$ is stable,
then there is a relation in the ideal in
$H^*_{\PG}(\C)$ generated by the slant products
\[
\{ c_r(-\pi_!(ParHom(\hat{\cal V}, {\cal V}))) \backslash \gamma: r > r(\n,\D,\hat{j}_1,\ldots,\hat{j}_m), 
\gamma \in
H_*^{{{{Par}\hat{\G}}}(\n,\D)}(\C(\n,\D)^{parss}) \}
\]
whose image in
$$H^*_\PG(S_{\n,\D},\hat{j}_1,\ldots, \hat{j}_m) \cong 
H^*_{\PG({\cal E}_1)} (\C({\cal E}_1)^{parss}) \otimes H^*_{\PG({\cal E}/{\cal E}_1)} 
(U(\n,\D,\hat{j}_1,\ldots, \hat{j}_m))$$
is zero when 
$$\frac{d_1+ \sum_{\ell =1}^m \alpha_\ell j_{\ell,1}}{n_1}<\frac{\D + \sum_{\ell=1}^m 
\alpha_\ell \hat{j}_\ell }{\n}$$
and also when equality holds and $n_1 \leq \n$, 
except for the case 
$(n_1,d_1,j_{1,1},\ldots,j_{m,1}) = (\n,\D,\hat{j}_1,\ldots, \hat{j}_m)$ 
when it equals $\eta e_{\n,\D}$.
\end{prop}

\section{The good case}

\renorm

We are aiming to prove the following generalisation of Theorem \ref{thm} which is stated
using the notation introduced in $\S$6.

\begin{thm} \label{Par2.1} If $n$ and $d$ are coprime 
and $\C(n,d)^{parss} = \C(n,d)^{pars}$
then the kernel of the restriction
map 
\[
\rho: H^*_{{Par}\G}(\C(n,d)) \to H^*_{\PG}(\C(n,d)^{parss})
\]
is generated as an ideal in $H^*_{\PG}(\C(n,d))$ by slant
products of the form 
\[
c_{r}(-\pi_!(ParHom(\hat{\cal V}, {\cal V}))) \backslash \gamma
\] 
for integers $\n,\D,\hat{j}_1,\ldots,\hat{j}_m$ and $r$ satisfying $\hat{j}_1 + \cdots + \hat{j}_m = \n$ 
with $0<\n<n$ and 
$0\leq \hat{j}_\ell \leq j_\ell$ for $1 \leq \ell \leq m$, and
\begin{equation} \label{(a)} 
\frac{d + \sum_{\ell=1}^m \alpha_\ell j_\ell}{n} <
\frac{\D + \sum_{\ell = 1}^m \alpha_\ell \hat{j}_\ell}{\n} < \frac{d + \sum_{\ell=1}^m \alpha_\ell j_\ell}{n}+1
\end{equation}
and
\begin{equation} \label{(b)}
n \n (g-1)  < r - 
\D n + d \n - \sum_{\ell =2}^m (j_1 + \cdots + j_{\ell - 1}) \hat{j}_\ell <
n \n (g+1), 
\end{equation}
where
$$c_r(- \pi_!(ParHom(\hat{{\cal V}},{\cal V})) \in H^*_\PG(\C(n,d)) \otimes H^*_{{Par}\hat{\G}}
(\C(\n,\D))$$
is the $r$th equivariant Chern class of the virtual bundle $- \pi_!(ParHom(\hat{{\cal V}},{\cal V}))$
and 
$\gamma$ lies in the image of the natural map
$$ H_*^{{{{Par}\hat{\G}}}(\n,\D)}(\C(\n,\D)^{parss})  \to
H_*^{{{{Par}\hat{\G}}}(\n,\D)}(\C(\n,\D)).$$
\end{thm}

The slightly weaker version of Theorem \ref{Par2.1} omitting the second inequalities
in (\ref{(a)}) and in (\ref{(b)}) follows immediately from Propositions \ref{MPCC} and
\ref{Par5.2} in the special case when the choices of $d$, $n$ and the parabolic data
$(\alpha_1,\ldots,\alpha_m)$ and $(j_1, \ldots, j_m)$ satisfy the following condition:

\begin{df} \label{defn7.2}
We shall call $d$, $n$ and the parabolic data
$(\alpha_1,\ldots,\alpha_m)$ and $(j_1, \ldots, j_m)$ a {\em good choice of parabolic data}
if every semistable parabolic bundle of rank $\n$, degree $\D$, parabolic weights
$(\alpha_1,\ldots,\alpha_m)$ and multiplicities $(\hat{j}_1, \ldots, \hat{j}_m)$
satisfying (\ref{(a)}) and $\hat{j}_1 + \cdots + \hat{j}_m = \n$, with $0<\n<n$ and 
$0\leq \hat{j}_\ell \leq j_\ell$ for $1 \leq \ell \leq m$, is stable.
\end{df}

\begin{rem} \label{crux}
In fact if we look carefully at this proof of the slightly weaker version of Theorem \ref{Par2.1} for good
parabolic data (in particular at the proofs of Propositions \ref{MCC} and  4.5), we find
that if $\eta$ lies in the kernel of the restriction map $\rho$ then we can write
$$\eta = \sigma_1 + \cdots + \sigma_k$$
where
$$\sigma_i = c_{r_i}(- \pi_!(ParHom(\hat{{\cal V}}^{(i)},{\cal V})) \backslash \gamma_i$$
for some $(\n_i, \D_i, \hat{j}_{1,i}, \ldots, \hat{j}_{m,i}, r_i)$ satisfying  $0<\n_i<n$ and 
$0\leq \hat{j}_{\ell,i} \leq j_\ell$ for $1 \leq \ell \leq m$ and $\hat{j}_{1,i} + \cdots + \hat{j}_{m,i} = \n_i$ and
$$\frac{d + \sum_{\ell=1}^m \alpha_\ell j_\ell}{n} <
\frac{\D_1 + \sum_{\ell = 1}^m \alpha_\ell \hat{j}_{\ell,1}}{\n_1} \leq \cdots
\leq \frac{\D_k + \sum_{\ell=1}^m \alpha_\ell \hat{j}_{\ell,k}}{\n_k}
$$
and
$$ r_i >
n \n_i (g-1) + \D_i n - d \n_i + \sum_{\ell =2}^m (j_{i} + \cdots + j_{\ell - 1}) \hat{j}_{\ell,i}, 
$$
where $\hat{{\cal V}}^{(i)}$ is a universal bundle while 
$\gamma_i$ is given by the image of a suitable restriction of $\eta$ under
the inverse of an appropriate Lefschetz duality map, which is an isomorphism since
semistability equals stability (cf. Remark 4.1).
\end{rem}

\begin{rem} \label{rem7.3}
If $L$ is any line bundle over $\Sigma$, then a parabolic bundle $E$ is semistable
(respectively stable) if and only if $E \otimes L$ is semistable (respectively stable).
Thus the condition that every semistable parabolic bundle of rank $\n$, degree
 $\D$, parabolic weights
$(\alpha_1,\ldots,\alpha_m)$ and multiplicities $(\hat{j}_1, \ldots, \hat{j}_m)$
is stable depends on $\D$ only through its remainder modulo $\n$. Thus for
given $n$ and $d$ there are effectively only finitely many choices of
$\n$, $\D$ and $(\hat{j}_1, \ldots, \hat{j}_m)$ to consider in Definition 
\ref{defn7.2}.
\end{rem}

\begin{rem} \label{rem7.4}
It follows from Remarks \ref{Parss=Pars} and \ref{rem7.3} that if $j_k=1$ for $1 \leq k \leq m=n$ then
an arbitrarily small perturbation can be made to any parabolic weights $(\alpha_1,\ldots,
\alpha_m)$ to give us good parabolic data in the sense of Definition \ref{defn7.2}.
\end{rem}

\begin{rem} \label{rem7.5} Suppose that every semistable parabolic
bundle of rank $n$, degree
 $d$, parabolic weights
$(\alpha_1,\ldots,\alpha_m)$ and multiplicities $({j}_1, \ldots, {j}_m)$ is stable.
Then it follows from Remark \ref{Parss=Pars} that if the parabolic weights 
$(\alpha_1, \ldots, \alpha_m)$ are perturbed very slightly then the semistability 
condition is unaffected (that is, $\C(n,d)^{parss}$ remains the same)
and in addition the conditions on the integers 
$\n$, $\D$, $(\hat{j}_1, \ldots, \hat{j}_m)$ and $r$ in (\ref{(a)}) and (\ref{(b)}) of
Theorem \ref{Par2.1} are unaltered. Thus if an arbitrarily small perturbation
of the parabolic weights can be chosen to give us good parabolic data in the
sense of Definition \ref{defn7.2}, then the weaker version of Theorem \ref{Par2.1}
omitting the second inequalities in (\ref{(a)}) and (\ref{(b)}) must hold. In 
particular Remark \ref{rem7.4} tells us that if $d$ and $n$ are coprime and
$j_\ell = 1$ for $1 \leq \ell \leq m=n$ then this weaker version of Theorem
\ref{Par2.1} holds.
\end{rem}

\section{Reduction to the good case}

\renorm

We have now proved the following weaker version of Theorem \ref{Par2.1} in
the special case when $m=n$ and $j_1=\ldots = j_m=1$ (see Remark \ref{rem7.5}).
Our aim in this section is to deduce the general case.

\begin{thm} \label{wPar2.1} If $n$ and $d$ are coprime 
and $\C(n,d)^{parss}=\C(n,d)^{pars}$ then the kernel of the restriction
map 
\[
\rho: H^*_{{Par}\G}(\C(n,d)) \to H^*_{\PG}(\C(n,d)^{parss})
\]
is generated as an ideal in $H^*_{\PG}(\C(n,d))$ by slant
products of the form 
\[
c_{r}(-\pi_!(ParHom(\hat{\cal V}, {\cal V}))) \backslash \gamma
\] 
for integers $\n,\D,\hat{j}_1,\ldots,\hat{j}_m$ and $r$ satisfying $0<\n<n$ and 
$0\leq \hat{j}_\ell \leq j_\ell$ for $1 \leq \ell \leq m$ and $\hat{j}_1 + \cdots + \hat{j}_m = \n$ and
\begin{equation} \label{(aa)} 
\frac{d + \sum_{\ell=1}^m \alpha_\ell j_\ell}{n} <
\frac{\D + \sum_{\ell = 1}^m \alpha_\ell \hat{j}_\ell}{\n} 
\end{equation}
and
\begin{equation} \label{(bb)}
n \n (g-1) + \D n - d \n + \sum_{\ell =2}^m (j_1 + \cdots + j_{\ell - 1}) \hat{j}_\ell < r 
\end{equation}
where
$$c_r(- \pi_!(ParHom(\hat{{\cal V}},{\cal V})) \in H^*_\PG(\C(n,d)) \otimes H^*_{{Par}\hat{\G}}
(\C(\n,\D))$$
is the $r$th equivariant Chern class of the equivariant virtual bundle $- \pi_!(ParHom(\hat{{\cal V}},{\cal V})$
and 
$\gamma$ lies in the image of the natural map
$$ H_*^{{{{Par}\hat{\G}}}(\n,\D)}(\C(\n,\D)^{parss})  \to
H_*^{{{{Par}\hat{\G}}}(\n,\D)}(\C(\n,\D)).$$
\end{thm}

To deduce this theorem from 
the special case when $m=n$ and $j_1=\ldots = j_m=1$,
we exploit the relationship between the kernel of $\rho$ and the kernel of
the corresponding restriction map 
$$
\rho^T: H^*_{{Par}_T\G}(\C(n,d)) \to H^*_{{Par}_T\G}(\C(n,d)^{{par}_T ss})
$$
where $\rho^T$, ${Par}_T\G$ and $\C(n,d)^{{par}_T ss}$ are defined
as $\rho$, $\PG$ and $\C(n,d)^{{Par} ss}$ except that $m$ is replaced by $n$,
the multiplicities $j_\ell$ all become $1$ and the weights $(\alpha_1, \ldots, \alpha_m)$ are
replaced by $(\alpha_1^T, \ldots, \alpha_n^T)$ with the following properties:
\begin{itemize}
\item[(i)]
$0 \leq \alpha_1^T < \alpha_2^T < \ldots < \alpha^T_n < 1$,
\item[(ii)] if $1 \leq i \leq m$ then
$\quad \displaystyle{j_i \alpha_i =  \sum_{j_1 + \cdots j_{i-1} < j \leq j_1 + \cdots + j_i} \alpha_j^T} \quad $
and 
\item[(iii)] $\alpha_j^T$ is a small
perturbation of $\alpha_i$ if $j_1 + \cdots j_{i-1} < j \leq j_1 + \cdots + j_i$.
\end{itemize}

Let $B$ be the standard Borel subgroup of $GL(n;\CC)$ and let $P$ be the parabolic
subgroup which preserves the flag
$$\CC^{j_1} \subset \CC^{j_1 + j_2} \subset \ldots \subset \CC^{j_1 + \cdots + j_{m-1}} \subset \CC^n.$$
Let ${\cal F} = GL(n;\CC)/P \cong U(n)/P\cap U(n)$ and ${\cal F}_T = GL(n;\CC)/B \cong U(n)/T$
where $T$ is the standard maximal torus of $U(n)$. Note that if 
$$\phi_p : \G \to U(n)$$
is the homomorphism which associates to any $g \in \G$ its action on the fibre of
${\cal E}$ at the marked point $p$ of $\Sigma$ then
$${Par}_T \G = \phi_p^{-1}(T)$$
and 
$${Par}\G = \phi_p^{-1}(P \cap U(n)).$$
We have 
\begin{equation} \label{decomp1}
H^*_\PG(\C) \cong H^*(B\G) \otimes H^*({\cal F})
\end{equation}
while
\begin{equation} \label{decomp2}
H^*_{{Par}_T\G} (\C) \cong H^*(B\G) \otimes H^*({\cal F}_T)
\end{equation}
and
\begin{equation} \label{decomp3}
H^*_\PG (\C) \cong [H^*_{{Par}_T\G}(\C)]^W
\end{equation}
where $W$ is the Weyl group of $P \cap U(n)$ (cf. \cite{Ni} Proposition 3.2 and 
\cite{Borel} Theorem 20.6).
By \cite{Knew} Lemma 7.10 and Remarks 7.11 and 7.12 
(see also Remark \ref{martin} below), if 
$\eta \in H^*_{\PG}(\C) \cong [H^*_{{Par}_T\G}(\C)]^W$
then
\begin{equation} \label{mart} \eta \in \ker \rho \mbox{ if and only if } {\cal D}\eta \in \ker \rho^T
\end{equation}
and, furthermore, multiplication by ${\cal D}$ defines an isomorphism from $\ker \rho$
to $\ker \rho^T \cap [H^*_{{Par}_T\G}(\C)]^{{anti}W}$, with inverse given by
$$\zeta \mapsto \frac{1}{|W| {\cal D}} \sum_{w \in W} (-1)^w w(\zeta).$$
Here ${\cal D}$ is the image in $
H^*_{{Par}_T\G} (\C) $ of the element of $H^*_T(\C) \cong H^*(BT)$ represented by the
product of the positive roots of $U(n)$ under the map
$$H^*_T(\C) \to 
H^*_{{Par}_T\G} (\C) $$
induced by the homomorphism $\phi_p: {Par}_T\G \to T$. Under the decomposition
(\ref{decomp2}) ${\cal D}$ corresponds to the fundamental class in $H^*({\cal F}_T)$.
Note that ${\cal D}$ is anti-invariant for the action of the Weyl group $W$ in the sense
that $w({\cal D}) = (-1)^w {\cal D}$ for all $w \in W$, where $(-1)^w$ is the determinant
of $w$ regarded as an automorphism of the Lie algebra of $T$. 
We denote by $ [H^*_{{Par}_T\G}(\C)]^{{anti}W}$
the set of elements of $ H^*_{{Par}_T\G}(\C)$ which are anti-invariant for the action of $W$; if
$ \eta \in H^*_{{Par}_T\G}(\C)$ then $\eta \in  H^*_{{Par}_T\G}(\C)]^{{anti}W}$
if and only if 
$$\eta = \frac{1}{|W|} \sum_{w \in W} (-1)^w w(\eta).$$

\bigskip

\newcommand{\PTG}{{{Par}_T\G}}
\newcommand{\PTGhat}{{{Par}_{\hat{T}}\hat{\G}}}

\noindent {\em Proof of Theorem \ref{wPar2.1}:} Suppose that $\eta \in \ker \rho$. Then by
(\ref{mart}) we have ${\cal D} \eta \in \ker \rho^T$. But by Remark \ref{rem7.5} we know that
Theorem \ref{wPar2.1} is true when $m=n$ and $j_1=\ldots = j_n = 1$. Therefore ${\cal D} \eta$
lies in the ideal in $H^*_{{Par}_T\G} (\C) $ generated by slant products of the form
$$
c_{r}(-\pi_!(ParHom_{\hat{T},T}(\hat{\cal V}, {\cal V}))) \backslash \gamma
$$
for integers $\n,\D$ and $r$ and a subset $J$ of $\{1,\ldots,n\}$ with $\hat{n}$ elements,
satisfying $0<\n<n$ and 
\begin{equation} \label{(aaa)} 
\frac{d + \sum_{\ell=1}^m \alpha_\ell j_\ell}{n} <
\frac{\D + \sum_{\ell \in J} \alpha_\ell }{\n} 
\end{equation}
and
\begin{equation} \label{(bbb)}
n \n (g-1) + \D n - d \n + \sum_{\ell =2}^n \left|J \cap \{\ell, \ldots , n\}\right| < r . 
\end{equation}
Here $\hat{T}$ is the standard maximal torus of $U(\n)$, while 
$$c_r(- \pi_!(ParHom_{\hat{T},T} (\hat{{\cal V}},{\cal V})) \in H^*_\PTG (\C(n,d)) \otimes H^*_{\PTGhat}
(\C(\n,\D))$$
is the $r$th equivariant Chern class of the virtual bundle $- \pi_!(ParHom_{\hat{T},T}
(\hat{{\cal V}},{\cal V})$, which is defined as in $\S$6 with $\PG$ and ${Par}\hat{\G}$ 
replaced by $\PTG$ and $\PTGhat$,
and 
$\gamma$ lies in the image of the natural map
\begin{equation} \label{bbb} H_*^{{{{Par}_{\hat{T}}\hat{\G}}}(\n,\D)}(\C(\n,\D)^{parss})  \to
H_*^{{{{Par}_{\hat{T}}\hat{\G}}}(\n,\D)}(\C(\n,\D)). \end{equation}
Moreover by Remark \ref{crux} we can take $\gamma$ to be the image of the restriction
of ${\cal D} \eta$  under the inverse of a suitable Lefshetz duality map.

There is an exact sequence of sheaves
$$0 \to Par_{\hat{T},T}Hom(\hat{\cal V}, {\cal V}) \to ParHom(\hat{\cal V},{\cal V}) \to \hat{\cal S} \to 0$$
over $\C(n,d) \times \C(\n,\D) \times \Sigma$ where $\hat{\cal S}$ is supported on 
$\C(n,d) \times \C(\n,\D) \times \{p\}$ and is the pullback of the skyscraper sheaf on $\Sigma$
supported at $p$ with fibre
$$\frac{ParHom(\hat{\cal V}_p,{\cal V}_p)}{ParHom_{\hat{T},T}
(\hat{\cal V}_p,{\cal V}_p)}.$$ 
For integers $\n,\D$ and $r$ and a subset $J$ of $\{1,\ldots,n\}$ with $\hat{n}$ elements
as above, define integers $\hat{j}_1, \ldots ,\hat{j}_m$ by
$$\hat{j}_k = \left| J \cap \{ j_1 + \cdots + j_{k-1}+1, \ldots, j_1+ \cdots + j_k\} \right|$$
for $1 \leq k \leq m$. Then $0 \leq \hat{j}_k \leq j_k$ and since $j_1 \cdots + j_m =n$ we
have
$$\hat{j}_1 + \cdots + \hat{j}_m = | J| = \n.$$
The weights of the action of $\hat{T} \times T$ on $ParHom(\hat{\cal V}_p,{\cal V}_p)$ are $\lambda_\ell
- \hat{\lambda}_i$ for $1 \leq i \leq \n$ and $1 \leq \ell \leq n$ such that
$h \leq k$
where $h$ and $k$ satisfy
$$\hat{j}_1+\cdots + \hat{j}_{h-1} < i \leq \hat{j}_1 + \cdots + \hat{j}_h$$
and
$$j_1+\cdots + j_{k-1} < \ell \leq j_1 + \cdots + j_k.$$
The weights of the action of $\hat{T} \times T$ on $ParHom_{\hat{T},T}
(\hat{\cal V}_p,{\cal V}_p)$ are $\lambda_\ell
- \hat{\lambda}_i$ for $1 \leq i \leq \n$ and $1 \leq \ell \leq n$ such that
$$i \leq \left| J \cap \{1,\ldots,\ell\}\right| .$$
Since $i > |J \cap \{1,\ldots,\ell\}|$ and $j_1+\cdots + j_{k-1} < \ell $
implies that $i > \hat{j}_1 + \cdots + \hat{j}_k$, it follows that
the weights of the action of $\hat{T} \times T$ on 
$$\frac{ParHom(\hat{\cal V}_p,{\cal V}_p)}{ParHom_{\hat{T},T}
(\hat{\cal V}_p,{\cal V}_p)}$$ 
are $\lambda_\ell
- \hat{\lambda}_i$ for $1 \leq i \leq \n$ and $1 \leq \ell \leq n$ such that
$$i > \left| J \cap \{1,\ldots,\ell\}\right| $$
and 
$$\hat{j}_1+\cdots + \hat{j}_{k-1} < i \leq \hat{j}_1 + \cdots + \hat{j}_k$$
and
$$j_1+\cdots + j_{k-1} < \ell \leq j_1 + \cdots + j_k$$
for some $k \in \{1, \ldots, m\}$. 
Thus the equivariant Chern polynomial of $\pi_!\hat{\cal S}$ is
$$c(\pi_!\hat{\cal S})(t) = \prod_{k=1}^m \,\,\, \prod_{\ell = j_1+\cdots + j_{k-1}+1}^{j_1+\cdots + j_k}
\,\,\,  \prod_{i= | J\cap\{1,\ldots,\ell\}|+1}^{\hat{j}_1 + \cdots + \hat{j}_k}
\left( 1 + (\lambda_\ell - \hat{\lambda}_i)t \right).$$
Therefore the equivariant Chern polynomials of 
$-\pi_!(ParHom(\hat{\cal V},{\cal V}))$ and $-\pi_!(ParHom_{\hat{T},T}
(\hat{\cal V},{\cal V}))$ are related by
$$c(-\pi_!(ParHom_{\hat{T},T}
(\hat{\cal V},{\cal V})))(t) = c(-\pi_!(ParHom(\hat{\cal V},{\cal V})))(t) \,\,\,c(\pi_!\hat{\cal S})(t)$$
where $c(\pi_!\hat{\cal S})(t)$ is a polynomial of degree
$$\sum_{k=1}^m \,\,\, \sum_{\ell = j_1+\cdots + j_{k-1}+1}^{j_1+\cdots + j_k}
\hat{j}_1 + \cdots + \hat{j}_k -  \left| J\cap\{1,\ldots,\ell\} \right| $$
$$ = \sum_{k=1}^m \,\,\, \sum_{\ell = j_1+\cdots + j_{k-1}+1}^{j_1+\cdots + j_k}
\left| J \cap \{\ell+1, \ldots, n\}\right| - (\hat{j}_{k+1} + \cdots + \hat{j}_m)$$
$$= \sum_{\ell = 1}^{n-1} \left| J \cap \{\ell+1, \ldots, n\}\right|  - \sum_{k=1}^{m-1} j_k (\hat{j}_{k+1}
+ \cdots + \hat{j}_m)$$
$$= \sum_{\ell = 2}^{n} \left| J \cap \{\ell, \ldots, n\}\right| - \sum_{k=2}^{m} (j_1 + \cdots j_{k-1})\hat{j}_{k}
.$$
Hence if the image of $c(-\pi_!(ParHom(\hat{\cal V},{\cal V})))(t)$
under some homomorphism is a polynomial in $t$ of degree at most
$$
n \n (g-1) + \D n - d \n + \sum_{k =2}^{m} (j_1 + \cdots + j_{k - 1}) \hat{j}_k
$$
then the image of $c(-\pi_!(ParHom_{\hat{T},T}
(\hat{\cal V},{\cal V})))(t)$ under the same homomorphism is a polynomial in $t$ of degree
at most
$$n \n (g-1) + \D n - d \n + \sum_{\ell =2}^n \left|J \cap \{\ell, \ldots , n\}\right|. $$
Thus the ideal in $H^*_\PTG(\C(n,d))$ generated by all elements of the form
$$
c_{r}(-\pi_!(ParHom_{\hat{T},T}(\hat{\cal V}, {\cal V}))) \backslash \gamma,
$$
for $r > 
n \n (g-1) + \D n - d \n + \sum_{\ell =2}^n |J \cap \{\ell, \ldots , n\}| $
and 
$\gamma$ in the image of the natural map from
$ H_*^{{{{Par}_{\hat{T}}\hat{\G}}}(\n,\D)}(\C(\n,\D)^{parss}) $ to $
H_*^{{{{Par}_{\hat{T}}\hat{\G}}}(\n,\D)}(\C(\n,\D))$, is contained in the 
ideal generated by all elements of the form
\[
c_{r}(-\pi_!(ParHom(\hat{\cal V}, {\cal V}))) \backslash \gamma
\] 
for $r > 
n \n (g-1) + \D n - d \n + \sum_{\ell =2}^m (j_1 + \cdots + j_{\ell - 1}) \hat{j}_\ell$
and 
$\gamma$ in the image of the natural map  from
$ H_*^{{{{Par}_{\hat{T}}\hat{\G}}}(\n,\D)}(\C(\n,\D)^{parss}) $ to $
H_*^{{{{Par}_{\hat{T}}\hat{\G}}}(\n,\D)}(\C(\n,\D))$, and so
${\cal D} \eta$ lies in this ideal. Recall from (\ref{decomp3}) that
$$
H^*_\PG (\C(n,d)) \cong [H^*_{{Par}_T\G}(\C(n,d))]^W
$$
and analogously we have 
$$
H^*_{{Par}\hat{\G}} (\C(\n,\D)) \cong [H^*_{\PTGhat}(\C(\n,\D))]^{\hat{W}}.
$$
By using the dual version of (\ref{mart}) for equivariant
homology, together with the observation after (\ref{bbb}), based on Remark \ref{crux}, 
that we can take $\gamma$ to be the image under a suitable duality map of a restriction
of ${\cal D} \eta$, and by replacing $\gamma$ by
$$\frac{1}{|\hat{W}|} \sum_{\hat{w} \in \hat{W}}  \hat{w} \gamma$$
if necessary, we find that 
we can assume that $\gamma$ lies in  the image of the natural map  from
$ H_*^{{{{Par}\hat{\G}}}(\n,\D)}(\C(\n,\D)^{parss}) $ to $
H_*^{{{{Par}\hat{\G}}}(\n,\D)}(\C(\n,\D)) 
\cong [H^*_{\PTGhat}(\C(\n,\D))]^{\hat{W}}$.
Thus ${\cal D} \eta$ lies in the ideal in 
$H^*_\PTG(\C(n,d))$  generated by all elements of the form
\[
c_{r}(-\pi_!(ParHom(\hat{\cal V}, {\cal V}))) \backslash \gamma
\] 
for $r > 
n \n (g-1) + \D n - d \n + \sum_{\ell =2}^m (j_1 + \cdots + j_{\ell - 1}) \hat{j}_\ell$
and 
$\gamma$ in the image of the natural map  from
$ H_*^{{{{Par}\hat{\G}}}(\n,\D)}(\C(\n,\D)^{parss}) $ to $
H_*^{{{{Par}\hat{\G}}}(\n,\D)}(\C(\n,\D))$.

Finally we consider the action of $W$ on 
$$
H^*_{{Par}_T\G} (\C) \cong H^*(B\G) \otimes H^*({\cal F}_T)
$$
(cf. (\ref{decomp2})). Since multiplication by ${\cal D}$ defines a bijection from
the set of $W$-invariant elements of $
H^*_{{Par}_T\G} (\C)$ to the set of elements which are anti-invariant for
the action of $W$ (see for example \cite{ES} Lemma 1.2), by writing
$${\cal D}\eta = \frac{1}{|{W}|} \sum_{{w} \in {W}} (-1)^{{w}} {w}({\cal D} \eta)$$
we deduce that $\eta$ lies in the ideal in 
$$
H^*_\PG (\C(n,d)) \cong [H^*_{{Par}_T\G}(\C(n,d))]^W
$$
generated by elements of the form
\[
c_{r}(-\pi_!(ParHom(\hat{\cal V}, {\cal V}))) \backslash \gamma
\] 
for $r > 
n \n (g-1) + \D n - d \n + \sum_{\ell =2}^m (j_1 + \cdots + j_{\ell - 1}) \hat{j}_\ell$ 
and 
$\gamma$ in the image of the natural map from 
$ H_*^{{{{Par}\hat{\G}}}(\n,\D)}(\C(\n,\D)^{parss}) $ to $
H_*^{{{{Par}\hat{\G}}}(\n,\D)}(\C(\n,\D))$, as required.

\section{Further refinements}

\renorm

As a special case of Theorem \ref{wPar2.1} we have

\begin{thm}
\label{a8} If $n$ and $d$ are coprime then 
the kernel of the restriction map 
\[
H^*_{\G(n,d)}(\C(n,d)) \to H^*_{\G(n,d)}(\C(n,d)^{ss})
\]
is generated as an ideal in $H^*_{\G(n,d)}(\C(n,d))$ by slant
products of the form 
\[
c_{r}(-\pi_!(\hat{\cal V}^* \otimes {\cal V})) \backslash \gamma
\] 
for integers $\n,\D$ and $r$ satisfying $0<\n<n$ and $\frac{\D}{\n} > \frac{d}{n}$ and 
\[
r > n\n(g-1)-d\n+\D n 
\]
where
\[
c_{r}(-\pi_!(\hat{\cal V}^* \otimes {\cal V})) \in
H^*_{\G(\n,\D)}(\C(\n,\D)) \otimes H^*_{\G(n,d)}(\C(n,d))
\]
is the $r$th equivariant Chern class of the virtual bundle
$-\pi_!(\hat{\cal V}^* \otimes {\cal V})$ over $\C(\n,\D) \times
\C(n,d)$, and $\gamma$ lies in the image of the natural map
\[
H_*^{\G(\n,\D)}(\C(\n,\D)^{ss}) \to H_*^{\G(\n,\D)}(\C(\n,\D)).
\]
\end{thm}

In order to complete the proof of Theorem \ref{thm} it remains only to restrict the range of 
$\D$ and $r$ in the statement of Theorem \ref{a8} to
\[
\frac{d}{n} < \frac{\D}{\n} < \frac{d}{n} + 1 \mbox{ and } \n n (g-1) - d \n + \D n < r < \n n (g+1) - d \n + \D n.
\]
For this recall that for any $\delta \in {\bf Z}$ there is an  isomorphism from ${\cal M}(\n,\D)$ to ${\cal M}(\n,\D 
+ \n \delta)$ given by tensoring bundles of rank $\n$ and $\D$ by a fixed line bundle $L$ over $\Sigma$ 
of degree $\delta$. Recall that $\C(\n,\D)$ is the space of all holomorphic structures on a fixed $C^\infty$ 
bundle ${\cal E}$ of rank $\n$ and degree $\D$ on $\Sigma$, so we may take $\C(\n,\D+\n \delta)$ to be 
the space of all holomorphic structures on $L \otimes {\cal E}$. Then the isomorphism from ${\cal M}(\n,\D)$
to ${\cal M}(\n,\D+ \n \delta)$ and its inverse come from maps between $\C(\n,\D)$ and $ \C(\n,\D+\n \delta)$ given by 
tensoring with $L$ and $L^{-1}$. This also gives us an isomorphism
\[
H^*_{\G(\n,\D+\n \delta)}(\C(\n,\D+\n \delta)^{ss}) \to H^*_{\G(\n,\D)}(\C(\n,\D)^{ss}).
\]
Under the map $\C(\n,\D) \to \C(\n,\D+\n \delta)$, the universal bundle on 
$$\C(\n,\D+\n \delta) \times \Sigma$$ 
pulls back to the tensor product $\hat{\cal V} \otimes L$ of the universal 
bundle $\hat{\cal V}$ over $\C(\n,\D) \times \Sigma$ and the line bundle 
$L$ over $\Sigma$, identified with its pullback to $\C(\n,\D) \times \Sigma$. 
Since $c_1(L) = \delta \omega$ where $\omega$ is the usual generator of 
$H^2(\Sigma)$, the Grothendieck-Riemann-Roch Theorem tells us that the 
equivariant Chern character of $-\pi!((\hat{\cal V} \otimes L)^* \otimes {\cal V})$ is 
\begin{eqnarray*}
\ch(-\pi_!((\hat{\cal V} \otimes L)^* \otimes {\cal V})) & = & \pi_* (((g-1)\omega-1) 
\ch (L^* \otimes \hat{\cal V}^* \otimes {\cal V}))\\
& = & \pi_* (((g-1)\omega-1) e^{\delta \omega} \ch (\hat{\cal V}^* \otimes {\cal V}))\\
& = & \ch (-\pi_!(\hat{\cal V}^* \otimes {\cal V})) + \left. \delta \ch \left(\hat{\cal V}^* \otimes 
{\cal V} \right|_{\C(\n,\D) \times \C(n,d) \times \{p\}} \right)
\end{eqnarray*}
for any $p \in \Sigma$, so its equivariant Chern polynomial satisfies
\[
c(-\pi_!((\hat{\cal V} \otimes L)^* \otimes {\cal V}))(t) = c(-\pi_!(\hat{\cal V}^* \otimes 
{\cal V}))(t) \left. (c\left(\hat{\cal V}^* \otimes {\cal V} \right|_{\C(\n,\D) \times \C(n,d) \times \{p\}} \right)(t))^\delta.
\]
Now 
\[
\left. c\left(\hat{\cal V}^* \otimes {\cal V} \right|_{\C(\n,\D) \times \C(n,d) \times \{p\}} \right)(t)
\]
is the Chern polynomial of a vector bundle of rank $n\n$ over $\C(\n,\D) \times \C(n,d)$, 
so it is a polynomial of degree at most $\n n$ in $t$ with coefficients in
\[
\HS_{\G(\n,\D)}(\C(\n,\D)) \otimes \HS_{\G(n,d)}(\C(n,d)).
\]
Moreover the difference between 
\[
n \n (g-1) - \n d + (\D + \n \delta)n \mbox{ and } n \n (g-1) - d \n + \D n
\]
is $\delta n \n$. It follows that if $\delta \geq 0$ and $r > n \n (g-1) - \n d + (\D + \n \delta)n$, 
then the $r$th equivariant Chern class of the virtual bundle 
$-\pi_!((\hat{\cal V} \otimes L)^* \otimes {\cal V})$ belongs to the ideal in
\[
\HS_{\G(\n,\D)}(\C(\n,\D)) \otimes \HS_{\G(n,d)}(\C(n,d)).
\]
generated by the equivariant Chern classes $c_{r'}(-\pi_!(\hat{\cal V}^* \otimes {\cal V}))$ such that
\[
r' > n \n (g-1) - d \n + \D n.
\]
This means we can replace the condition that $\D/\n > d/n$ in the statement of Theorem \ref{a8} by
\[
\frac{d}{n} < \frac{\D}{\n} \leq \frac{d}{n}+1.
\]
Since $d$ and $n$ are coprime and $0 < \n < n$, this is equivalent to the condition
\[
\frac{d}{n} < \frac{\D}{\n} < \frac{d}{n}+1.
\]

Next we use the Grothendieck-Riemann-Roch Theorem to obtain explicit formulas 
for the equivariant  Chern polynomial $c(-\pi_!(\hat{\cal V}^* \otimes {\cal V}))(t)$ in 
terms of the Atiyah-Bott generators $a_r,b_r^j,f_r$ of $H^*_{\G(n,d)}(\C(n,d))$ and
corresponding generators of $H^*_{\G(\n,\D)}(\C(\n,\D))$.

\begin{prop}
\label{r10}
The equivariant Chern polynomial $c(-\pi_!(\hat{\cal V}^* \otimes {\cal V}))(t)$ is equal to 
\begin{equation}
\label{EC}
\Omega(t)^{g-1} \prod_{k=1}^{n} \prod_{l=1}^{\n} (1+(\delta_k-\hat{\delta}_l)t)^{-W_{k,l}}
\exp \left\{ \frac{-\Xi_(k,l)t}{1+(\delta_k-\hat{\delta}_l)t} \right\}
\end{equation}
where $\delta_1,\ldots,\delta_n$ and $\hat{\delta}_1,...,\hat{\delta}_{\n}$ are formal degree two classes 
such that the $r$th elementary symmetric polynomial in $\delta_1,\ldots,\delta_n$ is $a_r$ and the $r$th 
elementary symmetric polynomial in $\hat{\delta}_1,...,\hat{\delta}_{\n}$  is $\hat{a}_r$. Here
\[
\Omega(t) = \prod_{k=1}^n \prod_{l=1}^{\n} (1+(\delta_k-\hat{\delta}_l)t)
\]
is the equivariant Chern polynomial of the restriction of $\hat{\cal V}^* \otimes {\cal V}$ restricted 
to $\C(n,d) \times \{p\}$ for any $p \in \Sigma$, and hence can be expressed in terms of the equivariant 
Chern polynomials
\[
\left. c({\cal V} \right|_{\C(n,d) \times \{p\}})(t) = 1 + a_1 t + a_2 t^2 + \cdots + a_n t^n
\]
and
\[
\left. c(\hat{\cal V} \right|_{\C(\n,\D) \times \{p\}})(t) = 1 + \hat{a}_1 t + \hat{a}_2 t^2 + \cdots + \hat{a}_{\n} t^{\n}.
\]
Also
$$ W_{k,l} = \sum_{i=1}^n f_i \frac{\partial \delta_k}{\partial a_i} - 
 \sum_{j=1}^{\n} \hat{f}_j \frac{\partial \hat{\delta}_k}
{\partial \hat{a}_j} +$$
$$\sum_{s=1}^g \left(\sum_{i=1}^n b_i^s \frac{\partial}{\partial a_i} + 
\sum_{j=1}^{\n} \hat{b}_j^s \frac{\partial}{\partial \hat{a}_j} \right)
\left(\sum_{i=1}^n b_i^{s+g} \frac{\partial}{\partial a_i} + 
\sum_{j=1}^{\n} \hat{b}_j^{s+g} \frac{\partial}{\partial \hat{a}_j} \right)(\delta_k - \hat{\delta}_l),
$$
and
\[
\Xi_{k,l} = \sum_{s=1}^g \left(\sum_{i=1}^n b_i^s \frac{\partial \delta_k}{\partial a_i} - 
\sum_{j=1}^{\n} \hat{b}_j^s \frac{\partial \hat{\delta_l}}{\partial \hat{a}_j} \right)
\left(\sum_{i=1}^n b_i^{s+g} \frac{\partial \delta_k}{\partial a_i} - 
\sum_{j=1}^{\n} \hat{b}_j^{s+g} \frac{\partial \hat{\delta_l}}{\partial \hat{a}_j} \right).
\]
Equivalently $c(-\pi_!(\hat{\cal V}^* \otimes {\cal V}))(t)$ equals
\begin{equation} 
\label{EC2}
\Omega(t)^{g-1} \exp \left\{
\int_{0}^{t} \left( \frac{d+\D}{u} - \left( \sum_{i=1}^n f_i \frac{\partial}{\partial a_i} - 
\sum_{j=1}^{\n} \hat{f}_j \frac{\partial}{\partial \hat{a}_j} + \right. \right. \right. 
\end{equation}
$$\left. \left. \left. \sum_{s=1}^g
\left(\sum_{i=1}^n b_i^s \frac{\partial}{\partial a_i} + \sum_{j=1}^{\n} \hat{b}_j^s \frac{\partial}{\partial \hat{a}_j} \right)
\left(\sum_{i=1}^n b_i^{s+g} \frac{\partial}{\partial a_i} + \sum_{j=1}^{\n} \hat{b}_j^{s+g} 
\frac{\partial}{\partial \hat{a}_j} \right) \right) \left( \frac{\log \Omega(u)}{u^2} \right) \right) \,\mbox{d}u \right\}.
$$
\end{prop}

{\em Proof:} We modify the proof of \cite[Prop. 10]{E}. We can write 
\[
\ch({\cal V}) = \sum_{k=1}^n e^{\gamma_k} \mbox { and } 
\ch(\hat{\cal V}^*) = \sum_{l=1}^{\n} e^{-\hat{\gamma}_l}
\]
where $\gamma_1,...,\gamma_n$ and $\hat{\gamma_1},...,\hat{\gamma}_{\n}$ are formal degree two classes such that the $r$th elementary symmetric polynomial of $\gamma_1,...,\gamma_n$
equals 
\[
c_r({\cal V}) = a_r \otimes 1 + \sum_{s=1}^{2g} b_r^s \otimes \alpha_s + f_r \otimes \omega
\]
for $1 \leq r \leq n$ and the $r$th elementary symmetric polynomial of $\hat{\gamma_1},...,\hat{\gamma}_{\n}$ equals
\[
c_r(\hat{\cal V}) = \hat{a}_r \otimes 1 + \sum_{s=1}^{2g} \hat{b}_r^s \otimes \alpha_s + \hat{f}_r \otimes \omega
\]
for $1 \leq r \leq \n$.
For each $k \geq 0$ there exist coefficients $\rho_{r_{1},...,r_{n}}^{(k)}$
such that
\[
(\gamma_{1})^{k} + \cdot \cdot \cdot + (\gamma_{n})^{k} = \sum
\rho_{r_{1},...,r_{n}}^{(k)} (c_{1}({\cal V}))^{r_{1}} \cdot \cdot \cdot
(c_{n}({\cal V}))^{r_{n}}
\]
where the sum is taken over all nonnegative $r_{1},...,r_{n}$ such that
$r_{1}+2r_{2}+ \cdots +nr_{n} = k$. Now
\[
(a_{1} \otimes 1 + \sum_{s=1}^{2g} b_{1}^{s} \otimes \alpha_{s} + f_{1} \otimes
\omega)^{r_{1}} \cdot \cdot \cdot (a_{n} \otimes 1 + \sum_{s=1}^{2g} b_{n}^{s}
\otimes \alpha_{s} + f_{n} \otimes \omega)^{r_{n}}
\]
equals
\[
(a_{1})^{r_{1}} \cdots (a_{n})^{r_{n}} \otimes 1 + \sum_{i=1}^{n}
\sum_{s=1}^{2g} b_{i}^{s} \dai (a_{1})^{r_{1}} \cdots
(a_{n})^{r_{n}} \otimes \alpha_{s}
\]
\[
+\sum_{i=1}^{n} f_{i} \dai (a_{1})^{r_{1}} \cdot \cdot \cdot (a_{n})^{r_{n}}
\otimes \omega + \sum_{i=1}^{n} \sum_{j=1}^{n} \sum_{s=1}^g b_i^s b_j^{s+g} \daij (a_{1})^{r_{1}}
\cdot \cdot \cdot (a_{n})^{r_{n}} \otimes \omega.
\]
Since
\[
\sum \rho_{r_{1},...,r_{n}}^{(k)} (a_{1})^{r_{1}} \cdot \cdot \cdot
(a_{n})^{r_{n}} = (\delta_{1})^{k} + \cdot \cdot \cdot + (\delta_{n})^{k}
\]
we find that $\ch({\cal V})$ equals
\[
\sum_{k=1}^{n} e^{\delta_{k}} \otimes 1 + \sum_{i=1}^{n} \sum_{s=1}^{2g}
\sum_{k=1}^{n} b_{i}^{s} \dai e^{\delta_{k}} \otimes \alpha_{s}
\]
\[
+\sum_{i=1}^{n} \sum_{k=1}^{n} f_{i} \dai e^{\delta_{k}} \otimes
\omega + \sum_{i=1}^{n} \sum_{j=1}^{n} \sum_{k=1}^{n} \sum_{s=1}^g b_i^s b_j^{s+g} \daij
e^{\delta_{k}} \otimes \omega. 
\]
We have a similar formula for $\ch({\cal V}^*)$ and so by the Grothendieck-Riemann-Roch 
Theorem (and substantial manipulation) obtain
\begin{eqnarray*}
\ch(-\pi_! (\hat{\cal V}^* \otimes {\cal V})) & = & \pi_* (\ch (\hat{\cal V}^*) \ch({\cal V}) ((g-1)\omega -1))\\
& = & \sum_{k=1}^n \sum_{l=1}^{\n} (g-1 - W_{k,l} - \Xi_{k,l}) e^{\delta_k - \hat{\delta}_l}.
\end{eqnarray*}
Since $W_{k,l}$ is a formal class of degree zero  and $\Xi_{k,l}$ is a formal class of degree two, it 
follows from Lemma 9 of \cite{E} that $c(-\pi_! (\hat{\cal V}^* \otimes {\cal V}))(t)$ equals (\ref{EC}) as required.
Further manipulation of this expression following the arguments of \cite[pp. 28-9]{E} give the second 
expression (\ref{EC2}). 

\begin{cor}
\label{C10}
The equivariant Chern polynomial $c(-\pi_! (\hat{\cal V}^* \otimes {\cal V}))(t)$ satisfies the equation
\[
(\Omega(t))^2 \frac{\mbox{d}}{\mbox{d}t} \left( c(-\pi_! (\hat{\cal V}^* \otimes {\cal V}))(t) \right) = Q(t) c(-\pi_! (\hat{\cal V}^* \otimes {\cal V}))(t)
\]
where 
\[
Q(t) = (g-1) \Omega'(t) - \sum_{k=1}^n \sum_{l=1}^{\n} \left( \frac{\Omega(t)}{1 + (\delta_k - \hat{\delta}_l)t} \right)^2
 (W_{k,l}(1 + (\delta_k - \hat{\delta}_l)t) + \Xi_{k,l} )
\]
is a polynomial of degree at most $2n\n-1$ in $t$ with coefficients in
\[
H^*_{\G(n,d)}(\C(n,d)) \otimes H^*_{\G(\n,\D)}(\C(\n,\D)).
\]
\end{cor}

{\em Proof:} This follows from Proposition \ref{r10}, together with two observations; firstly that $\Omega(t)$ is a 
polynomial in $t$ of degree $n\n$ and is divisible by $(1+(\delta_k - \hat{\delta}_l)t)$ for all $k,l$, and secondly 
that the polynomial 
\[
\sum_{k=1}^n \sum_{l=1}^{\n} \left( \frac{\Omega(t)}{1 + (\delta_k - \hat{\delta}_l)t} \right)^2
(W_{k,l}(1 + (\delta_k - \hat{\delta}_l)t) + \Xi_{k,l} )
\]
in $t$ is invariant under the actions of the permutation groups on $\{1,...,n\}$ and on $\{1,...,\n\}$, 
and hence its coefficients can be expressed as polynomial 
functions of the generators $a_r,b_r^j,f_r$ of $H^*_{\G(n,d)}(\C(n,d))$ and $\hat{a}_r, \hat{b}_j^r,\hat{f}_r$ 
of $H^*_{\G(\n,\D)}(\C(\n,\D)).$ 

\bigskip

Now we can complete the proof of Theorem \ref{thm}. For 
\[
c(-\pi_! (\hat{\cal V}^* \otimes {\cal V}))(t) = \sum_{r \geq 0}
c_r(-\pi_! (\hat{\cal V}^* \otimes {\cal V})) t^r
\]
and 
\[
\frac{\mbox{d}}{\mbox{d}t} \left( c(-\pi_! (\hat{\cal V}^* \otimes {\cal V}))(t) \right) = \sum_{r \geq 0}
(r+1)c_{r+1}(-\pi_! (\hat{\cal V}^* \otimes {\cal V})) t^r.
\] 
Therefore if $r \geq 0$ the coefficient of $t^r$ in 
\[
\left( \Omega(t) \right)^2  \frac{\mbox{d}}{\mbox{d}t} \left( c(-\pi_! (\hat{\cal V}^* \otimes {\cal V}))(t) \right) 
\]
is the sum of $(r+1)c_{r+1}(-\pi_! (\hat{\cal V}^* \otimes {\cal V}))$ and an element of the ideal in 
\[
H^*_{\G(n,d)}(\C(n,d)) \otimes H^*_{\G(\n,\D)}(\C(\n,\D)).
\]
generated by
\[
\{c_{r-j}(-\pi_! (\hat{\cal V}^* \otimes {\cal V})): 0 \leq j \leq 2 n \n -1\},
\]
whereas the coefficient of $t^r$ in
\[
Q(t) c(-\pi_! (\hat{\cal V}^* \otimes {\cal V}))(t)
\]
is an element of this same ideal. Thus it follows from Corollary \ref{C10} that if 
$r \geq 0$ then $c_{r+1}(-\pi_! (\hat{\cal V}^* \otimes {\cal V}))$ lies in the ideal generated by 
\[
\{c_{r-j}(-\pi_! (\hat{\cal V}^* \otimes {\cal V})): 0 \leq j \leq 2 n \n -1\}.
\]
Hence the ideal in $H^*_{\G(n,d)}(\C(n,d)) \otimes H^*_{\G(\n,\D)}(\C(\n,\D))$ generated by 
\[
\{c_{r}(-\pi_! (\hat{\cal V}^* \otimes {\cal V})): r > n \n (g-1) - d \n + \D n\}
\]
is also generated by 
\[
\{c_{r}(-\pi_! (\hat{\cal V}^* \otimes {\cal V})): n \n (g-1) - d \n + \D n < r <  n \n (g+1) - d \n + \D n\}
\]
and the ideal in $H^*_{\G(n,d)}(\C(n,d))$ generated by the slant products
\[
\{c_{r}(-\pi_! (\hat{\cal V}^* \otimes {\cal V})) \backslash \gamma: n \n (g-1) - d \n + \D n < r, 
\gamma \in H_*^{\G(\n,\D)}(\C(\n,\D)^{ss})\}
\]
is also generated by
those $c_{r}(-\pi_! (\hat{\cal V}^* \otimes {\cal V})) \backslash \gamma$ with
$$ n \n (g-1) - d \n + \D n < r < 
 n \n (g+1) - d \n + \D n, \gamma \in H_*^{\G(\n,\D)}(\C(\n,\D)^{ss})\}.
$$
This completes the strengthening of Theorem \ref{a8} to give us Theorem \ref{thm}.

\begin{rem} \label{R9.5}
Similar arguments in the parabolic case allow us to refine Theorem \ref{wPar2.1}
to get Theorem \ref{Par2.1}.
\end{rem}

\begin{rem} \label{R9.6}
In principle Corollary 10.4 gives us a recurrence relation for the relations
$c_{r}(-\pi_! (\hat{\cal V}^* \otimes {\cal V})) \backslash \gamma$
(cf. \cite{KN,Z} for the case when $n=2$).
\end{rem}

\section{Explicit formulas for relations}

\renorm

We can give explicit formulas for the relations 
$$c_r(-\pi_!(\hat{{\cal V}}^* \otimes {\cal V}))\backslash
\gamma$$
which appear in Theorem \ref{thm}
in terms of the generators $a_r, b_r^j$ and $f_r$ of $H^*_{\G(n,d)}(\C(n,d))$.
First let us consider the case when
$\n$ and $\D$ are coprime. In this case the Lefschetz duality map (\ref{ld})
is an isomorphism
$$LD: H_*^{\GG(\n,\D)}(\C(\n,\D)^{ss}) \cong H_*({\cal M}(\n,\D)) \cong 
H^{D(\n,\D)-*}({\cal M}(\n,\D))$$
$$ \cong H^{D (\n,\D) -*}_{\GG(\n,\D)}(\C(\n,\D)^{ss})$$
given by Poincar\'{e} duality on ${\cal M}(\n,\D)$, and if
$\gamma \in  H^{D (\n,\D) -*}_{\GG(\n,\D)}(\C(\n,\D)^{ss})$ and
$\eta = LD^{-1}(\gamma) \in  H_*^{\GG(\n,\D)}(\C(\n,\D)^{ss}) $ then
$$c_r(-\pi_!(\hat{{\cal V}}^* \otimes {\cal V}))\backslash
\gamma = \int_{\M} \Phi (\eta c_r(-\pi_!(\hat{{\cal V}}^* \otimes {\cal V})))
$$
where $\Phi:H^*_{\G(n,d)}(\C(n,d)) \otimes H^*_{\G(\n,\D)}(\C(\n,\D)) \to 
H^*_{\G(n,d)}(\C(n,d)) \otimes H^*(\M)$
is the natural map defined as at (\ref{14}). This integral over the moduli 
space $\M$ can be calculated using
Witten's principle of nonabelian localisation \cite{WI} and the results of \cite{JK3,JK2}, as
we shall now describe. 

In fact the moduli spaces considered in the papers \cite{JK3,JK2,WI} are
not the moduli spaces $\mnd$ but instead are the moduli spaces
${\cal M}_\Lambda (n,d)$ of semistable bundles over $\Sigma$ with coprime
rank $n$ and degree $d$ and fixed determinant line bundle $\Lambda$. However
if $Jac_d(\Sigma) = {\cal M}(1,d)$ denotes the Jacobian of degree $d$ line bundles over
$\Sigma$, then the fibration $\mbox{det}: \mnd \to Jac_d(\Sigma)$ with fibre
${\cal M}_\Lambda (n,d)$ induces an isomorphism
$$H^*(\mnd) \cong H^*({\cal M}_\Lambda (n,d)) \otimes H^*(Jac_d(\Sigma))$$
such that the generators $b_1^1,...,b_1^{2g}$ of $H^*(\mnd)$ correspond to the standard
generators of $H^*(Jac_d(\Sigma)) \cong H^*((S^1)^{2g})$ (see \cite{AB}), and thus it is easy to 
translate the results of \cite{JK3,JK2,WI} from $\mndl$ to $\mnd$.

In \cite{WI} Witten used physical methods to obtain explicit formulas for the
evaluations on the fundamental class $[\mndl]$ (or equivalently the integrals over
$\mndl$) of polynomial expressions in the generators $a_r, b_r^j$ and $f_r$. He
derived these formulas using an infinite-dimensional version of his principle of
nonabelian localisation. Using another version of nonabelian localisation \cite{JK3}
applied to $\mndl$ regarded as the symplectic quotient of an $SU(n)$ action on
an extended moduli space \cite{J}, equivalent formulas were obtained in \cite{JK2}.
This version of nonabelian localisation involves looking at the fixed points of the
action of the maximal torus $T$ of $SU(n)$ on the extended moduli space, and these
correspond, roughly speaking, to bundles of rank $n$ and degree $d$ which are
direct sums of line bundles. Technical difficulties arise because the extended
moduli space is noncompact, which is reflected in the fact that the fixed point set of $T$
has infinitely many components; this problem is overcome in \cite{JK2} by 
exploiting a certain periodicity in the situation. 
In order to calculate
$$ \int_{\M} \Phi (\eta c(-\pi_!(\hat{{\cal V}}^* \otimes {\cal V}))(t))
$$
(as a polynomial in $t$ with coefficients in $H^*_{\G(n,d)}(\C(n,d))$), we need to
consider the image of $\eta c(-\pi_!(\hat{{\cal V}}^* \otimes {\cal V}))(t)$ under the
restriction map
\begin{equation} \label{rs}
H^*_{\G(\n,\D)}(\C(\n,\D)) \to \oplus_{l=1}^{\n} H^*_{\G(1,\D_l)}(\C(1,\D_l))
\end{equation}
given by identifying our fixed $C^{\infty}$ bundle of rank $\n$ and degree $\D$ with
a direct sum of $C^{\infty}$ line bundles of degrees $\D_1$,...,$\D_{\n}$ where
$$\D_1 + ... + \D_{\n} = \D,$$
and then embedding $\C(1,\D_1) \times ... \times \C(1, \D_{\n})$ into $\C(\n,\D)$ via
$$(L_1,...,L_{\n}) \mapsto L_1 \oplus ... \oplus L_{\n}.$$
We also need to see what happens to this image when the degrees $\D_1$,...,$\D_{\n}$
of these line bundles are modified so that for some $p \in \{1,...,\n\}$ the degree
$\D_p$ is decreased by one, while $\D_{p+1}$ is increased by one and the other degrees
$\D_l$ for $l \neq p,p+1$ are all unchanged.

Let $\hat{a}_1^l$ and $\hat{b}_1^{s,l}$ (for $1 \leq s \leq 2g$) be the generators of 
$ H^*_{\G(1,\D_l)}(\C(1,\D_l))$ defined as at (\ref{genabf}). The map (\ref{rs}) is determined 
by the fact that the Chern polynomial
$$c(\hat{{\cal V}})(t) = 1 + \sum_{r=1}^{\n} (\hat{a}_r \otimes 1 + \sum_{s=1}^{2g} \hat{b}^s_r 
\otimes \alpha_s + \hat{f}_r\otimes \omega)t^r$$
maps to 
$$\prod_{l=1}^{\n} c(\hat{{\cal V}}_l)(t) = \prod_{l=1}^{\n}
(1 + (\hat{a}_1^l \otimes 1 + \sum_{s=1}^{2g} \hat{b}^{s,l}_1 
\otimes \alpha_s + \hat{d}_l\otimes \omega)t)$$
where $\hat{{\cal V}}_l$ is a universal bundle over $\C(1,\D_l) \times \Sigma$ for
$1 \leq l \leq \n$. Thus $\hat{a}_r$ is mapped to the $r$th elementary symmetric 
polynomial $\sigma_r^{\n}(\hat{a}_1^i: 1 \leq i \leq \n)$ in $\hat{a}_1^1,...,\hat{a}_1^{\n}$
(or equivalently the formal degree two classes $\hat{\delta}_1,...,\hat{\delta}_{\n}$ 
introduced in Proposition \ref{r10} are mapped to $\hat{a}_1^1,...,\hat{a}_1^{\n}$), while
$$\hat{b}_r^s \mapsto \sum_{l=1}^{\n} \hat{b}_1^{s,l} \sigma_{r-1}^{\n -1}(\hat{a}^i_1:i \neq l)$$
and
$$\hat{f}_r \mapsto  \sum_{l=1}^{\n} \left(\D_{l} \sigma_{r-1}^{\n -1}(\hat{a}^i_1:i \neq l)
+ \sum_{k\neq l} \sum_{s=1}^{2g} \hat{b}_1^{s,l} \hat{b}_1^{s+g,k} \sigma_{r-2}^{\n-2}(\hat{a}_1^i: i \neq k,l)\right).$$
The image of $c(-\pi_!(\hat{{\cal V}}^* \otimes {\cal V}))(t)$ under the map (\ref{rs}) is
\begin{equation} \label{cpwhole}  c(-\pi_!(\bigoplus_{l=1}^{\n} \hat{{\cal V}}_l^* \otimes {\cal V}))(t)
= \prod_{l=1}^{\n} c(-\pi_!(\hat{{\cal V}}_l^* \otimes {\cal V}))(t),\end{equation}
and  by Proposition \ref{r10} this is
\begin{equation}
\label{cpi}
\prod_{k=1}^{n} \prod_{l=1}^{\n} (1+(\delta_k-\hat{a}_1^l)t)^{g-1+\hat{d}_l-W_{k}}
\exp \left\{ \frac{-\Xi^{(k,l)}t}{1+(\delta_k-\hat{a}_1^l)t} \right\}
\end{equation}
where 
\begin{equation} \label{wk}
 W_{k} = \sum_{i=1}^n f_i \frac{\partial \delta_k}{\partial a_i} + 
\sum_{s=1}^g \sum_{i=1}^n b_i^s b_j^{s+g} \frac{\partial^2 \delta_k}{\partial a_i
\partial a_j} \end{equation}
and
\begin{equation} \label{xikl}
\Xi^{(k,l)} = \sum_{s=1}^g \left(\sum_{i=1}^n b_i^s \frac{\partial \delta_k}{\partial a_i} - 
 \hat{b}_1^{s,l}  \right)
\left(\sum_{i=1}^n b_i^{s+g} \frac{\partial \delta_k}{\partial a_i} - 
 \hat{b}_1^{s+g,l}  \right).
\end{equation}
If we formally modify the degrees $\D_1,...,\D_{\n}$ so that for some $p \in \{1,...,\n\}$ the degree
$\D_p$ is decreased by one, while $\D_{p+1}$ is increased by one and the other degrees
$\D_l$ for $l \neq p,p+1$ are all unchanged, then $c(-\pi_!(\hat{{\cal V}}^* \otimes {\cal V}))(t)$ 
is multiplied by the factor
\begin{equation} \label{change} \prod_{k=1}^n \frac{1 + ( \delta_k - \hat{a}_1^{p+1})t}{1 +
(\delta_k - \hat{a}^p_1)t}.  \end{equation}

In \cite{JK2} Theorem 8.1 (see also \cite{JK2} Remark 8.3(c)) it is proved that if
$\eta$ is any polynomial expression in the generators $\hat{a}_r$ and $b_r^j$ of
$H^*_{\GG(\n,\D)}(\C(\n,\D))$, not involving the generators $f_r$, then
\begin{equation} \label{residue}
\int_{\M} \Phi( \eta \mbox{ exp}(\epsilon \hat{f}_2)) = \frac{(-1)^{\n(\n -1)(g-1)/2}}{\n !}
\mbox{Res}_{Y_1 =0} ...\end{equation}
$$...\mbox{Res}_{Y_{\n-1} =0} \left( \sum_{w \in W_{\n -1}} \frac{\exp(\langle \epsilon
[[w \hat{c}]],X \rangle) \int_{(S^1)^{2g{\n}}} \eta(X) \exp(\epsilon \omega)}{{\cal D}_{\n}^{2g-2}
\prod_{l=1}^{{\n} -1} (\exp(\epsilon Y_l) - 1)} \right) $$
when $\epsilon \neq 0$, where
$(S^1)^{\n}$ is identified in the usual way with the maximal torus of $U(\n)$ consisting
of diagonal unitary matrices (which appears in (\ref{residue}) but not
in \cite{JK2}, because we are working with the moduli space $\M$ not ${\cal M}_\Lambda (\n,\D)$), 
while $T=(S^1)^{\n} \cap SU(\n)$ is the standard maximal torus of 
$SU(\n)$. Also any $X \in \mbox{Lie}(T)$ has coordinates
$$Y_1 = X_1 - X_2,...,Y_{{\n} -1} = X_{{\n} -1} - X_{\n}$$
defined by the simple roots of $SU(\n)$, and the Weyl group $W_{\n -1}$ of 
$SU(\n -1)$ is embedded in $SU(\n)$ in the standard way using the first $\n - 1$ coordinates.
The polynomial function ${\cal D}_{\n} : \mbox{Lie}(T) \to \RR$ is defined to be the product
of the positive weights of $SU(\n)$. Also if $\gamma \in \mbox{Lie}(T)$ then
$[[\gamma]]$ is the unique element of the fundamental domain defined by the
simple roots for the translation action on $\mbox{Lie}(T)$ of the integer lattice, such
that $[[\gamma]] - \gamma$ lies in the integer lattice, and
$$\hat{c} = [[(\D/\n,...,\D/\n, \D/\n - \D)]].$$
Moreover $\mbox{Res}_{Y_j = 0}$ is the usual residue at 0 of a meromorphic
function of the variable $Y_j$, and $\langle  ,  \rangle $ denotes the inner product
on $\mbox{Lie}(SU(\n))$ given by
$$\langle X,X \rangle = - \mbox{Trace}(X^2)/4 \pi^2 .$$
$\omega$ is the standard symplectic form on $(S^1)^{2g\n}$, normalised so that
$$\int_{(S^1)^{2g\n}} \exp(\omega) = (\n)^g,$$
and the appearance of $\eta$ in the integral
$$ \int_{(S^1)^{2g\n}} \eta(X) \exp(\epsilon \omega)$$
is interpreted as the image of $\eta$ under the map (\ref{rs}), via the
identification
\begin{equation} \label{iden} H^*_{\G(1,\D_l)}(\C(1,\D_l)) \cong
H^*(BS^1) \otimes H^*({\cal M}(1,\D_l))
\end{equation}
when we identify $(S^1)^{2g}$ with the Jacobian $Jac_{\D_l}(\Sigma)={\cal M}(1,\D_l)$.
It does not matter here which $\D_l$ (for $1 \leq l \leq \n$) we choose in the identification of
$(S^1)^{2g}$ with ${\cal M}(1,\D_l)$; tensoring by any fixed line bundle of degree $\delta$ gives
us an isomorphism from ${\cal M}((1,\D_l)$ to ${\cal M}(1, \D_l + \delta)$ which preserves the image
of $\eta$ since $\eta$ is a polynomial in the generators $\hat{a}_r$ and $\hat{b}_r^j$ (the
generators $f_r$ are not involved). However, in some sense the appropriate choice for $\D_l$ here
is $\D/\n - [[w \hat{c}]]_l$; this accounts for the occurrence of the term 
$\exp(\langle \epsilon
[[w \hat{c}]],X \rangle)$ in the formula (\ref{residue}).

\begin{rem} \label{szenes} 
The formula given by Witten in \cite{WI} for an integral of the form
$$
\int_{\M} \Phi( \eta \mbox{ exp}(\epsilon \hat{f}_2)) $$
is not expressed as an iterated residue. Instead it is an infinite sum
running over the intersection of the weight lattice of $SU(\n)$ with the interior
of a fundamental Weyl chamber. These two formulas can be reconciled 
using an argument of Szenes (see \cite{JK2} Proposition 2.2), based on the
elementary fact that the sum of all the residues of a meromorphic function on
$\CC$ which vanishes sufficiently fast at infinity is zero. This means that
an iterated residue  of the form (\ref{residue}) can be rewritten as a sum
of iterated residues at the nonzero points $Y_l \in 2 \pi i \epsilon^{-1} \ZZ \backslash \{ 0 \}$
where the terms
$ \exp (\epsilon Y_l) - 1$
vanish. There are of course infinitely many such points, but they contribute
simple poles whose residues are easy to calculate.
\end{rem}

The   product
$$\prod_{l=1}^{\n -1} ( \exp (\epsilon Y_l) - 1)$$
appears in the residue formula (\ref{residue}) because when we formally
modify  the degrees $\D_1,...,\D_{\n}$ so that for some $p \in \{1,...,\n\}$ the degree
$\D_p$ is decreased by one, while $\D_{p+1}$ is increased by one and the other degrees
$\D_l$ for $l \neq p,p+1$ are all unchanged, then the image of $\eta \exp(\epsilon \hat{f}_2)$
under the map (\ref{rs}) is multiplied by $\exp(\epsilon (\hat{a}^p_1 - \hat{a}_1^{p+1}))$.
The generators $\hat{a}_1^1,...,\hat{a}_1^{\n}$ correspond to the coordinates
$X_1,...,X_{\n}$ on the Lie algebra of the maximal torus of $U(\n)$, and 
$\exp(\epsilon (\hat{a}^p_1 - \hat{a}_1^{p+1}))$ then corresponds to
$$\exp(\epsilon (X_p - X_{p+1})) = \exp (\epsilon Y_p).$$
If we apply the proof of \cite{JK2} Theorem 8.1 to $\eta c(-\pi_!(\hat{{\cal V}}^* \otimes {\cal V}))(t)$ 
instead of $\eta \exp ( \epsilon \hat{f}_2)$, then by (\ref{change}) the term
$$\prod_{l=1}^{\n -1} ( \exp (\epsilon Y_l) - 1)$$
is replaced by
$$\prod_{l=1}^{\n -1} \left( \prod_{k=1}^n \frac{1+ (\delta_k - X_{l+1})t}{1+( \delta_k - X_l)t} -1\right),$$
and the same proof gives us

\begin{thm} \label{explicit} If
$\eta$ is any polynomial expression in the generators $\hat{a}_r$ and $b_r^j$ of
$H^*_{\GG(\n,\D)}(\C(\n,\D))$, not involving the generators $f_r$, then
$$ \int_{\M} \Phi( \eta c(-\pi_!(\hat{{\cal V}}^* \otimes {\cal V}))(t)
) = \frac{(-1)^{\n(\n -1)(g-1)/2}}{\n !}
\mbox{Res}_{Y_1 =0} ... $$
$$...\mbox{Res}_{Y_{\n-1} =0} \left( \sum_{w \in W_{\n -1}} 
\int_{(S^1)^{2g{\n}}} \eta(X) \prod_{l=1}^{\n} \prod_{k=1}^n (1 + (\delta_k - X_l)t)^{g-1 + \D/\n -[[w\hat{c}]]_l - W_k}\right.$$
$$ \left. \frac{ \exp\left( \frac{-\Xi^{(k,l)}}{1 + (\delta_k - X_l)t} \right)}{{\cal D}_{\n}^{2g-2}
\prod_{l=1}^{{\n} -1} ( \prod_{k=1}^n \frac{1+ (\delta_k - X_{l+1})t}{1+( \delta_k - X_l)t} 
 - 1)} \right) . $$
\end{thm}

\begin{rem} \label{bim}
Here we interpret 
$$( \prod_{k=1}^n \frac{1+ (\delta_k - X_{l+1})t}{1+( \delta_k - X_l)t} 
 - 1)^{-1}$$
as
$$\left( \prod_{k=1}^n (1+ (\delta_k - X_{l+1})t) \right) (nt Y_l)^{-1} \sum_{m \geq 0} (-1)^m \left(\frac{ b_l(t,X_1,...,X_{\n})}{n} \right)^m$$
where
$$ b_l(t,X_1,...,X_{\n}) = \sum_{i=0}^{n-1} ((1-tX_{l+1})^i (1-tX_l )^{n-i-1} - 1)$$
$$ + \sum_{j=1}^{n-1} t^j a_j \sum_{i=0}^{n-j-1} ((1-tX_{l+1})^i (1-tX_l )^{n-j-i-1} - 1)$$
using the expansion of
$$\prod_{k=1}^n (1 + (\delta_k - X_{l+1})t) - \prod_{k=1}^n (1 + (\delta_k - X_{l})t) $$
as 
$$t^n \left( \sum_{j=0}^n a_j (( t^{-1} - X_{l+1})^{n-j} - (t^{-1}-X_l)^{n-j}) \right) $$
$$ =
\sum_{j=0}^{n-1} a_j t^j \sum_{i=0}^{n-j-1} (1-tX_{l+1})^i (1-tX_l)^{n-j-i-1} t (X_l - X_{l+1})$$
$$= t(X_l - X_{l+1})(n + b_l(t,X_1,...,X_{\n})).$$
Notice also that $\D/\n - [[w \hat{c}]]_l$ is an integer for each $w \in W_{\n -1}$ and $1 \leq l \leq \n$.
\end{rem}

More generally if we want to allow $\eta$ to depend on the generators $\hat{f}_r$ as well as the
generators $\hat{a}_r$ and $\hat{b}_r^j$, then we can modify the proof of \cite{JK2} Theorem 9.12 to
obtain for any scalars $\epsilon_2,...,\epsilon_{\n}$ with $\epsilon_2 \neq 0$ the following result,
which will give us the formulas we require if we differentiate with respect to $\epsilon_2,...,\epsilon_{\n}$.

\begin{thm} \label{explicit2} Let $m_r$ and $p_{r,k_r}$ be nonnegative integers for
$2 \leq r \leq \n$ and $1 \leq k_r \leq 2g$, with $p_{r,k_r} \in \{0,1\}$. If $\epsilon_2 \neq 0$ then the integral
$$ \int_{\M} \Phi \left( \exp(\epsilon_2 \hat{f}_2 + ... + \epsilon_{\n} \hat{f}_{\n})\prod_{r=1}^{\n} \hat{a}_r^{m_r}
\prod_{k_r =1}^{2g} (\hat{b}_r^{k_r})^{p_{r,k_r}}
c(-\pi_!(\hat{{\cal V}}^* \otimes {\cal V}))(t)
\right) $$
is given by
$$ \frac{(-1)^{\n(\n -1)(g-1)/2}}{\n !}
\mbox{Res}_{Y_1 =0} ... 
...\mbox{Res}_{Y_{\n-1} =0} \left( \sum_{w \in W_{\n -1}} e^{-dq_X([[w \hat{c}]])} 
\prod_{r=2}^{\n} \sigma_r^{\n}(X)^{m_r} \right. $$
$$  \int_{(S^1)^{2g{\n}}} 
\prod_{k_r =1}^{2g} (\sum_{a=1}^{\n -1} (d\sigma_r^{\n})_X(\hat{e}_a) \zeta_a^{k_r})^{p_{r,k_r}}
 \prod_{l=1}^{\n} \prod_{k=1}^n (1 + (\delta_k - X_l)t)^{g-1 + \D/\n -[[w\hat{c}]]_l - W_k}$$
$$\left. \frac{ \exp\left( \frac{-\Xi^{(k,l)}}{1 + (\delta_k - X_l)t} - \sum_{a,b=1}^{\n -1} \sum_{s=1}^{2g}
\zeta_a^l \zeta_b^{l+g} \partial^2 q_X(\hat{e}_a,\hat{e}_b) \right)}{{\cal D}_{\n}^{2g-2}
\prod_{l=1}^{{\n} -1} ( \exp(-(dq)_X(\hat{e}_l))\prod_{k=1}^n \frac{1+ (\delta_k - X_{l+1})t}{1+( \delta_k - X_l)t} 
 - 1)} \right)  $$
where as before $\sigma_r^n$ denotes the $r$th elementary symmetric polynomial in $n$ variables,
and
$$q(X) = \epsilon_2 \sigma_2^{\n}(X) + ... + \epsilon_{\n} \sigma_{\n}^{\n}(X).$$
Also $\hat{e}_1, ..., \hat{e}_{\n -1}$ is the basis of $\mbox{Lie}(T)$ given by the
simple roots; if we write the Maurer-Cartan form $\theta$ on $T$ as
$\theta = \sum_{a=1}^{\n -1} \theta_a \hat{e}_a$ then $\theta_1,...,\theta_{\n-1}$ form
a set of generators for $H^1(T)$ and we let $\zeta^{l}_a$ be the pullback of
$\theta_a$ to $H^1(T^{2g})$ under the projection from $T^{2g}$ onto the $l$th 
copy of $T$.
\end{thm}

These theorems \ref{explicit} and \ref{explicit2} give us explicit formulas for the relations
$c_r (-\pi_!(\hat{{\cal V}}^* \otimes {\cal V})) \backslash \gamma$ which appear in
Theorem 2.1, in the case when $\n$ and $\D$ are coprime.

\begin{rem} \label{martin}

\newcommand{\liek}{{\frak k}}
\newcommand{\lieks}{{\frak k}^*}

\newcommand{\liet}{{\frak t}}
\newcommand{\liets}{{\frak t}^*}

Moduli spaces of parabolic bundles with fixed determinant can be 
regarded as symplectic quotients of the extended 
moduli spaces of \cite{J}, where the action is that of the subgroup $P \cap SU(n)$ of $SU(n)$.
In particular in the case when $m=n$ and $j_1 = \ldots = j_m =1$, which was studied in $\S$7 and $\S$8, we
have the symplectic quotient of the extended moduli space by the action of the maximal torus
$T$ of $SU(n)$. A crucial step in the application of nonabelian localisation in \cite{JK2} to
obtain formulas for intersection pairings on the moduli spaces ${\cal M}_\Lambda (\n,\D)$
of vector bundles with fixed determinant was an observation due to S. Martin \cite{M,M1}
concerning any Hamiltonian action of a compact group $K$ with Lie algebra
$\liek$ on a compact symplectic manifold $X$
when 0 is a regular value of the moment map $\mu:X \to \lieks$. 
If we assume for simplicity that 0 is also a regular value of the induced
moment map $\mu_T: X \to \liets$ for the maximal torus $T$ of $K$,
then Martin observed that the integral over the symplectic 
quotient $\mu^{-1}(0)/K$ of the image $\eta_0$ of any $\eta \in H^*_K(X)$ in 
$H^*_K(\mu^{-1}(0)) \cong H^*(\mu^{-1}(0)/K)$ can be
expressed as a fixed constant multiple of the integral
over $\mu_T^{-1}(0)/T$ of the image in 
$H^*_T(\mu_T^{-1}(0)) \cong H^*(\mu_T^{-1}(0)/T)$ of 
the product of ${\cal D} \in H^*_T(X)$ and $\eta \in H^*_K(X) \cong [H_T^*(X)]^W \subseteq H_T^*(X)$
where ${\cal D}$ is the product of the positive roots of $K$ and $W$ is the Weyl group.

With this in mind, the proof of Theorem 2.1 via Theorem \ref{wPar2.1} shows us that the explicit
formulas for relations provided by Theorems \ref{explicit} and \ref{explicit2}, in the case
when $\n$ and $\D$ are coprime, also give relations when $\n$ and $\D$ are not coprime, and
the set of all such relations for all $\n$ and $\D$ satisfying $0<\n<n$ and $\D/\n > d/n$ gives
a complete set of relations between the generators $a_r$, $b_r^j$ and $f_r$ of
$H^*_{\G(n,d)}(\C(n,d)^{ss})$.
\end{rem}

\end{document}